\documentstyle[preprint,aps,prb]{revtex}
\def\C{{C\kern-.647em C}}
\def\MC{M_{\C}}
\def\MAC{M_{{\cal{A}},\C}}
\def\MCA{M_{\C,\cal{A}}}
\def\MVC{M_{V,\C}}
\def\MCV{M_{\C,V}}
\def\MCVbar{M_{\C,\bar{V}}}
\def\sqtimes{\Box{\raise 0.1em\hbox{\kern-0.76em $\times$}}}
\def\twsqtimes{\widetilde{\raise-0.1em\hbox{$\sqtimes$}}}

\begin{document}

\title{Irreducible tensor operators in the\\ regular coaction formalisms of
compact\\ quantum group algebras}
\author{J.F.Cornwell} 
\address{School of Physics and Astronomy, University of 
St.Andrews,\\North Haugh, St.Andrews, Fife,
KY16 9SS, Scotland, U.K.} 
\date{23 May 1996} 
\maketitle
\begin{abstract} The defining conditions for the irreducible tensor 
operators associated with the
unitary irreducible corepresentions of compact quantum group algebras 
are deduced first in both the
right and left regular coaction formalisms. In each case it is shown that 
there are {\em{two}} types
of irreducible tensor operator, which may be called `ordinary' and 
`twisted'. The consistency of the
definitions is demonstrated, and various consequences are deduced, 
including generalizations of the
Wigner-Eckart theorem for both the ordinary and twisted operators. Also 
included are discussions
(within the regular coaction formalisms for compact quantum group 
algebras) of inner-products,
basis functions, projection operators, Clebsch-Gordan coefficients, and two 
types of tensor product
of corepresentations.

The formulation of quantum homogeneous spaces  for compact quantum 
group algebras is discussed, and
the defining conditions for the irreducible tensor operators associated with 
such quantum
homogeneous spaces and with the unitary irreducible corepresentions of 
the compact quantum group
algebras are then deduced. There are two versions, which correspond to 
restrictions of the right and
left regular coactions. In each case it is again shown that there are 
ordinary and twisted
irreducible tensor operators. Various consequences are deduced, including 
the corresponding
generalizations of the Wigner-Eckart theorem. 

\end{abstract}
\pacs{02.10.T, 02.20}

 \centerline{{\sl Running title:} Irreducible tensor operators}
\newpage

\section{Introduction}

It is well known that most of the applications to physics of the theories of 
groups and Lie algebras
depend on the Wigner-Eckart theorem. It is therefore not surprising that 
the question of the
generalization   of this theorem to Hopf algebras having the structure of a 
deformation of a Lie
algebra has also been the subject of a number of studies$^{1-15}$. The 
present paper is intended to
complement and extend these investigations in various important 
respects. Its detailed relationship
to previous work will be indicated in the appropriate places.

The perspective of the present communication is best introduced by 
considering matters first in the
very well established and  familiar context of a compact Lie group 
$\cal{G}$ (c.f. Refs.16,17). Even
in this context, one can distinguish {\em{three}} distinct forms of the 
Wigner-Eckart theorem: 

\begin{enumerate}

\item The original form$^{18}$ involves the situation in which $\cal{G}$ 
is a group of
tranformations that act on an external manifold $\cal{M}$, the classic 
example being the case in
which $\cal{M}$ is three-dimensional Euclidean space $\Re^{3}$, and 
$\cal{G}$ is the group of all
rotations in this space about some fixed point, which may be taken to be 
the origin $O$ of
$\Re^{3}$. Associated with every such rotation
$T$ there exists a $3 \times 3$ real orthogonal matrix ${\bf{R}}(T)$, so 
that the effect of $T$ is
to transform each position vector ${\bf{r}}$ into another position vector 
${\bf{r}}^{\prime}$, where
\begin{equation} {\bf{r}}^{\prime} = {\bf{R}}(T){\bf{r}} . 
\label{eq:1.2}\end{equation}
Also
associated with every rotation $T$ is a unitary operator $P(T)$ whose 
effect on any function
$f({\bf{r}})$ is defined by 
\begin{equation} P(T)\,f({\bf{r}}) = f({\bf{R}}(T)^{-1}{\bf{r}}). 
\label{eq:1.17} \end{equation} 
Let  ${\bf{\Gamma}}^{p}$ be a unitary irreducible representation of
dimension
$d_{p}$ of the group $\cal{G}$. If there exists a set of functions 
$\psi^{p}_{1}({\bf{r}}),
\psi^{p}_{2}({\bf{r}}),\ldots, \psi^{p}_{d_{p}}({\bf{r}})$ such that 
\begin{equation} P(T)\;\psi^{p}_{n}({\bf{r}}) = \sum_{m=1}^{d_{p}}
\Gamma^{p}(T)_{mn}\psi^{p}_{m}({\bf{r}}) \label{eq:1.26}\end{equation} 
for all $T \in \cal{G}$ and
all $n = 1, 2,
\ldots, d_{p}$, then these are said to form a set of basis functions for  
${\bf{\Gamma}}^{p}$. Similarly, if there exists a set of $d_{p}$ operators 
$Q^{p}_{1}, Q^{p}_{2},
\ldots, Q^{p}_{d_{p}}$ that act on functions $f({\bf{r}})$  in such a way 
that 
\begin{equation} P(T)\;Q^{p}_{n}\;P(T)^{-1} = \sum_{m=1}^{d_{p}}
\Gamma^{p}(T)_{mn}Q^{p}_{m} \label{eq:5.34} \end{equation}     for all 
$T \in \cal{G}$ and all $n =
1, 2,
\ldots, d_{p}$, then these are said to form a set of irreducible tensor 
operators for  
${\bf{\Gamma}}^{p}$. Finally, if the inner product for the Hilbert space of 
functions
$f({\bf{r}})$ is defined by
\begin{equation} (f,g) = \int_{-\infty}^{\infty}\int_{-
\infty}^{\infty}\int_{-\infty}^{\infty}
\overline{f({\bf{r}})}g({\bf{r}}) 
\;\mbox{d}x\;\mbox{d}y\;\mbox{d}z\;,\label{eq:1.19}
\end{equation}
where ${\overline{f({\bf{r}})}}$ denotes the complex conjugate of
$f({\bf{r}})$,  then the Wigner-Eckart theorem for this situation states 
that the $j,
k,$ and
$\ell$ dependence of 
$(\psi^{r}_{\ell},Q^{q}_{k}(\phi^{p}_{j}))$ depends only on Clebsch-Gordan 
coefficients for the
reduction of the tensor product ${\bf{\Gamma}}^{p} \otimes 
{\bf{\Gamma}}^{q}$ into its irreducible
constituents ${\bf{\Gamma}}^{r}$.

\item In this form the role of the manifold $\cal{M}$ is played by 
$\cal{G}$ itself, so that one is
concerned with the space of complex-valued continuous functions defined 
on $\cal{G}$. Let this be denoted by
$C({\cal{G}})$. The inner product of $C({{\cal{G}}})$ may be taken to be
\begin{equation} (f,g) = \int_{{\cal{G}}} \overline{f(T)}g({T})
\;\mbox{d}T\;,\label{eq:Haar} \end{equation}   where the integral is the 
left and right invariant
normalised Haar integral of $\cal{G}$, and $\overline{f(T)}$ is the 
complex conjugate of ${f(T)}$. 
In the {\em{right regular formalism}}, for each $T \in \cal{G}$ there 
exists an operator
$\widehat{R}(T)$ that is defined by
$\widehat{R}(T)\,f(T^{\prime}) = f(T^{\prime}T)$  for all $f$ and for all 
$T,T^{\prime} \in
\cal{G}$. If $f$ is a
member of $C({\cal{G}})$ such that $\widehat{R}(T)f$ spans a finite-
dimensional subspace of
$C({\cal{G}})$, then $f$ is said to be a {\em{representative function}} on 
$\cal{G}$. The subspace
of $C({\cal{G}})$ consisting of representative functions will be denoted by 
$R({\cal{G}})$. If there
exists a set of functions
$\psi^{p}_{1}(T),
\psi^{p}_{2}(T),\ldots, \psi^{p}_{d_{p}}(T)$ such that $ 
\widehat{R}(T)\,\psi^{p}_{n}(T^{\prime}) =
\sum_{m=1}^{d_{p}}
\Gamma^{p}(T)_{mn}\psi^{p}_{m}(T^{\prime})$ for all $T,T^{\prime} \in 
\cal{G}$ and all
$n = 1, 2, \ldots, d_{p}$, then these are said to form a set of basis functions 
for 
${\bf{\Gamma}}^{p}$. Similarly, if there exists a set of $d_{p}$ operators 
$Q^{p}_{1}, Q^{p}_{2},
\ldots, Q^{p}_{d_{p}}$ that act on functions $f(T)$  in such a way that 
$\widehat{R}(T)\,Q^{p}_{n}\,\widehat{R}(T)^{-1} = \sum_{m=1}^{d_{p}}
\Gamma^{p}(T)_{mn}Q^{p}_{m}$ for all $T \in \cal{G}$ and all $n = 1, 2, 
\ldots, d_{p}$, then this
set is said to form a set of irreducible tensor operators for  
${\bf{\Gamma}}^{p}$. The Wigner-Eckart theorem for this case states 
that the $j, k,$ and $\ell$
dependence of 
$(\psi^{r}_{\ell},Q^{q}_{k}(\phi^{p}_{j}))$ again depends only on Clebsch-
Gordan coefficients for the
reduction of the tensor product ${\bf{\Gamma}}^{p} \otimes 
{\bf{\Gamma}}^{q}$ into its irreducible
constituents ${\bf{\Gamma}}^{r}$. In the {\em{left regular formalism}} 
the situation is the same,
except only that the operators $\widehat{R}(T)$ are replaced by operators 
$\widehat{L}(T)$ that are
defined by
$\widehat{L}(T)\,f(T^{\prime}) = f(T^{-1}T^{\prime})$ for all $f$ and for 
all
$T,T^{\prime} \in \cal{G}$.

\item The final form involves using the abstract carrier spaces of the 
unitary irreducible
representations of $\cal{G}$. Let $V^{p}$ be such a carrier space for
${\bf{\Gamma}}^{p}$, with ortho-normal basis $\psi^{p}_{1}, 
\psi^{p}_{2},\ldots, \psi^{p}_{d_{p}}$,
and define for each $T \in \cal{G}$ a linear operator $\Phi^{p}(T)$ that 
acts on $V^{p}$ by the
requirement that $\Phi^{p}(T)\,\psi^{p}_{n} = \sum_{m=1}^{d_{p}}
\Gamma^{p}(T)_{mn}\psi^{p}_{m}$  for all $T \in \cal{G}$ and all $n = 1, 
2, \ldots, d_{p}$. Let
${\bf{\Gamma}}^{p}$,
${\bf{\Gamma}}^{q}$, and ${\bf{\Gamma}}^{r}$ be any three unitary 
irreducible resentations of
$\cal{G}$. Then one can consider a set of irreducible tensor operators 
$Q^{q}_{1}, Q^{q}_{2},
\ldots, Q^{q}_{d_{q}}$ that each map
$V^{p}$ into $V^{r}$ and which are such that
$\Phi^{r}(T)\,Q^{q}_{n}\,\Phi^{p}(T)^{-1} = \sum_{m=1}^{d_{q}}
\Gamma^{q}(T)_{mn}Q^{q}_{m}$ for all $T \in \cal{G}$ and all $n = 1, 2, 
\ldots, d_{q}$. In this case
the Wigner-Eckart theorem deals with inner products $\langle\;,\; 
\rangle$ defined on
$V^{r}$ and states that the $j, k,$ and $\ell$ dependence of 
$\langle \psi^{r}_{\ell},Q^{q}_{k}(\phi^{p}_{j}) \rangle$ also depends only 
on  Clebsch-Gordan
coefficients for the reduction of the tensor product ${\bf{\Gamma}}^{p} 
\otimes {\bf{\Gamma}}^{q}$
into its irreducible constituents ${\bf{\Gamma}}^{r}$. In a minor 
extension of this formalism, one
could introduce an inner product space $V$ that is a direct sum of carrier 
spaces of certain unitary
irreducible representations of $\cal{G}$ and which contains at least 
$V^{p}\oplus V^{r}$ (and which,
in the extreme case,  may contain one carrier space for every inequivalent 
irreducible
representation of $\cal{G}$).  Then, for each $T \in \cal{G}$ an operator 
$\Phi(T)$ can be defined
which maps elements of $V$ into $V$, and which acts as
$\Phi^{p}(T)$ on $V^{p}$, as $\Phi^{r}(T)$ on $V^{r}$, and so on.  The 
irreducible tensor operators
are then required to each map
$V$ into $V$ and to be such that
$\Phi(T)\,Q^{q}_{n}\;\Phi(T)^{-1} = \sum_{m=1}^{d_{q}}
\Gamma^{q}(T)_{mn}Q^{q}_{m}$ for all $T \in \cal{G}$ and all $n = 1, 2, 
\ldots, d_{q}$. In this case
the Wigner-Eckart theorem deals with inner products $\langle\;,\; 
\rangle$ defined on
$V$, but is otherwise the same as above. 

\end{enumerate}

The developments that will be described in the present paper up and 
including
Section~VIII are essentially within the spirit of the second of these 
formulations,
but deal with a more general Hopf algebra structure. The generalization of 
the first formulation
in terms of quantum homogeneous spaces then follows in Section~IX.  (It 
is
intended to extend this analysis to the remaining formulation in a 
subsequent paper).

One most important lesson that can be drawn from these simple group 
theoretical considerations
concerns the {\em{consistency}} of the definitions of the basis functions (or 
basis vectors)
{\em{and}} of the irreducible tensor operators. The essential point will be 
illustrated in the first
of the above formulations, but similar considerations apply in the others.  
As $P(T)P(T^{\prime}) =
P(TT^{\prime})$ and 
${\bf{\Gamma}}^{p}(T){\bf{\Gamma}}^{p}(T^{\prime}) =
{\bf{\Gamma}}^{p}(TT^{\prime})$ for all $T,T^{\prime} \in \cal{G}$, it 
follows that if
(\ref{eq:1.26}) is valid for $T$ and for $T^{\prime}$, then it is also valid 
for their product
$TT^{\prime}$. Similarly, and very significantly, by defining for each $T
\in \cal{G}$ an operator $\Psi(T)$ by $\Psi(T)Q = P(T)QP(T)^{-1}$ for 
every operator $Q$ that acts
on functions $f({\bf{r}})$, the definition (\ref{eq:5.34}) can be recast as  
\begin{equation} \Psi(T)(Q^{p}_{n}) = \sum_{m=1}^{d_{p}}
\Gamma^{p}(T)_{mn}Q^{p}_{m} \label{eq:5.34a} \end{equation}      for all 
$T \in \cal{G}$ and all $n
= 1, 2, \ldots, d_{p}$. As $\Psi(T)\Psi(T^{\prime}) =
\Psi(TT^{\prime})$  for all
$T,T^{\prime} \in \cal{G}$, it follows that if (\ref{eq:5.34a}) is valid for 
$T$ and for
$T^{\prime}$, then it is also valid for their product $TT^{\prime}$. Put 
another way, because of the
similarity in form between (\ref{eq:1.26}) and (\ref{eq:5.34a}), the 
{\em{consistency}} of the
definition (\ref{eq:5.34}) of the irreducible tensor operators $Q^{p}_{n}$ is 
ensured by the fact
that they too form a basis for a carrier space of ${\bf{\Gamma}}^{p}$. In 
the analysis that follows
(cf. Section~VI), essentially this argument will be used to justify the
definitions that will be given for the irreducible tensor operators of the 
compact quantum group
algebras in the regular corepresentation formalisms, the only essential 
difference being that the
argument has to be cast in terms of {\em{co}}representations instead of 
representations.

It is well known that the set of functions defined on a Lie group $\cal{G}$ 
form a Hopf algebra,
$\cal{A}$, and that the dual ${\cal{A}}^{\prime}$ of $\cal{A}$ is the 
universal enveloping algebra
of the Lie algebra $\cal{L}$ of $\cal{G}$. Moreover, the structure of 
$\cal{G}$ can be encoded into
the structure of $\cal{A}$, and, in particular,
$\cal{A}$ is commutative.  A `deformation' (or `quantization') of  
${\cal{A}}^{\prime}$ induces a
corresponding deformation of $\cal{A}$, and will make $\cal{A}$ non-
commutative as well as being
non-cocommutative. Most of the previous work on irreducible tensor 
operators has been focused on the
deformed Hopf algebras ${\cal{A}}^{\prime}$, with $su_{q}(2)$ receiving 
the most attention. However,
as has been demonstrated by the pioneering work of Woronowicz$^{19-
21}$, which itself has been
refined and developed by Dijkhuizen and Koornwinder$^{22-26}$, it is of 
very great interest to
produce a self-contained and direct study of generalizations of the Hopf 
algebras $\cal{A}$, which
can be done by assuming that they have certain characteristic properties. 
The resulting structures
have been  called {\em{compact matrix pseudogroups}} by 
Woronowicz$^{19-21}$, and {\em{compact
quantum group algebras}} by Dijkhuizen and Koornwinder$^{22-26}$. 
These provide the framework for
the  present paper, which is devoted to the study of the irreducible tensor 
operators for compact
quantum group algebras. As explained above, this analysis will be given in 
the regular
corepresentation formalisms. (The only previous investigation of 
irreducible tensor operators within
the general compact matrix pseudogroup theory has been by 
Bragiel$^{5}$, who looked at the analogue
of the carrier space formalism (3) above, but with certain restrictive 
assumptions on
multiplicities, though some of the work of Klimyk$^{9}$ involves a 
discussion of special cases,
again in the carrier space formalism). 

 The structure of the present paper is as follows. Section~II contains a 
brief
summary of the essential preliminaries, starting in Subsection~II.A
 with the properties of Hopf $*$-algebras, and continuing in 
Subsection~II.B with the main features of their right
comodules.  The definition and relevant properties of a compact quantum 
group algebra $\cal{A}$
follow in Subsection~II.C. (Of course the developments of
Woronowicz of and Dijkhuizen and Koornwinder extend far beyond what 
is mentioned here, particularly
in their  their invocation of quantum Tannaka-Krein duality.)  This 
section is concluded in
Subsection~II.D with some new lemmas concerning the
Haar functional of $\cal{A}$. The right and left regular comodules of 
$\cal{A}$ are described in
Section~III, and these are employed in Section~IV to
introduce and develop the concept of basis functions for right 
corepresentations of $\cal{A}$. In
Section~V the tensor products (both `ordinary' and `twisted') of
corepresentations of $\cal{A}$ are discussed, along with their associated 
Clebsch-Gordan
coefficients. The heart of the paper is reached in Section~VI, where the
irreducible tensor operators are defined and some of their immediate 
properties are deduced. In
particular, it will be shown there that in both the right and left regular 
coaction formulations
there are {\em{two}} types of irreducible tensor operators, which will be 
described as being
{\em{ordinary}} and {\em{twisted}} respectively. The motivations for the 
definitions of
Section~VI are deliberately relegated to Appendix B
 in order to emphasize that the treatment of given for the compact 
quantum
group algebras in Sections~II to IX is entirely self-contained.  In
Section~VII it is shown that there are {\em{two}} theorems of the Wigner-
Eckart
type, one for `ordinary' and one for the `twisted' irreducible tensor 
operators. Likewise, in
Section~VIII, it is demonstrated that these two types of irreducible tensor 
operator
behave differently under multiplication. Finally, in Section~IX it is
shown how all developments generalize when one considers operators 
associated with the
corresponding homogeneous spaces. In particular, it emerges that there 
are again two formulations,
one associated with the right regular representation and the other with 
the left regular
represntation. The vital algebraic quantity that appears in each version is 
a $\star$-subalgebra
$\cal{B}$ of $\cal{A}$, which is  a right coideal of $\cal{A}$ in the right 
regular
formulation, but is a left coideal of $\cal{A}$ in the left regular
formulation. In Subsection~IX.B attention is focused on the right coactions
$\pi^{R}_{\cal{B}}$ and $\pi^{L}_{\cal{B}}$ of  $\cal{A}$ that are 
obtained by restricting the right
and left regular coactions of $\cal{A}$ to its subalgebra $\cal{B}$. As 
these are the
{\em{transitive $\star$-coactions}} that correspond to the transitive 
action of a quantum group on a
quantum homogeneous space in the sense of Dijkhuizen and 
Koornwinder$^{23,24}$, they play the key
role in the analysis. In particular the properties of basis functions, as 
defined in
terms of  these restricted coactions, are presented in Subsection~IX.C, and
in Section~IX.D the irreducible tensor operators are also defined in terms
of these coactions. It is shown there that, associated with
both $\pi^{R}_{\cal{B}}$ and $\pi^{L}_{\cal{B}}$, there are two types of 
irreducible tensor
operator, which are again called {\em{ordinary}} and {\em{twisted}}, and 
the immediate properties of
all these irreducible tensor operators are described. In Subsection~IX.E it 
is shown that
the irreducible tensor operators satisfy theorems of the Wigner-Eckart 
type, and the analysis is
concluded in Section~IX.F with a demonstration that the products of these
irreducible tensor operators are themselves expressible as linear 
combinations of irreducible tensor
operators that involve the relevant Clebsch-Gordan coefficients. 

Because the space of functions defined on a compact Lie group $\cal{G}$ is 
a special example of a
compact quantum group algebra, all the well-known results for compact 
Lie groups naturally reappear
in this particular case. However, as the detailed analysis shows, the 
theory in the general
situation is rather more subtle, and exhibits various complications.

\section{Properties of compact quantum group algebras}

\subsection{Hopf $*$-algebras}
	The purpose of this subsection is mainly to establish notations, and 
summarize the essential
properties. For further details see, for example, Sweedler$^{27}$, 
Majid$^{28}$, and Chari and
Pressley$^{29}$.

 A {\em{Hopf algebra}} $\cal{A}$ over the field of complex numbers $\C$ 
is a complex vector space
with an identity element $1_{\cal{A}}$ that possesses a multiplication 
operator $M$ (which maps
$\cal{A}
\otimes\cal{A}$ into
$\cal{A}$), a unit operator $u$ (which maps
$\C$ into $\cal{A}$), a comultiplication operator $\Delta$ (which maps 
$\cal{A}$ into $\cal{A}
\otimes\cal{A}$), a counit operator $\epsilon$ (which maps
$\cal{A}$ into $\C$), and an antipode operator $S$ (which maps
$\cal{A}$ into $\cal{A}$). These are assumed to be linear in all their 
arguments and to have the
following properties:
\begin{equation}  M \circ (M \otimes id) = M \circ (id \otimes M) ,
\label{eq:QD8} \end{equation} 
\begin{equation}    (\Delta \otimes id) \circ \Delta =   (id \otimes \Delta) 
\circ \Delta
,\label{eq:QD11} \end{equation}
\begin{equation} \Delta \circ M = (M \otimes M) \circ (id \otimes \sigma 
\otimes id) \circ (\Delta
\otimes
\Delta) ,\label{eq:QD17} \end{equation}
\begin{equation} \epsilon \circ M =  \MC \circ (\epsilon \otimes 
\epsilon), \label{eq:QD15}
\end{equation}
\begin{equation} \MCA \circ (\epsilon \otimes id) \circ \Delta = \MAC 
\circ (id \otimes \epsilon)
\circ \Delta =  id ,
\label{eq:QD13,14a} \end{equation}
\begin{equation}  u(1_{\C}) = 1_{\cal{A}} , \;\; \epsilon(1_{\cal{A}}) = 
1_{\C}, \;\; S(1_{\cal{A}})
= 1_{\cal{A}},
\label{eq:QD7,27a,27g}
\end{equation} 
\begin{equation}  M(a \otimes 1_{\cal{A}}) = M(1_{\cal{A}} \otimes a ) = 
a\;,\;
\mbox{for all $a \in {\cal{A}}$},\label{eq:QD10}
\end{equation}
\begin{equation}  \Delta(1_{\cal{A}}) = 1_{\cal{A}} \otimes 1_{\cal{A}},
\label{eq:QD19}
\end{equation}
\begin{equation}  S \circ M = M \circ \sigma \circ (S \otimes S),
\label{eq:QD21}
\end{equation}
\begin{equation}  \Delta \circ S = (S \otimes S) \circ \sigma \circ 
\Delta,
\label{eq:QD23}
\end{equation} 
\begin{equation}  M \circ ( S \otimes id ) \circ \Delta = M \circ ( id 
\otimes S ) \circ \Delta = u
\circ \epsilon,
\label{eq:QD26}
\end{equation} 
\begin{equation}  \epsilon \circ S = \epsilon.
\label{eq:QD27d}
\end{equation}  Here $\sigma$ is the transposition operator which 
interchanges the order of its
arguments, so that, for example, when acting on $\cal{A} \otimes 
\cal{A}$,
$\sigma(a
\otimes b) = b
\otimes a$ for all $a,b \in \cal{A}$. Also $\MC$, $\MAC$, and $\MCA$ 
are the mutliplication operators
defined by $\MC (w \otimes z) = wz$ for all $w,z \in \C$, and $\MAC (a 
\otimes z) = \MCA (z \otimes
a)= za$ for all $z \in \C$ and all $a \in \cal{A}$. The product $M(a 
\otimes b)$ will sometimes be
written more concisely as $ab$, and the coproduct $\Delta$ will 
sometimes be expressed as
\begin{equation}  \Delta (a) = \sum_{(a)} a_{(1)} \otimes a_{(2)} .
\label{eq:extra1}
\end{equation}

If  $\cal{A}$ is finite-dimensional, with basis elements $a_{1}, a_{2},
\ldots$, the structure constants $m_{jk}^{\ell}$, $\mu^{jk}_{\ell}$, 
$s_{j}^{k}$, $\epsilon_{j}$ and
$\epsilon^{j}$ may be defined by $M(a_{j} \otimes a_{k}) = \sum_{\ell} 
m_{jk}^{\ell}\, a_{\ell}$,
$\Delta(a_{\ell}) = \sum_{j,k} \mu^{jk}_{\ell}\, a_{j} \otimes a_{k}$,
$S(a_{j}) = \sum_{k} s_{j}^{k}\, a_{k}$,
$\epsilon(a_{j}) =  \epsilon_{j}$, and
$1_{\cal{A}} =  \sum_{j} \epsilon^{j} a_{j}$. Then (\ref{eq:QD8}) to 
(\ref{eq:QD27d}) imply that
\begin{equation}  \sum_{s} m_{jk}^{s} m_{s\ell}^{t} = \sum_{s} m_{js}^{t} 
m_{k\ell}^{s},
\label{eq:QD9} \end{equation}
\begin{equation}  \sum_{j} \mu^{jk}_{\ell} \mu^{st}_{j} = \sum_{j} 
\mu^{sj}_{\ell} \mu^{tk}_{j},
\label{eq:QD12} \end{equation}
\begin{equation}  \sum_{p,q,s,t} \mu^{pq}_{j} \mu^{st}_{k} m_{ps}^{r} 
m_{qt}^{u} = \sum_{p}
m_{jk}^{p} \mu^{ru}_{p},
\label{eq:QD18} \end{equation}
\begin{equation}  \sum_{\ell} m_{jk}^{\ell} \epsilon_{\ell}
 = \epsilon_{j} \epsilon_{k},
\label{eq:QD16} \end{equation}
\begin{equation}  \sum_{j} \mu^{jk}_{\ell} \epsilon_{j}
 = \sum_{j} \mu^{kj}_{\ell} \epsilon_{j} = \delta_{\ell}^{k},
\label{eq:QD14,14c} \end{equation}
\begin{equation}  \sum_{j} \epsilon^{j} s_{j}^{k}
 = \epsilon^{k}, \;\sum_{j} \epsilon^{j}\epsilon_{j} =1_{\C},
\label{eq:QD27b,27h} \end{equation}
\begin{equation}  \sum_{k} \epsilon^{k} m_{jk}^{\ell}
 = \sum_{k} \epsilon^{k} m_{kj}^{\ell} = \delta_{j}^{\ell},
\label{eq:QD10a} \end{equation}
\begin{equation}  \sum_{j} \epsilon^{j} \mu^{k\ell}_{j}
 = \epsilon^{k} \epsilon^{\ell},
\label{eq:QD20} \end{equation}
\begin{equation}  \sum_{q} m_{jk}^{q} s_{q}^{p} = \sum_{q,r} m_{rq}^{p} 
s_{j}^{q} s_{k}^{r},
\label{eq:QD22} \end{equation}
\begin{equation}  \sum_{k} \mu^{pq}_{k} s_{j}^{k} = \sum_{k,\ell} 
\mu^{k\ell}_{j} s_{\ell}^{p}
s_{k}^{q} ,
\label{eq:QD25} \end{equation}
\begin{equation}  \sum_{k,\ell,r} \mu^{k\ell}_{j} s_{k}^{r} m_{r\ell}^{t} = 
\sum_{k,\ell,r}
\mu^{k\ell}_{j} s_{\ell}^{r} m_{kr}^{t} = \epsilon_{j} \epsilon^{t},
\label{eq:QD27} \end{equation} and
\begin{equation}  \sum_{j} \epsilon_{j} s_{k}^{j}
 = \epsilon_{k}.
\label{eq:QD27e} \end{equation}

A {\em{Hopf $*$-algebra}} ${\cal{A}}$ is defined to be a Hopf algebra that 
possesses an additional
$*$-operation that maps ${\cal{A}}$ into ${\cal{A}}$. The effect of the $*$ 
operation on $a \in
\cal{A}$ will sometimes be denoted by
$a^{*}$. In particular
\begin{equation}  1_{\cal{A}}^{*} = 1_{\cal{A}}\; . 
\label{eq:extra9} \end{equation} The other properties are:
\begin{equation}  (* \circ \MCA)(z \otimes a) = \MCA (\bar{z} \otimes 
a^{*}) = (* \circ \MAC)(a
\otimes z) = \MAC (a^{*} \otimes \bar{z}) =   \bar{z} a^{*} , 
\label{eq:extra2} \end{equation}  (for all $z \in \C$ and all $a \in 
{\cal{A}}$, where $\bar{z}$
denotes the complex conjugate of $z$),
\begin{equation}  * \circ \;* = id ,
\label{eq:extra3} \end{equation}
\begin{equation}  * \circ\; M = M \circ (* \otimes *) \circ \;\sigma , 
\label{eq:extra4} \end{equation}
\begin{equation}  \Delta \circ * = (* \otimes *) \circ \;\Delta , 
\label{eq:extra5} \end{equation}
\begin{equation}  (\epsilon \circ *)(a) = \overline{\epsilon (a)}\;\;  
\mbox{for all $a \in
\cal{A}$,} \label{eq:extra6} \end{equation}
\begin{equation}  S \circ * \circ S \circ * = id ,  
\label{eq:extra7} \end{equation} which implies that $S$ is invertible with 
inverse given by 
\begin{equation}  S^{-1} = * \circ S \circ *\; .  
\label{eq:extra8} \end{equation} 

If  $\cal{A}$ is finite-dimensional, its linear dual will be denoted by
$\cal{A}^{\prime}$, the prime being used instead of the usual star to 
avoid any confusion with the
$*$-operation that has just been defined. The effect of $a^{\prime} \in
\cal{A}^{\prime}$ on
$a
\in
\cal{A}$ will be denoted by
$\langle a^{\prime},a
\rangle$, and the evaluation map $ev$ (from $\cal{A}^{\prime} \otimes
\cal{A}$ to $\C$) will be defined by
\begin{equation} ev(a^{\prime} \otimes a) = \langle a^{\prime},a \rangle  
\label{eq:extra10} \end{equation} for all $a^{\prime} \in 
\cal{A}^{\prime}$ and all $a \in \cal{A}$.
In the case in which ${\cal{A}}$ is of finite dimension n, the dual basis of
$\cal{A}^{\prime}$ will be denoted by $a^{1}, a^{2}, \ldots, a^{n}$, and 
will be assumed to be such
that
\begin{equation} \langle a^{j},a_{k} \rangle = \delta^{j}_{k} , 
\label{eq:extra11} \end{equation} for all $j,k = 1,2,\ldots, n.$   
 
\subsection{Right comodules of Hopf $*$-algebras}

A {\em{right $\cal{A}$-comodule}} consists of a vector space $V$ and a 
linear mapping $\pi_{V}$ from
$V$ to $V \otimes \cal{A}$ such that
\begin{equation}  (\pi_{V} \otimes id) \circ \pi_{V} = (id \otimes \Delta) 
\circ \pi_{V}
\label{eq:R1} \end{equation} and \begin{equation}  \MVC \circ (id 
\otimes \epsilon) \circ \pi_{V} =
id,
\label{eq:R2} \end{equation}  where $\MVC(v \otimes z) = zv$ for all $v 
\in V$ and all $z \in \C$.
The operation $\pi_{V}$ is then said to be a {\em{right coaction}} and 
provides a
{\em{corepresentation}} of $\cal{A}$ with carrier space $V$.  The present 
subsection will be devoted
to a very brief account of the essential features of the corepresentations of 
$\cal{A}$. (For the
intimate connection between the {\em{corepresentation}} theory of 
$\cal{A}$ and the
{\em{representation}} theory of ${\cal{A}}^{\prime}$, see Appendix A.)

 If $V$ is of finite dimension $d$, with basis $v_{1}, v_{2}, \ldots, v_{d}$, 
then there exists a
uniquely determined set of elements $\pi_{jk}^{V}$ of $\cal{A}$ (for
$j,k = 1,2,\ldots,d$), called the {\em{matrix coefficients}} of $\pi_{V}$, 
which are such that
\begin{equation}  \pi_{V}(v_{j}) = \sum_{k=1}^{d} v_{k} \otimes 
\pi_{kj}^{V}
\label{eq:R5} \end{equation} for all $j = 1,2,\ldots,d$. (In this situation 
the corepresentation is
said to have dimension $d$). The requirements (\ref{eq:R1}) and 
(\ref{eq:R2}) then imply that
\begin{equation}  \Delta(\pi_{jk}^{V}) = \sum_{\ell=1}^{d} \pi_{j\ell}^{V} 
\otimes \pi_{\ell k}^{V}
\label{eq:R6} \end{equation} and
\begin{equation}  \epsilon(\pi_{jk}^{V}) = \delta_{jk}
\label{eq:R7} \end{equation} (for $j,k = 1,2,\ldots,d$). It is sometimes 
convenient to write
\begin{equation}  \pi_{V}(v) = \sum_{[v]} v_{[1]} \otimes v_{[2]} ,
\label{eq:R12} \end{equation} where $v_{[1]} \in V$ and $v_{[2]} \in 
\cal{A}$.

Two right $\cal{A}$-comodules, with carrier spaces $V$ and $W$, 
coactions $\pi_{V}$ and $\pi_{W}$,
and matrix coefficients $\pi^{V}_{jk}$ and $\pi^{W}_{jk}$, are said to be 
{\em{equivalent}} if there
exists a one-to-mapping $\Phi$ from $V$ to $W$ such that
\begin{equation}  \pi_{W} \circ \Phi = (\Phi \otimes id) \circ \pi_{V} .
\label{eq:D.eq1.1.22} \end{equation} If $V$ and $W$ have bases $v_{1}, 
v_{2}, \ldots, v_{d}$ and
$w_{1}, w_{2}, \ldots, w_{d}$ respectively, then to the mapping $\Phi$ 
there corresponds a $d \times
d$ non-singular matrix ${\bf{\Phi}}$ such that
\begin{equation}  \sum_{\ell=1}^{d} \Phi_{j\ell} \pi^{V}_{\ell k} = 
\sum_{\ell=1}^{d}
\pi^{W}_{j\ell} \Phi_{\ell k}
\label{eq:D.eq1.1.23} \end{equation} for all $j,k = 1,2,\ldots,d$.

A subspace $W \subset V$ is said to be {\em{invariant}} under $\pi_{V}$ 
if $\pi_{V}(w) \subset W
\otimes \cal{A}$ for all $w \in W$, and a corepresentation is described as 
being {\em{irreducible}}
if $V$ and ${0}$ are the only invariant subspaces of $V$. If $V$ is the 
direct sum of two
invariant subspaces of $V$, then the
corepresentation $\pi_{V}$ is said to be {\em{completely reducible}}. 

If $V$ is endowed with an inner product $\langle \; ,\; \rangle_{V}$ (such 
that
$\langle zw,z^{\prime}v \rangle_{V} = \bar{z}z^{\prime}\langle w,v 
\rangle_{V}$ for all
$z,z^{\prime} \in \C$ and all $v,w \in V$), then
$\pi_{V}$ is  said to give a {\em{unitary}} corepresentation if
\begin{equation} \sum_{[v]} \langle w, v_{[1]} \rangle_{V}\; S(v_{[2]}) = 
\sum_{[w]} \langle
w_{[1]},v \rangle_{V}\; w_{[2]}^{*}
\label{eq:D.sect1.1.45i} \end{equation} for all $v,w \in V$. It can be 
shown$^{22-26}$ that if
$v_{1}, v_{2},
\ldots, v_{d}$ is an ortho-normal basis of $V$ then
\begin{equation} S(\pi_{jk}^{V}) = \pi_{kj}^{V*} ,
\label{eq:D.sect1.1.45ii} \end{equation}
\begin{equation} \sum_{\ell=1}^{d} M(\pi_{\ell j}^{V*} \otimes \pi_{\ell 
k}^{V}) = \delta _{jk}
1_{\cal{A}}  ,
\label{eq:D.sect1.1.45iii} \end{equation} and
\begin{equation} \sum_{\ell=1}^{d} M(\pi_{j\ell}^{V} \otimes 
\pi_{k\ell}^{V*}) = \delta _{jk}
1_{\cal{A}}  
\label{eq:D.sect1.1.45iv} \end{equation} (for all $j,k = 1,2,\ldots,d$).

 Corresponding to a right $\cal{A}$-comodule with carrier space $V$ and 
coaction
$\pi_{V}$ from $V$ to $V \otimes \cal{A}$ there exist two other right 
$\cal{A}$-comodules formed
from the same carrier space. Firstly, there is the coaction 
$\pi_{V}^{\ddagger}$, which is said to be
{\em{doubly contragredient}} to $\pi_{V}$, and which is defined (as a 
mapping from $V$ to $V
\otimes \cal{A}$) by
\begin{equation} \pi_{V}^{\ddagger} = (id \otimes S^{2}) \circ \pi_{V} .
\label{eq:R144} \end{equation} With the matrix coefficients 
$\pi_{jk}^{V\ddagger}$ of
$\pi_{V}^{\ddagger}$ being defined  by 
\begin{equation}  \pi_{V}^{\ddagger}(v_{j}) = \sum_{k=1}^{d} v_{k} 
\otimes \pi_{kj}^{V\ddagger}
\label{eq:R146} \end{equation} for all $j = 1,2,\ldots,d$, it follows from 
(\ref{eq:R5}) that 
\begin{equation} \pi^{V\ddagger}_{jk} = S^{2}(\pi_{jk}^{V}) .
\label{eq:R147} \end{equation} for all $j,k = 1,2,\ldots,d$. Secondly, let 
$\bar{V}$ be the
conjugate space to $V$ (so that as an Abelian group $\bar{V}$ is 
isomorphic to $V$, but the
scalar multiplication operator $\MCVbar$ for $\bar{V}$ is defined in 
terms of the corresponding
operator $\MCV$ for $V$ by $\MCVbar(z \otimes v) = \MCV(\bar{z} 
\otimes v))$. Then the coaction
$\bar{\pi}_{\bar{V}}$, which is said to be {\em{conjugate}} to $\pi_{V}$, is 
defined (as a mapping
from
$\bar{V}$ to $\bar{V} \otimes \cal{A}$) by
\begin{equation} \bar{\pi}_{\bar{V}} = (id \otimes *) \circ \pi_{V} ,
\label{eq:extra17} \end{equation} so its matrix coefficients 
$\bar{\pi}^{\bar{V}}_{jk}$ are given by
\begin{equation} \bar{\pi}^{\bar{V}}_{jk} = \pi_{jk}^{V*} 
\label{eq:extra18} \end{equation} for all $j,k = 1,2,\ldots,d$.

\subsection{Compact quantum group algebras}

 A {\em{compact quantum group algebra}} (or {\em{CQG algebra}} for 
short) may be defined$^{22-26}$
as a Hopf $*$-algebra that is spanned by the matrix coefficients of its non-
equivalent
finite-dimensional unitary irreducible corepresentations.  Koornwinder 
and Dijkhuizen$^{22-26}$ have
shown that if
$\cal{A}$ is a CQG algebra then every finite-dimensional corepresentation 
of $\cal{A}$ is equivalent
to a unitary corepresentation, and that every finite-dimensional reducible 
corepresentation of
$\cal{A}$ is completely reducible. Moreover$^{22-26}$ if
$\cal{A}$ is a CQG algebra then $\cal{A}$ possesses a {\em{Haar 
functional}}, $h$, which is a
mapping of
$\cal{A}$ into $\C$ such that
\begin{equation} h(1_{\cal{A}}) = 1_{\C} ,  
\label{eq:D.def2.1.9a} \end{equation}
\begin{equation} h(M(a^{*} \otimes a) > 0 ,  
\label{eq:D.th2.1.16} \end{equation}
\begin{equation} h(a^{*}) = \overline{h(a)} ,  
\label{eq:D.lemma2.1.10} \end{equation}
\begin{equation} h(S(a)) = h(a) ,  
\label{eq:D.lemma2.1.11} \end{equation} and 
\begin{equation} (\MCA \circ (h \otimes id) \circ \Delta)(a) = (\MAC 
\circ (id \otimes h) \circ
\Delta)(a) = h(a)\;1_{\cal{A}}  
\label{eq:D.def2.1.9b} \end{equation} for all $a \in \cal{A}$.

Koornwinder and Dijkhuizen$^{22-26}$ have also shown that if $\pi^{p}$ 
and $\pi^{q}$ are two
non-equivalent irreducible corepresentations of a CQG algebra
$\cal{A}$ with dimensions $d_{p}$ and $d_{q}$ and matrix coefficients 
$\pi^{p}_{jk}$ and
$\pi^{q}_{mn}$ respectively, then
\begin{equation} h(M(\pi^{p}_{jk} \otimes S(\pi^{q}_{mn})) = 
0,\;\;h(M(S(\pi^{p}_{jk}) \otimes
\pi^{q}_{mn}) = 0,  
\label{eq:D.prop2.1.3} \end{equation} (for all $j,k = 1,2,\ldots,d_{p}$ and 
for all $m,n =
1,2,\ldots,d_{q}$). Moreover every $\pi^{p}$ irreducible corepresentation 
of $\cal{A}$ is equivalent
to its doubly contragredient partner
$\pi^{p\ddagger}$, so in each such case there exists a non-singular  
$d_{p} \times d_{p}$ matrix
${\bf{F}}^{p}$ such that
\begin{equation} \sum_{k=1}^{d_{p}} F_{jk}^{p} \pi_{k\ell}^{p} = 
\sum_{k=1}^{d_{p}}
\pi_{jk}^{p\ddagger} F_{k\ell}^{p}
\label{eq:R148} \end{equation} (for all $j,\ell = 1,2,\ldots,d_{p}$). Then, if 
$\pi^{p}$ is a
unitary irreducible corepresentation of $\cal{A}$,
\begin{equation} h(M(\pi^{p}_{jk} \otimes S(\pi^{p}_{mn}))) =
\delta_{jn}F_{mk}^{p}/tr({\bf{F}}^{p}) ,  
\label{eq:D.eq2.1.15} \end{equation} and
\begin{equation} h(M(S(\pi^{p}_{jk}) \otimes \pi^{p}_{mn})) =
\delta_{jn}(({\bf{F}}^{p})^{-1})_{mk}/tr(({\bf{F}}^{p})^{-1})   
\label{eq:D.eq2.1.16} \end{equation} (for all $j,k,m,n = 1,2,\ldots,d_{p}$).

Of course in the special case in which $\cal{A}$ is the space of functions 
defined on a compact
group ${\cal{G}}$, $\cal{A}$ is commutative (i.e. $M = M \circ \sigma$), 
$S^{2} = id$, $h$ is the
Haar integral
\begin{equation} h(a) = \int_{\cal{G}} a(x)\;\mbox{d}x ,  
\label{eq:extra13} \end{equation} and (\ref{eq:D.def2.1.9b}) express the 
invariance properties
\begin{equation} \int_{\cal{G}} a(yx)\;\mbox{d}x = \int_{\cal{G}} 
a(x)\;\mbox{d}x= \int_{\cal{G}}
a(xy)\;\mbox{d}x ,  
\label{eq:extra14} \end{equation} (for all $y \in {\cal{G}}$). Moreover in 
this case each
corepresentation of $\cal{A}$ is {\em{identical}} to its doubly 
contragredient partner, and
(\ref{eq:D.prop2.1.3}), (\ref{eq:D.eq2.1.15}), and (\ref{eq:D.eq2.1.16}) 
correspond to the
orthogonality theorems for the unitary irreducible representations of 
${\cal{G}}$.

\subsection{Lemmas concerning the Haar functional}

It will now be shown that
\begin{equation} \sum_{(b)} h(M(a \otimes b_{(1)}))\, S(b_{(2)}) = 
\sum_{(a)} h(M(a_{(1)} \otimes
b))\, a_{(2)}   
\label{eq:R43} \end{equation} for all $a,b \in \cal{A}$.

To prove this consider the operation $u \circ h \circ M$. As the left-hand 
equality of
(\ref{eq:D.def2.1.9b}) can be re-expressed as 
\begin{equation} \MCA \circ (h \otimes id) \circ
\Delta = u \circ h ,\label{eq:R47} \end{equation} it follows from 
(\ref{eq:QD17}) that
\begin{equation}  u \circ h \circ M = \MCA \circ (h \otimes id) \circ (M 
\otimes M)
\circ (id \otimes \sigma \otimes id) \circ (\Delta \otimes \Delta) .
\label{eq:R50} \end{equation} However, by (\ref{eq:QD13,14a}), 
(\ref{eq:QD26}), and (\ref{eq:QD11}),
it also follows that
\begin{equation}  u \circ h \circ M = \MCA \circ (h \otimes M) \circ (M 
\otimes S
\otimes id)
\circ (id \otimes \Delta \otimes id) \circ (id \otimes \Delta) .
\label{eq:R53} \end{equation} Comparing (\ref{eq:R50}) and 
(\ref{eq:R53}) then gives
\[ \begin{array}{l} \MCA \circ (h \otimes M) \circ (M \otimes S
\otimes id) \circ (id \otimes \Delta \otimes id) \circ (id \otimes \Delta) 
\\
 =\MCA \circ (h \otimes id) \circ (M \otimes M) \circ (id \otimes \sigma 
\otimes id) \circ (\Delta
\otimes \Delta) .
 \end{array} \] 
which can be re-expressed as
\begin{equation}  M \circ (U \otimes id) \circ (id \otimes \Delta) = M 
\circ (V \otimes id) \circ
(id \otimes \Delta) \;,\label{eq:R54a}\end{equation}
where
\begin{equation}  U = \MCA \circ (h \otimes id) \circ (M \otimes id) \circ 
(id \otimes \sigma)
\circ (\Delta \otimes id)\end{equation}
and 
\begin{equation}  V = \MCA \circ (h \otimes S) \circ (M \otimes id) \circ 
(id
\otimes \Delta)\;.\end{equation}
However, (\ref{eq:R54a}) implies that
\begin{equation} \begin{array}{l}  M \circ (M \otimes id) \circ (U \otimes 
id \otimes S) \circ (id
\otimes
\Delta
\otimes id) \circ (id \otimes \Delta) \\
= M \circ (M \otimes id) \circ (V \otimes id \otimes S) \circ
(id \otimes \Delta
\otimes id) \circ (id \otimes \Delta) \;.\end{array} \label{eq:R54c} 
\end{equation}
But
\[ \begin{array}{l}  M \circ (M \otimes id) \circ (U \otimes id \otimes S) 
\circ (id \otimes \Delta
\otimes id) \circ (id \otimes \Delta)\\= 
M \circ (U \otimes \{M \circ (id \otimes S) \circ \Delta\}) 
\circ (id \otimes \Delta)\\=
M \circ (U \otimes \{u \circ \epsilon\}) 
\circ (id \otimes \Delta) = U.\end{array} \]
As a similar result is true with $U$ replace by $V$, it follows from 
(\ref{eq:R54c}) that $U = V$,
which is an equivalent way of expressing (\ref{eq:R43}). 

It can shown by a similar argument that
\begin{equation} \sum_{(b)} h(M(b_{(2)} \otimes a))\, S(b_{(1)}) = 
\sum_{(a)} h(M(b \otimes
a_{(2)}))\, a_{(1)}   
\label{eq:R60} \end{equation} for all $a,b \in \cal{A}$.

\section{The right and left regular comodules}

The {\em{right regular comodule}} of $\cal{A}$ is defined to have 
$\cal{A}$ itself as its carrier
space, with $\Delta$ providing the coaction $\pi^{R}_{\cal{A}}$. That is
\begin{equation}  V = {\cal{A}} , \;\pi^{R}_{\cal{A}} = \Delta .
\label{eq:R8} \end{equation}  In this case the conditions (\ref{eq:R1}) and 
(\ref{eq:R2}) for
$\pi^{R}_{\cal{A}}$ to be a right coaction are  immediately satisfied by 
virtue of the assumptions 
(\ref{eq:QD11}) and (\ref{eq:QD13,14a}).

The {\em{left regular comodule}} of $\cal{A}$ is also defined to have 
$\cal{A}$ itself as its
carrier space, but has $\sigma \circ (S \otimes id) \circ \Delta$ as its 
coaction
$\pi^{L}_{\cal{A}}$. That is
\begin{equation}  V = {\cal{A}} ,\;\pi^{L}_{\cal{A}} = \sigma \circ (S 
\otimes id) \circ \Delta .
\label{eq:L1} \end{equation}  For this case the condition (\ref{eq:R1}) for 
$\pi^{L}_{\cal{A}}$ to
be a right coaction is satisfied by the assumptions  (\ref{eq:QD11}) and 
(\ref{eq:QD23}), while the
condition (\ref{eq:R2}) is again satisfied as a result of the assumption  
(\ref{eq:QD13,14a}). It
should be noted that $\pi^{R}_{\cal{A}}$ and $\pi^{L}_{\cal{A}}$ are both 
{\em{right}} coactions,
for, as discussed in  Appendix A, the designation `left' of
$\pi^{L}_{\cal{A}}$ comes from its relation to the left regular action of a 
group in the special 
case in which the dual ${\cal{A}}^{\prime}$ is a group algebra. It is also 
useful to note that 
(\ref{eq:L1}) implies that
\begin{equation}  \Delta = (S^{-1} \otimes id) \circ \sigma \circ
\pi^{L}_{\cal{A}} .
\label{eq:L25} \end{equation}

The notation of (\ref{eq:R12}) can be developed further by writing
\begin{equation}  \pi^{X}_{\cal{A}}(a) = \sum_{[a]} a^{X}_{[1]} \otimes 
a^{X}_{[2]} ,
\label{eq:extra15} \end{equation} for $X = R$ and $X = L$, where 
$a^{X}_{[1]}$ and $a^{X}_{[1]}$ are
elements of
$\cal{A}$. Then (\ref{eq:extra1}) and (\ref{eq:R8}) imply that
\begin{equation}  a^{R}_{[1]} = a_{(1)} ,\; a^{R}_{[2]} = a_{(2)} ,
\label{eq:extra15a} \end{equation}  but (\ref{eq:extra1}) and (\ref{eq:L1}) 
give
\begin{equation}  a^{L}_{[1]} = a_{(2)} ,\;a^{L}_{[2]} = S(a_{(1)}) .
\label{eq:L5} \end{equation}

The right and left regular corepresentations are both {\em{unitary}}, 
provided that the inner
products on the carrier space $\cal{A}$ are chosen in the following way:
\begin{enumerate}
\item for the {\em{right}} regular corepresentation take
\begin{equation}  \langle a , b \rangle_{\cal{A}} = (a,b)^{R} = h(M(a^{*} 
\otimes b))\;\;{\mbox{for
all $a,b
\in
\cal{A}$}};
\label{eq:R40} \end{equation}
\item for the {\em{left}} regular corepresentation take
\begin{equation}  \langle a , b \rangle_{\cal{A}} = (a,b)^{L} = h(M(b 
\otimes (S^{2}(a))^{*}))\;\;
{\mbox{for all $a,b \in \cal{A}$}} .
\label{eq:L14} \end{equation}
\end{enumerate}

In outline the proofs of these statements are as follows. For the 
{\em{right}} regular
corepresentation, the unitary condition (\ref{eq:D.sect1.1.45i}) with the 
choice (\ref{eq:R40}) for
inner product and with (\ref{eq:extra15a}) becomes 
\[ \sum_{(v)} h(M(w^{*} \otimes v_{(1)}))\, S(v_{(2)}) = \sum_{(w)} 
h(M(w_{(1)}^{*}
\otimes v))\, w_{(2)}^{*} , \]    which in turn reduces to (\ref{eq:R43}) with 
the substitutions $w =
a^{*}$ and $v = b$. Similarly, for the {\em{left}} regular corepresentation, 
the unitary condition
(\ref{eq:D.sect1.1.45i}) with the choice (\ref{eq:L14}) for inner product and 
with (\ref{eq:L5})
becomes 
\[ \sum_{(v)} h(M(v_{(2)}  \otimes (S^{2}(w))^{*}))\, S^{2}(v_{(1)}) = 
\sum_{(w)} h(M(v
\otimes (S^{2}(w_{(2)}))^{*}))\, (S(w_{(1)}))^{*} . \]    With the 
substitutions $w = S^{-1}(a^{*})$
and $v = S^{-1}(b)$, and the application of (\ref{eq:QD23}), 
(\ref{eq:extra5}), and
(\ref{eq:extra7}), this reduces again to (\ref{eq:R43}).
 
With the choices (\ref{eq:R40}) and (\ref{eq:L14}), both $(a,a)^{R}$ and 
$(a,a)^{L}$ are real and
positive for all non-zero $a \in \cal{A}$. For $(a,a)^{R}$ this is an 
immediate consequence of
(\ref{eq:D.th2.1.16}), while for $(a,a)^{L}$ the demonstration requires 
(\ref{eq:D.lemma2.1.11}),
(\ref{eq:QD21}), (\ref{eq:D.lemma2.1.10}), and (\ref{eq:extra7}) as well.

The inner products (\ref{eq:R40}) and (\ref{eq:L14}) will be used 
throughout this paper. In the
special case in which $\cal{A}$ is the space of functions defined on a 
compact group ${\cal{G}}$,
both $(a,b)^{R}$ and
$(a,b)^{L}$ reduce to the usual inner product
\[ \int_{{\cal{G}}} \overline{a(x)} b(x) \;\mbox{d}x . \]

It is worth noting at this stage that the invariance properties 
(\ref{eq:D.def2.1.9b}) imply that 
\begin{equation} (\MCA \circ (h \otimes id) \circ \pi^{R}_{\cal{A}})(a) = 
(\MAC \circ (id \otimes h)
\circ
\pi^{R}_{\cal{A}})(a) = h(a)\;1_{\cal{A}}  
\label{eq:R78a,R78b} \end{equation} and 
\begin{equation} (\MCA \circ (h \otimes id) \circ \pi^{L}_{\cal{A}})(a) = 
(\MAC \circ (id \otimes h)
\circ
\pi^{L}_{\cal{A}})(a) = h(a)\;1_{\cal{A}}  
\label{eq:L19a,L19b} \end{equation} for all $a \in \cal{A}$. Acting with 
$h$ again in
(\ref{eq:R78a,R78b}) and (\ref{eq:L19a,L19b}), and using 
(\ref{eq:D.def2.1.9a}), gives
\begin{equation} (\MC \circ (h \otimes h) \circ \pi^{X}_{\cal{A}})(a)  = 
h(a)  
\label{eq:R78c,L19d} \end{equation} for both $X = R$ and $L$ and for all 
$a \in \cal{A}$. In terms
of the notation of (\ref{eq:extra15}), this can be re-expressed as
\begin{equation}  h(a) = \sum_{[a]} h(a^{X}_{[1]}) \,h(a^{X}_{[2]}) ,
\label{eq:extra15b} \end{equation} for $X = R$ and $L$ and for all $a \in 
\cal{A}$.

The effects of the right and left regular coactions on products are quite 
different. For the right
regular coaction, (\ref{eq:R8}) and (\ref{eq:QD17}) imply immediately that
\begin{equation}   \pi^{R}_{\cal{A}} \circ M = (M \otimes M) \circ (id 
\otimes \sigma \otimes id)
\circ (\pi^{R}_{\cal{A}} \otimes \pi^{R}_{\cal{A}}) ,
\label{eq:R164v} \end{equation} whereas for the left regular coaction, 
(\ref{eq:L1}),
(\ref{eq:QD17}), and (\ref{eq:QD21}) show that
\begin{equation}   \pi^{L}_{\cal{A}} \circ M = (M \otimes M) \circ (id 
\otimes id \otimes \sigma)
\circ (id \otimes \sigma
\otimes id) \circ (\pi^{L}_{\cal{A}} \otimes \pi^{L}_{\cal{A}}) ,
\label{eq:L73} \end{equation} which contains an extra twist factor 
$\sigma$.
  
\section{Basis functions}

\subsection{Definitions and properties}

Suppose that $\pi^{p}_{jk}$ are the matrix coefficients of a 
corepresentation $\pi^{p}$ of
$\cal{A}$ of finite dimension $d_{p}$. Then the {\em{basis functions}} 
$\psi^{pR}_{j}$ {\em{of}}
$\pi^{p}$ {\em{with respect to the right regular coaction}} may be defined 
to be a set of $d_{p}$
elements of
$\cal{A}$ that have the property that   
\begin{equation} \pi^{R}_{{\cal{A}}}(\psi^{pR}_{j}) = \sum_{k=1}^{d_{p}} 
\psi^{pR}_{k} \otimes
\pi^{p}_{kj}
\label{eq:R79} \end{equation} for all $j = 1, 2, \ldots, d_{p}$. Similarly the 
{\em{basis
functions}} $\psi^{pL}_{j}$ {\em{of}} $\pi^{p}$ {\em{with respect to the 
left regular coaction}} may
be defined as a set of $d_{p}$ elements of $\cal{A}$ that have the property 
that   
\begin{equation} \pi^{L}_{{\cal{A}}}(\psi^{pL}_{j}) = \sum_{k=1}^{d_{p}} 
\psi^{pL}_{k} \otimes
\pi^{p}_{kj}
\label{eq:L20} \end{equation} for all $j = 1, 2, \ldots, d_{p}$.

In the right regular coaction case, an example of a set of basis functions is 
provided (for any
fixed choice of $\ell = 1, 2, \ldots$) by
\begin{equation} \psi^{pR}_{j} = \pi^{p}_{\ell j}
\label{eq:R80} \end{equation} for all $j = 1, 2, \ldots, d_{p}$. Likewise, in 
the left regular
coaction case, an example is provided (for any fixed choice of $\ell = 1, 2, 
\ldots$) by
\begin{equation} \psi^{pL}_{j} = S^{-2}(\pi^{p*}_{j\ell})
\label{eq:L22} \end{equation} for all $j = 1, 2, \ldots, d_{p}$ (provided that 
the copresentation
$\pi^{p}$ is unitary).

One very useful result, which comes from applying (\ref{eq:L20}), 
(\ref{eq:L25}), (\ref{eq:QD23}),
(\ref{eq:extra5}), and (\ref{eq:extra7}), is that
\begin{equation} \pi^{L}_{{\cal{A}}}((S^{2}(\psi^{qL}_{k}))^{*}) = 
\sum_{t=1}^{d_{p}}
(S^{2}(\psi^{qL}_{t}))^{*} \otimes \pi^{q*}_{tk} 
\label{eq:L90} \end{equation} for all $k = 1, 2, \ldots, d_{p}$.

In spite of the fact that the inner products (\ref{eq:R40}) and 
(\ref{eq:L14}) for the right and
left regular coactions are different, in both the cases the basis functions 
possess the {\em{same}}
orthogonality properties, which are as follows: If
$\psi^{qX}_{k}$ and
$\phi^{pX}_{j}$ are basis functions of the unitary irreducible 
corepresentations $\pi^{q}$ and
$\pi^{p}$ of $\cal{A}$, then  
\begin{equation} (\psi^{qX}_{k},\phi^{pX}_{j})^{X} = 0\;\;{\mbox{unless 
$p = q$ and $j = k$,}}
\label{eq:R81c} \end{equation} and  
\begin{equation} (\psi^{pX}_{j},\phi^{pX}_{j})^{X} \;\;{\mbox{is 
independent of $j$, for $j =
1,2,\ldots,d_{p}$.}}
\label{eq:R81d} \end{equation} Here $X$ denotes both $R$ and $L$, and 
in (\ref{eq:R81d})
$\psi^{pX}_{j}$ and $\phi^{pX}_{j}$ need not be identical sets. Indeed, 
with $X = R$, if the
functions $\psi^{pR}_{j}$ and
$\phi^{pR}_{j}$ are defined by 
\begin{equation} \psi^{pR}_{j} = \pi^{p}_{sj} \; ,\; \phi^{pR}_{j} = 
\pi^{p}_{tj} ,
\label{eq:R81e,R81f} \end{equation} and with $X = L$, if the functions 
$\psi^{pL}_{j}$ and
$\phi^{pL}_{j}$ are similarly defined by 
\begin{equation} \psi^{pL}_{j} = S^{-2}(\pi^{p*}_{js}) \; ,\; \phi^{pL}_{j} = 
S^{-2}(\pi^{p*}_{jt}) ,
\label{eq:L31,L32} \end{equation} then in {\em{both}} cases
\begin{equation} (\psi^{pX}_{j},\phi^{pX}_{j})^{X} =
(({\bf{F}}^{p})^{-1})_{ts}/tr(({\bf{F}}^{p})^{-1}) 
\label{eq:R81g} \end{equation} for $j =1, 2, \ldots, d_{p}$. 

The proofs of (\ref{eq:R81c}), (\ref{eq:R81d}), and (\ref{eq:R81g}) will now 
be outlined. Applying
(\ref{eq:R164v}) and (\ref{eq:R78a,R78b}) to $M(\psi^{qR*}_{k} \otimes 
\phi^{pR}_{j})$ gives
\[ h(M(\psi^{qR*}_{k} \otimes \phi^{pR}_{j}))\; 1_{{\cal{A}}} = \sum_{s,t} 
h(M(\psi^{qR*}_{t} \otimes
\phi^{pR}_{s}))\; M(\pi^{q*}_{tk} \otimes \pi^{p}_{sj}) .\] A further 
application of $h$ to both
sides gives
\begin{equation} (\psi^{qX}_{k} , \phi^{pX}_{j})^{X} = \sum_{s,t} 
(\psi^{qX}_{t} ,
\phi^{pX}_{s})^{X} \; h(M(\pi^{q*}_{tk} \otimes \pi^{p}_{sj})) 
\label{eq:R81prime} \end{equation} with $X = R$, where (\ref{eq:R40}) 
and (\ref{eq:D.def2.1.9a})
have been invoked. Similarly, applying (\ref{eq:L73}) and (\ref{eq:L90}) to 
$M(\phi^{pL}_{j} \otimes
(S^{2}(\psi^{qL}_{k}))^{*})$, and then applying (\ref{eq:L19a,L19b}) to the 
result gives
\[ h(M(\phi^{pL}_{j} \otimes (S^{2}(\psi^{qL}_{k}))^{*}))\; 1_{{\cal{A}}} = 
\sum_{s,t}
h(M(\phi^{pL}_{s}
\otimes (S^{2}(\psi^{qL}_{t}))^{*}))\; M(\pi^{q*}_{tk} \otimes \pi^{p}_{sj}) 
.\] A further
application of $h$ to both sides gives (\ref{eq:R81prime}) with $X = L$, 
where this time
(\ref{eq:L14}) and (\ref{eq:D.def2.1.9a}) have been used. Thus, in both 
cases, it follows from
(\ref{eq:R81prime}), (\ref{eq:D.sect1.1.45ii}), (\ref{eq:D.prop2.1.3}), and 
(\ref{eq:D.eq2.1.16})
that
\begin{equation} (\psi^{qX}_{k},\phi^{pX}_{j})^{X} =
\sum_{s,t} (\psi^{qX}_{t},\phi^{pX}_{s})^{X}\delta^{qp} \delta_{kj}
(({\bf{F}}^{p})^{-1})_{st}/tr(({\bf{F}}^{p})^{-1}) .
\label{eq:R81a} \end{equation} This implies that (\ref{eq:R81c}) is true, 
and if $j = k$ and $p =
q$, then (\ref{eq:R81a}) and (\ref{eq:R81c}) together give
\begin{equation} (\psi^{pX}_{j},\phi^{pX}_{j})^{X} =
\sum_{s} (\psi^{pX}_{s},\phi^{pX}_{s})^{X} (({\bf{F}}^{p})^{-
1})_{ss}/tr(({\bf{F}}^{p})^{-1}) .
\label{eq:R81b} \end{equation} As the right-hand side of (\ref{eq:R81b}) is 
independent of $j$, so
too must be the left-hand side, which thereby establishes (\ref{eq:R81d}). 
Finally the combination of
(\ref{eq:R40}) and (\ref{eq:R81e,R81f}) and the combination of 
(\ref{eq:L14}) and (\ref{eq:L31,L32})
both produce the  result
\[ (\psi^{pX}_{j},\phi^{pX}_{j})^{X} = h (M(S(\pi^{p}_{js}) \otimes 
\pi^{p}_{tj}))
  \] for $X = R$ and for $X = L$, which gives (\ref{eq:R81g}) when 
(\ref{eq:D.eq2.1.16}) is used.

\subsection{Projection operators}

 The argument in  Appendix C suggests the following
definition. Let $\pi^{p}$ be a unitary irreducible corepresentation of 
${\cal{A}}$ of dimension
$d_{p}$ with matrix coefficients
$\pi^{p}_{mn}$. Then the projection operators ${\cal{P}}^{pR}_{mn}$ and 
${\cal{P}}^{pL}_{mn}$ are
defined by   
\begin{equation} {\cal{P}}^{pX}_{mn}(a) = d_{p} \sum_{[a]}\; a^{X}_{[1]}\; 
h(M(\pi^{p*}_{mn} \otimes
a^{X}_{[2]})
\label{eq:R143a} \end{equation} for $X = R$ and $X = L$, for all $m,n = 1, 
2, \ldots, d_{p}$, and
all $a \in \cal{A}$.

These projection operators have the following two useful properties: Let 
$\pi^{p}$ and $\pi^{q}$ be
two unitary irreducible corepresentations of ${\cal{A}}$ of dimensions 
$d_{p}$ and $d_{q}$ with
matrix coefficients $\pi^{p}_{mn}$ and $\pi^{q}_{jk}$.
\begin{enumerate}
\item  Then
\begin{equation} {\cal{P}}^{pX}_{mn}\;{\cal{P}}^{qX}_{jk} = d_{p}\; 
\delta^{pq}\;
\{(({\bf{F}}^{p})^{-1})_{nj}/tr(({\bf{F}}^{p})^{-1})\}\;{\cal{P}}^{pX}_{mk}
\label{eq:R152} \end{equation} for $X = R$ and $X = L$, for all $m,n = 1, 
2, \ldots, d_{p}$, and for
all $j,k = 1, 2, \ldots, d_{q}$.
\item Also, if $\psi^{qX}_{k}$ are basis functions for $\pi^{q}$, then 
\begin{equation} {\cal{P}}^{pX}_{mn}(\psi^{qX}_{k}) = d_{p}\; 
\delta^{pq}\; \delta_{nk}\;
\sum_{\ell =1}^{d_{p}} \psi^{qX}_{\ell}\;(({\bf{F}}^{p})^{-1})_{\ell 
m}/tr(({\bf{F}}^{p})^{-1})
\label{eq:R153} \end{equation} for $X = R$ and $X = L$, for all $m,n = 1, 
2, \ldots, d_{p}$, and for
all $k = 1, 2, \ldots, d_{q}$.
\end{enumerate}

The proof of (\ref{eq:R152}) is as follows. For any $f \in {\cal{A}}$, 
(\ref{eq:R143a})  and
(\ref{eq:QD11}) imply that
\[ {\cal{P}}^{pX}_{mn}({\cal{P}}^{qX}_{jk} (f)) =
d_{p}d_{q}\;\sum_{[f]}\;f^{X}_{[1]}\;h(M(\pi^{q*}_{jk} \otimes 
(f^{X}_{[2]})_{(1)}))
\;h(M(\pi^{p*}_{mn}
\otimes (f^{X}_{[2]})_{(2)})) .
\] On applying (\ref{eq:R60}) with $a = f^{X}_{[2]}$ and $b = 
\pi^{q*}_{jk}$, this reduces  
\[ {\cal{P}}^{pX}_{mn}({\cal{P}}^{qX}_{jk} (f)) = 
d_{p}d_{q}\;\sum_{[f]}\;\sum_{\ell =
1}^{d_{q}}\;f^{X}_{[1]}\;h(M(\pi^{q*}_{\ell k} \otimes f^{X}_{[2]}))
\;h(M(\pi^{p*}_{mn}
\otimes S(\pi^{q*}_{j\ell})) ,
\] so
\[ {\cal{P}}^{pX}_{mn}({\cal{P}}^{qX}_{jk} (f)) = d_{p}\;\sum_{\ell =
1}^{d_{q}}\;{\cal{P}}^{qX}_{\ell k}(f)
\;h(M(\pi^{p*}_{mn}
\otimes S(\pi^{q*}_{j\ell})) ,
\] which, by (\ref{eq:D.sect1.1.45ii}), (\ref{eq:QD21}), and 
(\ref{eq:D.eq2.1.16}), gives
(\ref{eq:R152}).

To prove (\ref{eq:R153}) it suffices to note that (\ref{eq:R143a}), 
(\ref{eq:R79}) and
(\ref{eq:L14}) imply that
\[ {\cal{P}}^{pX}_{mn}(\psi^{qX}_{k}) = d_{p}\;
\sum_{\ell =1}^{d_{p}} \psi^{qX}_{\ell}\;h(M(\pi^{p*}_{mn} \otimes 
\pi^{q}_{\ell k})) ,
\] which, by (\ref{eq:D.eq2.1.16}), leads immediately to (\ref{eq:R152}).

\section{Tensor products and Clebsch-Gordan coefficients}

\subsection{Ordinary and twisted tensor products}

The {\em{tensor product}} of two corepresentations
$\pi_{V}$ and $\pi_{W}$ of $\cal{A}$ (with carrier spaces $V$ and $W$ 
respectively) is the mapping
$\pi_{V} \sqtimes \pi_{W}$ from $V \otimes W$ to $V \otimes W \otimes 
\cal{A}$ that is defined by
\begin{equation} (\pi_{V} \sqtimes \pi_{W})(v \otimes w) = ((id \otimes id 
\otimes M) \circ (id
\otimes \sigma \otimes id) \circ (\pi_{V} \otimes \pi_{W}))(v \otimes w)
\label{eq:R164j} \end{equation} for all $v \in V$ and all $w \in W$. It is 
easily checked that the
conditions (\ref{eq:R1}) and (\ref{eq:R2}) are satisfied with $\pi_{V}$ 
replaced by $\pi_{V}
\sqtimes \pi_{W}$ and $V$ replaced by $V
\otimes W$, so $\pi_{V} \sqtimes \pi_{W}$ is indeed a coaction of 
$\cal{A}$ with carrier space $V
\otimes W$.  If $V$ and $W$ are of finite dimensions $d_{V}$ and 
$d_{W}$, with bases
$v_{1}, v_{2}, \ldots, v_{d_{V}}$, and $w_{1}, w_{2}, \ldots, w_{d_{W}}$, 
then (\ref{eq:R5}) and
(\ref{eq:R164j}) give
\begin{equation} (\pi_{V} \sqtimes \pi_{W})(v_{j} \otimes w_{k}) = 
\sum_{s = 1}^{d_{V}}\sum_{t =
1}^{d_{W}} (v_{s} \otimes w_{t}) \otimes (M(\pi_{sj}^{V} \otimes 
\pi_{tk}^{W})) ,
\label{eq:extra16} \end{equation} which implies that the matrix 
coefficients of $\pi_{V} \sqtimes
\pi_{W}$ are given by 
\begin{equation} (\pi_{V} \sqtimes \pi_{W})_{st,jk} = M(\pi_{sj}^{V} 
\otimes \pi_{tk}^{W}) .
\label{eq:R155} \end{equation} Henceforth it will always be assumed that 
in tensor product matrices
such as $\pi_{V} \sqtimes
\pi_{W}$ the pair of indices that specify the columns have the ordering: 
\begin{equation}(j,k) = (1,1), (1,2),
\ldots, (1,d_{W}), (2,1), (2,2), \ldots, \label{eq:extra19} \end{equation} and 
that the same
ordering applies to the rows.  

There exists a second tensor product of $\pi_{V}$ and $\pi_{W}$ that  has 
the same carrier space $V
\otimes W$. This will be called the  {\em{twisted tensor product}}, and is 
defined as the mapping
$\pi_{V} \twsqtimes \pi_{W}$ from $V \otimes W$ to $V
\otimes W \otimes \cal{A}$ that is given by
\begin{equation} (\pi_{V} \twsqtimes \pi_{W})(v \otimes w) = ((id \otimes 
id \otimes M) \circ (id
\otimes id \otimes \sigma) \circ (id
\otimes \sigma \otimes id) \circ (\pi_{V} \otimes \pi_{W}))(v \otimes w)
\label{eq:R164n} \end{equation} for all $v \in V$ and all $w \in W$. It is 
again easily checked that
the conditions (\ref{eq:R1}) and (\ref{eq:R2}) are satisfied with $\pi_{V}$ 
replaced by $\pi_{V}
\twsqtimes \pi_{W}$ and $V$ replaced by $V \otimes W$, so $\pi_{V} 
\twsqtimes \pi_{W}$ is also a
coaction of $\cal{A}$ with carrier space $V \otimes W$.  Then 
(\ref{eq:R5}) and (\ref{eq:R164n}) give
\begin{equation} (\pi_{V} \twsqtimes \pi_{W})(v_{j} \otimes w_{k}) = 
\sum_{s = 1}^{d_{V}}\sum_{t =
1}^{d_{W}} (v_{s} \otimes w_{t}) \otimes (M(\pi_{tk}^{W} \otimes 
\pi_{sj}^{V})) ,
\label{eq:R164r} \end{equation} which implies that the matrix coefficients 
of $\pi_{V} \twsqtimes
\pi_{W}$ are given by 
\begin{equation} (\pi_{V} \twsqtimes \pi_{W})_{st,jk} = M(\pi_{tk}^{W} 
\otimes
\pi_{sj}^{V}) .
\label{eq:R164s} \end{equation}
   
The tensor product $\pi_{W} \sqtimes \pi_{V}$, whose carrier space is
$\pi_{W} \otimes \pi_{V}$, has matrix coefficients that are given 
(according to (\ref{eq:R155})) by 
\begin{equation} (\pi_{W} \sqtimes \pi_{V})_{ts,kj} = M(\pi_{tk}^{W} 
\otimes
\pi_{sj}^{V}) .
\label{eq:R164t} \end{equation} As the matrix coefficients of 
(\ref{eq:R164s}) and (\ref{eq:R164t})
differ only in their ordering of the pairs of  indices that label their rows 
(and, in the
corresponding manner, their columns), it follows that $\pi_{V} \twsqtimes 
\pi_{W}$ and
$\pi_{W} \sqtimes \pi_{V}$ are {\em{equivalent}} corepresentations. If 
$\cal{A}$ is
{\em{coquasitriangular}} (that is, if  ${\cal{A}}^{\prime}$ is 
quasitriangular (c.f
Drinfel'd$^{30}$, then $\pi_{V} \sqtimes \pi_{W}$ and $\pi_{W} \sqtimes 
\pi_{V}$ are equivalent, so
in this case $\pi_{V} \sqtimes \pi_{W}$ and $\pi_{V} \twsqtimes 
\pi_{W}$ are equivalent. (Of course,
in the special case in which $\cal{A}$ is {\em{commutative}},  the 
corepresentations $\pi_{V}
\sqtimes \pi_{W}$ and $\pi_{V} \twsqtimes \pi_{W}$ are 
{\em{identical}}).

Applying these considerations to the special case in which $\pi_{V} = 
\pi^{p}$ and
$\pi_{W} = \pi^{q}$ are two irreducible unitary corepresentations of 
$\cal{A}$,  (\ref{eq:R155}) and
(\ref{eq:R164s}) become
\begin{equation} (\pi^{p} \sqtimes \pi^{q})_{st,jk} = M(\pi_{sj}^{p} 
\otimes \pi_{tk}^{q}) 
\label{eq:R155a} \end{equation} and
\begin{equation} (\pi^{p} \twsqtimes \pi^{q})_{st,jk} = M(\pi_{tk}^{q} 
\otimes
\pi_{sj}^{p}) 
\label{eq:R164sa} \end{equation} respectively.

\subsection{Characters}

The {\em{character}} of a corepresentation $\pi_{V}$ of $\cal{A}$ of 
dimension $d_{V}$ is defined in
terms of its matrix coefficients by
\begin{equation} \chi^{V} = \sum_{j=1}^{d_{V}} \pi^{V}_{jj} ,
\label{eq:R164a} \end{equation} so $\chi^{V}$ is also an element of 
$\cal{A}$. Clearly equivalent
corepresentations have identical characters.

If $\pi^{p}$ and $\pi^{q}$ are two irreducible unitary corepresentations of 
$\cal{A}$ (assumed to be
inequivalent if $p \neq q$) and if $\chi^{p}$ and $\chi^{q}$ are their 
corresponding characters,  
then (\ref{eq:D.eq2.1.15}) and (\ref{eq:D.eq2.1.16}) imply that
\begin{equation} h(M(\chi^{p*} \otimes \chi^{q})) = h(M(\chi^{q} \otimes 
\chi^{p*})) =
\delta^{pq} .
\label{eq:R164b} \end{equation}

If $\pi_{V}$ is a (completely) reducible corepresentation of $\cal{A}$ that 
is equivalent to the
direct sum $\sum \oplus \;n^{p}\pi^{p}$, then the number of times 
$n^{p}$ that the irreducible
corepresentation $\pi^{p}$ (or a corepresentation equivalent to $\pi^{p}$) 
appears in the reduction
of  $\pi_{V}$ is given by
\begin{equation} \chi^{V}  = \sum_{p} n^{p}\chi^{p} ,
\label{eq:R164c} \end{equation} so, by (\ref{eq:R164b}),
\begin{equation} n^{p} = h(M(\chi^{V} \otimes \chi^{p*})) = 
h(M(\chi^{p*} \otimes
\chi^{V})) .
\label{eq:R164d} \end{equation}
  
For the tensor product $\pi^{p} \sqtimes \pi^{q}$ of two irreducible 
unitary corepresentations
$\pi^{p}$ and $\pi^{q}$ the character is just $\chi^{p}\chi^{q}\; ( = 
M(\chi^{p} \otimes
\chi^{q}))$, so the number of times $n_{pq}^{r}$ that the irreducible 
corepresentation $\pi^{r}$ (or
a corepresentation equivalent to $\pi^{r}$) appears in the reduction of  
$\pi^{p} \sqtimes \pi^{q}$
is given by
\begin{equation} n_{pq}^{r} = h(\chi^{p}\chi^{q}\chi^{r*}) = 
h(\chi^{r*}\chi^{p}\chi^{q}) .
\label{eq:R164f} \end{equation}

If $\bar{\pi}^{p}$ is the irreducible unitary corepresentation that is 
{\em{conjugate}} to
$\pi^{p}$, then (\ref{eq:extra18}) and (\ref{eq:D.lemma2.1.10}) show that
\begin{equation} n_{pq}^{r} = n_{\bar{p}r}^{q}\; ,\; n_{r\bar{p}}^{q} = 
n_{qp}^{r} ,
\label{eq:R164h,R164i} \end{equation} where $n_{\bar{p}r}^{q}$ is the 
number of times that $\pi^{q}$
(or its equivalent) appears  in the reduction of  $\bar{\pi}^{p} \sqtimes 
\pi^{r}$, and
$n_{r\bar{p}}^{q}$ is the number of times that
$\pi^{q}$ (or its equivalent) appears  in the reduction of  $\pi^{r} 
\sqtimes \bar{\pi}^{p}$. 

\subsection{Clebsch-Gordan coefficients}

As above, suppose that the the tensor product $\pi^{p} \sqtimes \pi^{q}$ 
of two irreducible unitary
corepresentations
$\pi^{p}$ and $\pi^{q}$ is reducible, and that $n_{pq}^{r}$ is the number 
of times that the
irreducible corepresentation
$\pi^{r}$ (or a corepresentation equivalent to it) appears in the reduction 
of  $\pi^{p} \sqtimes
\pi^{q}$. If $\pi^{p}$ has carrier space $V^{p}$ with basis elements
$v^{p}_{1},v^{p}_{1}, \ldots, v^{p}_{d_{p}}$ and $\pi^{q}$ has carrier space 
$V^{q}$ with basis
elements
$v^{q}_{1},v^{q}_{1}, \ldots, v^{q}_{d_{q}}$, then the set of elements 
$v^{p}_{j} \otimes v^{q}_{k}$
form a basis for $V^{p} \otimes V^{q}$, the carrier space of $\pi^{p} 
\sqtimes \pi^{q}$, and
consequently appropriate linear combinations form bases for all the 
irreducible corepresentations
$\pi^{r}$ that appear in the reduction of the tensor product. Let
$w^{r,\alpha}_{\ell}$ be such a combination, so that 
\begin{equation} w^{r,\alpha}_{\ell} = \sum_{j=1}^{d_{p}} 
\sum_{k=1}^{d_{q}} \left( \begin{array}{cc}
p & q\\ j &  k \end{array} \right| \left. \begin{array}{ccc}r&,& 
\alpha\\\ell & &  
\end{array}\right) v^{p}_{j} \otimes v^{q}_{k} ,
\label{eq:C43} \end{equation} for $\ell = 1, 2, \ldots, d_{r}$, and $\alpha 
= 1, 2, \ldots,
n_{pq}^{r}$, and 
\begin{equation} (\pi^{p} \sqtimes \pi^{q})(w^{r,\alpha}_{\ell})  = 
\sum_{u=1}^{d_{r}}
w^{r,\alpha}_{u} \otimes \pi^{r}_{ul} ,
\label{eq:C44} \end{equation} for $u = 1, 2, \ldots, d_{r}$, and $\alpha = 
1, 2, \ldots,
n_{pq}^{r}$. The inverse of (\ref{eq:C43}) is
\begin{equation}  v^{p}_{j} \otimes v^{q}_{k} = \sum_{r} 
\sum_{\alpha=1}^{n_{pq}^{r}} \sum_{\ell
=1}^{d_{r}} \left(
\begin{array}{ccc} r&,& \alpha\\\ell & &     \end{array} \right| \left. 
\begin{array}{cc}p & q\\ j &
k   \end{array}\right) w^{r,\alpha}_{\ell} ,
\label{eq:C45} \end{equation} for $j = 1, 2, \ldots, d_{p}$ and $k = 1, 2, 
\ldots, d_{q}$.  The
{\em{Clebsch-Gordan coefficients}} defined in (\ref{eq:C43}) form the 
elements of a $d_{p} \times
d_{q}$ matrix ${\bf{C}}$, while the inverse coefficients defined in 
(\ref{eq:C45}) form the elements
of 
${\bf{C}}^{-1}$, where
\begin{equation}{\bf{C}}^{-1} ({\bf{\pi}}^{p} \sqtimes 
{\bf{\pi}}^{q}){\bf{C}} = \sum_{r} \oplus
n_{pq}^{r} {\bf{\pi}}^{r} .
\label{eq:R159} \end{equation} That is,
\begin{equation} \sum_{j=1}^{d_{p}} \sum_{t=1}^{d_{q}} (\pi^{p} \sqtimes 
\pi^{q})_{is,jt} \left(
\begin{array}{cc} p & q\\ j &  t \end{array} \right| \left. 
\begin{array}{ccc}r&,& \alpha\\\ell & &  
\end{array}\right) =  \sum_{u=1}^{d_{r}} \left(
\begin{array}{cc} p & q\\ i &  s \end{array} \right| \left. 
\begin{array}{ccc}r&,& \alpha\\u & &  
\end{array}\right) \pi^{r}_{u\ell}
\label{eq:C45q} \end{equation} for $i = 1, 2, \ldots, d_{p}$, $s = 1, 2, 
\ldots, d_{q}$, $\ell = 1,
2, \ldots, d_{r}$, and $\alpha = 1, 2, \ldots, n_{pq}^{r}$.

\subsection{Products of basis functions and Clebsch-Gordan coefficients}

Consider first the {\em{right}} regular corepresentation. If $\phi^{pR}_{j}$ 
and $\psi^{qR}_{k}$ are
basis functions of the unitary irreducible corepresentations $\pi^{p}$ and 
$\pi^{q}$ of $\cal{A}$,
then (\ref{eq:R164v}) and (\ref{eq:R155a}) imply that 
\begin{equation} \pi^{R}_{\cal{A}}(M(\phi^{pR}_{j} \otimes 
\psi^{qR}_{k})) = \sum_{s=1}^{d_{p}} 
\sum_{t=1}^{d_{q}} M(\phi^{pR}_{s} \otimes \psi^{qR}_{t}) \otimes 
(\pi^{p} \sqtimes
\pi^{q})_{st,jk}\; ,
\label{eq:R154} \end{equation} so that if the set of products 
$M(\phi^{pR}_{j} \otimes
\psi^{qR}_{k})$ (for $j = 1, 2, \ldots, d_{p}$, and $k = 1, 2, \ldots, d_{q}$) 
form a linearly
independent set, then they form a basis for the tensor product 
corepresentation $\pi^{p} \sqtimes
\pi^{q}$. Thus, by (\ref{eq:C43}) and (\ref{eq:C44}), there exists a set of 
elements
$\theta^{r,\alpha R}_{\ell}$, all members of $\cal{A}$, that are defined 
by
\begin{equation} \theta^{r,\alpha R}_{\ell} = \sum_{j=1}^{d_{p}} 
\sum_{k=1}^{d_{q}}
\left(
\begin{array}{cc} p & q\\ j &  k \end{array} \right| \left. 
\begin{array}{ccc}r&,& \alpha\\\ell & &  
\end{array}\right) M(\phi^{pR}_{j} \otimes \psi^{qR}_{k}) ,
\label{eq:R156} \end{equation} for $\ell = 1, 2, \ldots, d_{r}$, and $\alpha 
= 1, 2, \ldots,
n_{pq}^{r}$, and which have the property that 
\begin{equation}
\pi^{R}_{\cal{A}}(\theta^{r,\alpha R}_{\ell})  =
\sum_{u=1}^{d_{r}}
\theta^{r,\alpha R}_{u} \otimes \pi^{r}_{ul} ,
\label{eq:R157} \end{equation} for $u = 1, 2, \ldots, d_{r}$, and $\alpha = 
1, 2, \ldots,
n_{pq}^{r}$. By (\ref{eq:C45}), the inverse of (\ref{eq:R156}) is then
\begin{equation}  M(\phi^{pR}_{j} \otimes \psi^{qR}_{k}) = \sum_{r} 
\sum_{\alpha=1}^{n_{pq}^{r}}
\sum_{\ell =1}^{d_{r}} \left(
\begin{array}{ccc} r&,& \alpha\\\ell & &     \end{array} \right| \left. 
\begin{array}{cc}p & q\\ j &
k   \end{array}\right) \theta^{r,\alpha R}_{\ell} ,
\label{eq:R158} \end{equation} for $j = 1, 2, \ldots, d_{p}$ and $k = 1, 2, 
\ldots, d_{q}$.  On
applying the projection operator
${\cal{P}}^{rR}_{ul}$ of (\ref{eq:R143a}) to $M(\phi^{pR}_{j}
\otimes
\psi^{qR}_{k})$, and using (\ref{eq:R153}), (\ref{eq:R158}) and 
(\ref{eq:R156}), it follows that
\begin{eqnarray} \begin{array}{lr}{\cal{P}}^{rR}_{ul}(M(\phi^{pR}_{j} 
\otimes \psi^{qR}_{k})) & 
\\ =  \sum_{\alpha=1}^{n_{pq}^{r}} \sum_{v=1}^{d_{r}} 
\sum_{s=1}^{d_{p}}
\sum_{t=1}^{d_{q}} \left(\begin{array}{ccc} r&,&
\alpha\\\ell & &    
\end{array} \right| \left. \begin{array}{cc}p & q\\ j & k   
\end{array}\right) \left(
\begin{array}{cc} p & q\\ s &  t \end{array} \right| \left. 
\begin{array}{ccc}r&,& \alpha\\v & &  
\end{array}\right) &  \\
 \times d_{r} \{(({\bf{F}}^{r})^{-1})_{vu}/tr(({\bf{F}}^{r})^{-1})\}
\; M(\phi^{pR}_{s}
\otimes \psi^{qR}_{t})  & \end{array}
\label{eq:R161} \end{eqnarray} for $j = 1, 2, \ldots, d_{p}$ and $k = 1, 2, 
\ldots, d_{q}$. However,
the definition (\ref{eq:R143a}) taken in conjunction with  (\ref{eq:R154}) 
and (\ref{eq:R155a}) gives
\begin{equation} {\cal{P}}^{rR}_{ul}(M(\phi^{pR}_{j} \otimes 
\psi^{qR}_{k})) =
d_{r}\sum_{s=1}^{d_{p}}
\sum_{t=1}^{d_{q}} M(\phi^{pR}_{s} \otimes
\psi^{qR}_{t})\; h(\pi^{r*}_{ul}\pi^{p}_{sj}\pi^{q}_{tk}) ,
\label{eq:R162} \end{equation} so equating coefficients of 
$M(\phi^{pR}_{s} \otimes
\psi^{qR}_{t})$ in  (\ref{eq:R161}) and (\ref{eq:R162}) yields
\begin{eqnarray} \begin{array}{lcl}h(\pi^{r*}_{ul}\pi^{p}_{sj}\pi^{q}_{tk})
 &= & \sum_{\alpha=1}^{n_{pq}^{r}} \sum_{v=1}^{d_{r}} 
\left(\begin{array}{ccc} r&,&
\alpha\\\ell & &    
\end{array} \right| \left. \begin{array}{cc}p & q\\ j & k   
\end{array}\right) \left(
\begin{array}{cc} p & q\\ s &  t \end{array} \right| \left. 
\begin{array}{ccc}r&,& \alpha\\v & &  
\end{array}\right)   \\
  & &\times  \{(({\bf{F}}^{r})^{-1})_{vu}/tr(({\bf{F}}^{r})^{-1})\}\end{array}
\label{eq:R163} \end{eqnarray} for all $j = 1, 2, \ldots, d_{p}$, $k = 1, 2, 
\ldots, d_{q}$, and $l
= 1, 2, \ldots, d_{r}$.   Of course, this implies that
\begin{eqnarray} \begin{array}{lcl}h(\pi^{r*}_{ul}\pi^{q}_{tk}\pi^{p}_{sj})
 &= & \sum_{\alpha=1}^{n_{qp}^{r}} \sum_{v=1}^{d_{r}} 
\left(\begin{array}{ccc} r&,&
\alpha\\\ell & &    
\end{array} \right| \left. \begin{array}{cc}q & p\\ k & j   
\end{array}\right) \left(
\begin{array}{cc} q & p\\ t &  s \end{array} \right| \left. 
\begin{array}{ccc}r&,& \alpha\\v & &  
\end{array}\right)   \\
  & &\times  \{(({\bf{F}}^{r})^{-1})_{vu}/tr(({\bf{F}}^{r})^{-1})\}\end{array}
\label{eq:R164} \end{eqnarray} for all $j = 1, 2, \ldots, d_{p}$, $k = 1, 2, 
\ldots, d_{q}$, and $l
= 1, 2, \ldots, d_{r}$.

Although the conclusions (\ref{eq:R163}) and (\ref{eq:R164}) also follow 
from consideration of the
{\em{left}} regular coaction, some of the intermediate results are 
significantly different in this
case. Firstly (\ref{eq:L73}) implies that 
\begin{equation} \pi^{L}_{\cal{A}}(M(\phi^{pL}_{j} \otimes 
\psi^{qL}_{k})) = \sum_{s=1}^{d_{p}} 
\sum_{t=1}^{d_{q}} M(\phi^{pL}_{s} \otimes \psi^{qL}_{t}) \otimes 
M(\pi^{q}_{tk}
\otimes \pi^{p}_{sj}) ,
\label{eq:L77} \end{equation} so, by (\ref{eq:R164sa}),
\begin{equation} \pi^{L}_{\cal{A}}(M(\phi^{pL}_{j} \otimes 
\psi^{qL}_{k})) =
\sum_{s=1}^{d_{p}} 
\sum_{t=1}^{d_{q}} M(\phi^{pL}_{s} \otimes \psi^{qL}_{t}) \otimes 
(\pi^{p} \twsqtimes
\pi^{q})_{st,jk}\; ,
\label{eq:L78} \end{equation} so that if the set of products 
$M(\phi^{pL}_{j} \otimes
\psi^{qL}_{k})$ (for $j = 1, 2, \ldots, d_{p}$, and $k = 1, 2, \ldots, d_{q}$) 
form a linearly
independent set, then they form a basis for the {\em{twisted}} tensor 
product corepresentation
$\pi^{p} \twsqtimes
\pi^{q}$. However, on writing $\Phi^{qp}_{kj} = M(\phi^{pL}_{j} \otimes 
\psi^{qL}_{k})$,
(\ref{eq:L78}) can be re-expressed as
\begin{equation} \pi^{L}_{\cal{A}}(\Phi^{qp}_{kj}) = \sum_{s=1}^{d_{p}} 
\sum_{t=1}^{d_{q}} \Phi^{qp}_{ts} \otimes (\pi^{q} \sqtimes 
\pi^{p})_{ts,kj} \; .
\label{eq:L79} \end{equation} That is, the set $\Phi^{qp}_{kj}$ (for $k = 1, 
2, \ldots, d_{q}$ and
$j = 1, 2, \ldots, d_{p}$) form a basis for the {\em{ordinary}} tensor 
product corepresentation
$\pi^{q} \sqtimes
\pi^{p}$. Consequently, there exists a set of elements
$\theta^{r,\alpha L}_{\ell}$, all members of
$\cal{A}$, that are defined by
\begin{equation} \theta^{r,\alpha L}_{\ell} = \sum_{j=1}^{d_{p}} 
\sum_{k=1}^{d_{q}}
\left(
\begin{array}{cc} q & p\\ k &  j \end{array} \right| \left. 
\begin{array}{ccc}r&,& \alpha\\\ell & &  
\end{array}\right) \Phi^{qp}_{kj} ,
\label{eq:L80} \end{equation} for $\ell = 1, 2, \ldots, d_{r}$, and $\alpha 
= 1, 2, \ldots,
n_{qp}^{r}$, and which have the property that 
\begin{equation}
\pi^{R}_{\cal{A}}(\theta^{r,\alpha L}_{\ell})  =
\sum_{u=1}^{d_{r}}
\theta^{r,\alpha L}_{u} \otimes \pi^{r}_{ul} ,
\label{eq:L83} \end{equation} for $u = 1, 2, \ldots, d_{r}$, and $\alpha = 
1, 2, \ldots,
n_{qp}^{r}$. Thus
\begin{equation} \theta^{r,\alpha L}_{\ell} = \sum_{j=1}^{d_{p}} 
\sum_{k=1}^{d_{q}}
\left(
\begin{array}{cc} q & p\\ k &  j \end{array} \right| \left. 
\begin{array}{ccc}r&,& \alpha\\\ell & &  
\end{array}\right) M(\phi^{pL}_{j} \otimes \psi^{qL}_{k}) ,
\label{eq:L81} \end{equation} whose inverse is
\begin{equation}  M(\phi^{pL}_{j} \otimes \psi^{qL}_{k}) = \sum_{r} 
\sum_{\alpha=1}^{n_{qp}^{r}}
\sum_{\ell =1}^{d_{r}} \left(
\begin{array}{ccc} r&,& \alpha\\\ell & &     \end{array} \right| \left. 
\begin{array}{cc}q & p\\ k &
j   \end{array}\right) \theta^{r,\alpha L}_{\ell} ,
\label{eq:L82} \end{equation} for $j = 1, 2, \ldots, d_{p}$ and $k = 1, 2, 
\ldots, d_{q}$.  
Applying the projection operator ${\cal{P}}^{rL}_{ul}$ of (\ref{eq:R143a}) 
to $M(\phi^{pL}_{j}
\otimes
\psi^{qL}_{k})$, and using (\ref{eq:R153}), (\ref{eq:L82}) and 
(\ref{eq:L81}), it follows that
\begin{eqnarray} \begin{array}{lr}{\cal{P}}^{rL}_{ul}(M(\phi^{pL}_{j} 
\otimes \psi^{qL}_{k})) & 
\\ =  \sum_{\alpha=1}^{n_{qp}^{r}} \sum_{v=1}^{d_{r}} 
\sum_{s=1}^{d_{p}}
\sum_{t=1}^{d_{q}} \left(\begin{array}{ccc} r&,&
\alpha\\\ell & &    
\end{array} \right| \left. \begin{array}{cc}q & p\\ k & j   
\end{array}\right) \left(
\begin{array}{cc} q & p\\ t &  s \end{array} \right| \left. 
\begin{array}{ccc}r&,& \alpha\\v & &  
\end{array}\right) &  \\
 \times d_{r} \{(({\bf{F}}^{r})^{-1})_{vu}/tr(({\bf{F}}^{r})^{-1})\}
\; M(\phi^{pL}_{s}
\otimes \psi^{qL}_{t})  & \end{array}
\label{eq:L85} \end{eqnarray} for $j = 1, 2, \ldots, d_{p}$ and $k = 1, 2, 
\ldots, d_{q}$. However,
the definition (\ref{eq:R143a}) taken in conjunction with  (\ref{eq:L77}) 
gives
\begin{equation} {\cal{P}}^{rL}_{ul}(M(\phi^{pL}_{j} \otimes 
\psi^{qL}_{k})) = d_{r}
\sum_{s=1}^{d_{p}}
\sum_{t=1}^{d_{q}} M(\phi^{pL}_{s} \otimes
\psi^{qL}_{t})\; h(\pi^{r*}_{ul}\pi^{q}_{tk}\pi^{p}_{sj}) ,
\label{eq:L86} \end{equation} so equating coefficients of $M(\phi^{pL}_{s} 
\otimes
\psi^{qL}_{t})$ in  (\ref{eq:L85}) and (\ref{eq:L86}) yields (\ref{eq:R164}) 
again.

It is also possible to obtain (\ref{eq:R163}) (and hence (\ref{eq:R164})) 
{\em{without}} making any
linear independence assumptions about the set of products 
$M(\phi^{pR}_{j} \otimes
\psi^{qR}_{k})$, for by (\ref{eq:extra16}), (\ref{eq:C43}), (\ref{eq:C44}), 
and
(\ref{eq:C45}),
\begin{eqnarray*} \begin{array}{l}
 \sum_{r}\sum_{\alpha=1}^{n_{pq}^{r}} \sum_{s=1}^{d_{p}}
\sum_{t=1}^{d_{q}} \sum_{v,w=1}^{d_{r}} \left(\begin{array}{ccc} r&,&
\alpha\\w & &    
\end{array} \right| \left. \begin{array}{cc}p & q\\ j & k   
\end{array}\right) \left(
\begin{array}{cc} p & q\\ s &  t \end{array} \right| \left. 
\begin{array}{ccc}r&,& \alpha\\v & &  
\end{array}\right)   \;v_{s}^{p} \otimes v_{t}^{q} \otimes \pi_{vw}^{r}\\
= \sum_{s=1}^{d_{p}}
\sum_{t=1}^{d_{q}} v_{s}^{p} \otimes v_{t}^{q} \otimes M(\pi_{sj}^{p} 
\otimes
\pi_{tk}^{q}) \end{array} \end{eqnarray*} 
for all
$j = 1, 2,
\ldots, d_{p}$, and $k = 1, 2, \ldots, d_{q}$. As the set $v_{s}^{p} \otimes 
v_{t}^{q}$ (for $s = 1,
2,
\ldots, d_{p}$, and $t = 1, 2, \ldots, d_{q}$) is
certainly linearly independent, it follows that
\begin{eqnarray*} \begin{array}{c}
 \sum_{r}\sum_{\alpha=1}^{n_{pq}^{r}}  \sum_{v,w=1}^{d_{r}} 
\left(\begin{array}{ccc} r&,&
\alpha\\w & &    
\end{array} \right| \left. \begin{array}{cc}p & q\\ j & k   
\end{array}\right) \left(
\begin{array}{cc} p & q\\ s &  t \end{array} \right| \left. 
\begin{array}{ccc}r&,& \alpha\\v & &  
\end{array}\right)   \; \pi_{vw}^{r}
= M(\pi_{sj}^{p} \otimes
\pi_{tk}^{q}) \end{array} \end{eqnarray*}
(for $j,s = 1,
2,
\ldots, d_{p}$, and $k,t = 1, 2, \ldots, d_{q}$). On replacing $r$ by 
$r^{\prime}$ in the sums of
on the left-hand side, multiplying through by $\pi^{r*}_{u\ell}$ from the 
left, applying the Haar
functional $h$, and using (\ref{eq:D.prop2.1.3}) and (\ref{eq:D.eq2.1.16}), 
one regains
(\ref{eq:R163}).

\section{The irreducible tensor operators}

\subsection{Introduction}

Let $\pi^{q}$ be a unitary irreducible right coaction of $\cal{A}$ of 
dimension
$d_{q}$ with matrix coefficients $\pi^{q}_{jk}$. It will be shown that 
within both the right and the
left regular corepresentation formalisms there exist two types of 
irreducible tensor operators that
both belong to this corepresentation
$\pi^{q}$. These will be denoted by
$Q^{qX}_{j}$ and ${\widetilde{Q}}^{qX}_{j}$ (for $j = 1, 2, \ldots, d_{q}$ 
and for $X = R$ or $L$),
and will be called {\em{ordinary}} and {\em{twisted}} irreducible tensor 
operators respectively.
Naturally the two types of irreducible tensor operators coincide in the 
special case in which
$\cal{A}$ is commutative.  Moreover, it should be noted that if 
${\cal{A}}^{op}$ is the  Hopf
algebra in which the multiplication operator $M$ and antipode $S$ of 
$\cal{A}$ are replaced by $M
\circ \sigma$ and
$S^{-1}$ respectively, then the `twisted' irreducible
tensor operators for $\cal{A}$ are the `ordinary' irreducible
tensor operators for ${\cal{A}}^{op}$ and the `ordinary' irreducible tensor 
operators for
$\cal{A}$ are the `twisted' irreducible tensor operators for 
${\cal{A}}^{op}$.

\subsection{Definitions in the right regular corepresentation formalism}

\subsubsection{Definition of the ordinary irreducible tensor operators
$Q^{qR}_{j}$}

The {\em{ordinary irreducible tensor operators}}
$Q^{qR}_{j}$ belonging to the unitary irreducible right coaction $\pi^{q}$ 
of
$\cal{A}$ are {\em{defined}} to be members of $\cal{L}(\cal{A})$ that 
satisfy the condition 
\begin{equation}  ((id \otimes M) \circ (\Delta \otimes id) \circ 
(Q^{qR}_{j}
\otimes S) \circ \Delta)(a) = \sum_{k=1}^{d_{q}} Q^{qR}_{k}(a) \otimes
\pi^{q}_{kj}
\label{eq:R95} \end{equation}  for all $a \in \cal{A}$ and all $j = 1, 2, 
\ldots, d_{q}$. Hereafter
$\cal{L}(\cal{A})$ denotes the set of linear operators that map $\cal{A}$ 
into $\cal{A}$.

Clearly this definition involves {\em{only}} quantities defined on 
$\cal{A}$, for the right-hand side
is a member of
$\cal{A}\otimes \cal{A}$. By virtues of (\ref{eq:R8}) this definition can be 
written equivalently in
terms of the right regular coaction
$\pi^{R}_{\cal{A}}$ as
\begin{equation}  ((id \otimes M) \circ (\pi^{R}_{\cal{A}} \otimes id) \circ 
(Q^{qR}_{j}
\otimes S) \circ \pi^{R}_{\cal{A}})(a) = \sum_{k=1}^{d_{q}} Q^{qR}_{k}(a) 
\otimes
\pi^{q}_{kj}
\label{eq:L46} \end{equation}  for all $a \in \cal{A}$ and all $j = 1, 2, 
\ldots, d_{q}$.  (The
motivation behind the definition (\ref{eq:R95}) is explained in  Appendix 
B.2.a).

It will now be shown that (\ref{eq:R95}) provides a {\em{consistent}} 
definition, in that it can be
re-expressed by saying that the operators $Q^{qR}_{j}$ (for $j = 1, 2, 
\ldots, d_{q}$) form the
basis of a carrier space for a certain right coaction of $\cal{A}$. The 
motivation for this
definition is given in  Appendix B.2.b., where the special case in which 
$\cal{A}$ is
finite-dimensional is considered in detail, but for general case it is 
necessary to apply
additional conditions to the domain of this coaction.  The analysis of 
Appendix B.2.b implies that if
$\cal{A}$ is {\em{finite}}-dimensional then,  for {\em{every}} operator $Q 
\in
{\cal{L}}(\cal{A})$, there exist operators $Q_{i} \in
{\cal{L}}(\cal{A})$ and elements $q_{i}$, both indexed by the same 
{\em{finite}} index set $I$,
such that 
\begin{equation}  \sum_{i \in I} Q_{i}(a) \otimes q_{i} = ((id \otimes M) 
\circ (\Delta \otimes id)
\circ (Q \otimes S) \circ \Delta)(a) \label{eq:newT}\end{equation} for all 
$a \in \cal{A}$, but
this is not necessarily true for every $Q \in
{\cal{L}}(\cal{A})$ if $\cal{A}$ is infinite-dimensional. However, one can 
always define a subspace 
${\cal{T}}(\cal{A})$ of ${\cal{L}}(\cal{A})$ by the requirement that $Q \in
{\cal{T}}(\cal{A})$ if (i) $Q \in
{\cal{L}}(\cal{A})$, (ii)  $Q$ satisfies (\ref{eq:newT}) with $I$ finite, and 
(iii) each $Q_{i}$
appearing on the left-hand side of (\ref{eq:newT}) also satisfies a condition 
of the same form. This subspace ${\cal{T}}(\cal{A})$ is certainly not 
empty, for the identity operator belongs
to it (c.f. (\ref{eq:R117})), as does every irreducible tensor operator 
$Q^{qR}_{j}$
(c.f.(\ref{eq:R95})), and, as just noted, if $\cal{A}$ is finite-dimensional 
then
${\cal{T}}(\cal{A})$ is identical to ${\cal{L}}(\cal{A})$. The
{\em{definition}} of the required right coaction, which will be denoted
by $\pi^{R}_{{\cal{T}}(\cal{A})}$, is then that it is the mapping of 
${\cal{T}}({\cal{A}})$ into
${\cal{T}}({\cal{A}}) \otimes \cal{A}$ that given by
\begin{equation} \pi^{R}_{{\cal{T}}(\cal{A})}(Q) = \sum_{[Q]} Q_{[1]} 
\otimes
Q_{[2]}\;,\label{eq:R115a} \end{equation} where  $Q_{[1]} \in 
{\cal{T}}(\cal{A})$ and $Q_{[1]} \in
\cal{A}$ are such that
\begin{equation}  \sum_{[Q]} Q_{[1]}(a) \otimes Q_{[2]} = ((id \otimes M) 
\circ (\Delta \otimes id)
\circ (Q \otimes S) \circ \Delta)(a) \label{eq:R115b}\end{equation} for all 
$Q \in
{\cal{T}}({\cal{A}})$ and all $a \in \cal{A}$. It is then quite easily shown 
(using the identities (\ref{eq:QD8}) to
(\ref{eq:QD27d})) that
$\pi^{R}_{{\cal{T}}(\cal{A})}$ satisfies (\ref{eq:R1}) and (\ref{eq:R2}) 
(with $\pi_{V}$ and $V$
replaced by
$\pi^{R}_{{\cal{T}}(\cal{A})}$ and 
${\cal{T}}(\cal{A})$ respectively, and with all operators of 
${\cal{T}}(\cal{A})$ acting on any
member of $\cal{A}$ according to the prescription of (\ref{eq:R115b})). 
That is,  
\begin{equation} \sum_{[Q]} \sum_{[Q_{[1]}]} (Q_{[1]})_{[1]}(a) \otimes 
(Q_{[1]})_{[2]} \otimes
Q_{[2]} = \sum_{[Q]} Q_{[1]}(a) \otimes \Delta(Q_{[2]})\end{equation} and
\[ \sum_{[Q]}  (Q_{[1]})(a) \;\epsilon(Q_{[2]}) = Q(a) \] for all $Q \in
{\cal{T}}({\cal{A}})$ and all $a \in \cal{A}$. Hence
$\pi^{R}_{{\cal{T}}(\cal{A})}$ is indeed a right coaction with carrier space
${\cal{T}}({\cal{A}})$. Thus (\ref{eq:R115a}) and (\ref{eq:R115b}) imply 
that the definition
(\ref{eq:R95}) can be written equivalently as
\begin{equation} \pi^{R}_{{\cal{T}}(\cal{A})}(Q^{qR}_{j}) = 
\sum_{k=1}^{d_{q}} Q^{qR}_{k} \otimes
\pi^{q}_{kj}  
\label{eq:R115} \end{equation}  (for all $j = 1, 2, \ldots $). Because 
(\ref{eq:R115}) is similar in
form to (\ref{eq:R5}), and as
$\pi^{R}_{{\cal{T}}(\cal{A})}$ is a right coaction with carrier space 
${\cal{T}}({\cal{A}})$, the
consistency of the definition (\ref{eq:R95}) is now ensured.

\subsubsection{Definition of the twisted irreducible tensor operators
${\widetilde{Q}}^{qR}_{j}$}

The {\em{twisted irreducible tensor operators}} ${\widetilde{Q}}^{qR}_{j}$ 
belonging to the unitary
irreducible right coaction $\pi^{q}$ of $\cal{A}$ are {\em{defined}} to be 
members of
${\cal{L}}({\cal{A}})$ that satisfy the condition 
\begin{equation}  ((id \otimes M) \circ (id \otimes \sigma) \circ (\Delta 
\otimes id) \circ
({\widetilde{Q}}^{qR}_{j}
\otimes S^{-1}) \circ \Delta)(a) = \sum_{k=1}^{d_{q}} 
{\widetilde{Q}}^{qR}_{k}(a) \otimes
\pi^{q}_{kj}
\label{eq:R132} \end{equation}  for all $a \in \cal{A}$ and all $j = 1, 2, 
\ldots, d_{q}$. This can
be written equivalently in terms of the right regular coaction 
$\pi^{R}_{\cal{A}}$ as
\begin{equation}  ((id \otimes M) \circ (id \otimes \sigma) \circ 
(\pi^{R}_{\cal{A}} \otimes id)
\circ ({\widetilde{Q}}^{qR}_{j}
\otimes S^{-1}) \circ \pi^{R}_{\cal{A}})(a) = \sum_{k=1}^{d_{q}} 
{\widetilde{Q}}^{qR}_{k}(a) \otimes
\pi^{q}_{kj}
\label{eq:extra21} \end{equation}  for all $a \in \cal{A}$ and all $j = 1, 2, 
\ldots, d_{q}$.  Both
(\ref{eq:R132}) and (\ref{eq:extra21}) differ from the corresponding 
definitions (\ref{eq:R95}) and
(\ref{eq:L46}) only in the replacement of $M$ by $M \circ \sigma$ 
{\em{and}} $S$ by $S^{-1}$
(neither of which have any effect in the special case in which $\cal{A}$ is 
commutative). (See
 Appendix B.2.a for further discussion of this
pair of substitutions. It should be recorded that Rittenberg and 
Scheunert$^{8}$ noted previously,
in the context of what was essentially the `abstract carrier space 
formalism' (form (3) of
Section~I) as generalized to irreducible {\em{representations}}  of
${\cal{A}}^{\prime}$,  that these substitutions do produce another type of 
irreducible tensor
operator, but they did not pursue  this observation at all.)  

The demonstration that (\ref{eq:R132}) provides a {\em{consistent}} 
definition again involves
showing  that it can be re-expressed by saying that the operators 
${\widetilde{Q}}^{qR}_{j}$ (for $j
= 1, 2,
\ldots, d_{q}$) form the basis of a carrier space for a  right coaction
 of
$\cal{A}$. This right coaction $\widetilde{\pi}^{R}_{{\cal{T}}(\cal{A})}$ 
(and its associated space
${\cal{T}}(\cal{A})$) are essentially obtained from
$\pi^{R}_{{\cal{T}}(\cal{A})}$ by replacing $M$ by $M \circ \sigma$ and 
$S$ by
$S^{-1}$, so $\widetilde{\pi}^{R}_{{\cal{T}}(\cal{A})}$ is {\em{defined}} as 
the mapping of 
${\cal{T}}({\cal{A}})$ into
${\cal{T}}({\cal{A}}) \otimes \cal{A}$ that given by
\begin{equation} \widetilde{\pi}^{R}_{{\cal{T}}(\cal{A})}(Q) = \sum_{[Q]} 
Q_{[1]} \otimes
Q_{[2]}\;,\label{eq:R115amod} \end{equation}
where  $Q_{[1]} \in {\cal{T}}(\cal{A})$ and $Q_{[1]} \in \cal{A}$ are such 
that
\begin{equation}  \sum_{[Q]} Q_{[1]}(a) \otimes Q_{[2]} = ((id \otimes (M 
\circ \sigma)) \circ
(\Delta \otimes id) \circ (Q \otimes S^{-1}) \circ \Delta)(a) 
\label{eq:R115bmod}\end{equation}
for all $Q \in {\cal{T}}({\cal{A}})$ and all $a \in \cal{A}$. Then
(\ref{eq:R115amod}) and (\ref{eq:R115bmod}) imply that the definition 
(\ref{eq:R132}) can be written
equivalently as
\begin{equation} 
\widetilde{\pi}^{R}_{{\cal{T}}(\cal{A})}({\widetilde{Q}}^{qR}_{j}) =
\sum_{k=1}^{d_{q}} {\widetilde{Q}}^{qR}_{k} \otimes
\pi^{q}_{kj}  
\label{eq:R133} \end{equation}  (for all $j = 1, 2, \ldots $), which then 
ensures the consistency of
the definition (\ref{eq:R132}).

\subsection{Definitions in the left regular corepresentation formalism}

\subsubsection{Definition of the ordinary irreducible tensor operators
$Q^{qL}_{j}$}

The {\em{ordinary irreducible tensor operators}} $Q^{qL}_{j}$ belonging to 
the unitary irreducible
right coaction $\pi^{q}$ of $\cal{A}$ are {\em{defined}} to be members of 
${\cal{L}}({\cal{A}})$
that satisfy the condition
\begin{equation}  ((id \otimes M) \circ (\pi^{L}_{\cal{A}} \otimes id) \circ 
(Q^{qL}_{j}
\otimes S) \circ \pi^{L}_{\cal{A}})(a) = \sum_{k=1}^{d_{q}} Q^{qL}_{k}(a) 
\otimes
\pi^{q}_{kj}
\label{eq:L45} \end{equation}  for all $a \in \cal{A}$ and all $j = 1, 2, 
\ldots, d_{q}$. (This is
can be obtained from (\ref{eq:L46}) by replacing $X = R$ by $X = L$, the 
justification being
discussed in more detail in  Appendix B.3.a). 
In terms of the elementary operations of $\cal{A}$, (\ref{eq:L45}) can be 
re-expressed using
(\ref{eq:L1}) as 
\begin{equation}  (\sigma \circ (S \otimes id) \circ (M \otimes id) \circ 
(id \otimes \Delta) \circ
(S \otimes Q^{qL}_{j}) \circ \Delta)(a)  =
\sum_{k=1}^{d_{q}} Q^{qL}_{k}(a)
\otimes
\pi^{q}_{kj}
\label{eq:L36} \end{equation}  for all $a \in \cal{A}$ and all $j = 1, 2, 
\ldots, d_{q}$. 

The demonstration that (\ref{eq:L45})  provides a {\em{consistent}} 
definition proceeds in the same
way as above. In this case the appropriate right coaction will be denoted 
by
$\pi^{L}_{{\cal{T}}(\cal{A})}$, and will be {\em{defined}} as the mapping 
of 
${\cal{T}}({\cal{A}})$ into
${\cal{T}}({\cal{A}}) \otimes \cal{A}$ that given by
\begin{equation} \pi^{L}_{{\cal{T}}(\cal{A})}(Q) = \sum_{[Q]} Q_{[1]} 
\otimes
Q_{[2]}\;, \label{eq:R115aL} \end{equation}
where  $Q_{[1]} \in {\cal{T}}(\cal{A})$ and $Q_{[1]} \in \cal{A}$ are such 
that
\begin{equation}  \sum_{[Q]} Q_{[1]}(a) \otimes Q_{[2]} = (\sigma \circ (S
\otimes id)
\circ (M \otimes id) \circ (id \otimes \Delta) \circ (S \otimes Q) \circ 
\Delta)(a)
\label{eq:R115bL} \end{equation} for all $Q \in {\cal{T}}({\cal{A}})$ and 
all $a \in \cal{A}$.
(Here the subspace ${\cal{T}}({\cal{A}})$ of ${\cal{L}}({\cal{A}})$ is defined 
as in
Section~VI.B.1, but with the right-hand side of (\ref{eq:R115bL}) 
replacing the right-hand side of
(\ref{eq:R115b}) in (\ref{eq:newT}). The motivation for this definition is 
given in  Appendix
B.3.b).     Then the definition (\ref{eq:L45}) can be written equivalently as
\begin{equation} \pi^{L}_{{\cal{T}}(\cal{A})}(Q^{qL}_{j}) = 
\sum_{k=1}^{d_{q}} Q^{qL}_{k} \otimes
\pi^{q}_{kj}  
\label{eq:L41a} \end{equation}  (for all $j = 1, 2, \ldots $), which then 
guarantees  its
consistency.

\subsubsection{Definition of the twisted irreducible tensor operators
${\widetilde{Q}}^{qL}_{j}$}

The {\em{twisted irreducible tensor operators}} ${\widetilde{Q}}^{qL}_{j}$ 
belonging $\pi^{q}$ of
$\cal{A}$  are {\em{defined}} to be members of ${\cal{L}}({\cal{A}})$ that 
satisfy the condition 
\begin{equation}  ((id \otimes M) \circ (id \otimes \sigma) \circ 
(\pi^{L}_{\cal{A}} \otimes id)
\circ ({\widetilde{Q}}^{qL}_{j}
\otimes S^{-1}) \circ \pi^{L}_{\cal{A}})(a) = \sum_{k=1}^{d_{q}} 
{\widetilde{Q}}^{qL}_{k}(a) \otimes
\pi^{q}_{kj}
\label{eq:L63} \end{equation}  for all $a \in \cal{A}$ and all $j = 1, 2, 
\ldots, d_{q}$. In terms of
the elementary operations of $\cal{A}$, (\ref{eq:L63}) can be re-expressed 
using (\ref{eq:L1}) as 
\begin{equation}  ((id \otimes M) \circ (\sigma \otimes S) \circ (id 
\otimes \sigma) \circ (id
\otimes \Delta) \circ (id \otimes {\widetilde{Q}}^{qL}_{j}) \circ \Delta)(a)  
=
\sum_{k=1}^{d_{q}} {\widetilde{Q}}^{qL}_{k}(a) \otimes \pi^{q}_{kj}
\label{eq:L52} \end{equation}  for all $a \in \cal{A}$ and all $j = 1, 2, 
\ldots, d_{q}$. The
definition (\ref{eq:L63})  differs from the corresponding definition 
(\ref{eq:L45}) only in the
replacement of
$M$ by $M \circ \sigma$ {\em{and}} $S$ by $S^{-1}$ (neither of which 
have any effect in the special
case in which $\cal{A}$ is commutative). However, because the two $S$ 
factors in (\ref{eq:L36}) have
different origins, these substitutions do {\em{not}} convert (\ref{eq:L36}) 
into (\ref{eq:L52}).
(See  Appendix B.3.a for a further discussion
of this point).

The consistency of the definition (\ref{eq:L63})  is again shown in the 
same way as above. In this
case the appropriate right coaction will be denoted by
$\widetilde{\pi}^{L}_{{\cal{T}}(\cal{A})}$, and will be {\em{defined}} as 
the mapping of 
${\cal{T}}({\cal{A}})$ into
${\cal{T}}({\cal{A}}) \otimes \cal{A}$ that given by
\begin{equation} \widetilde{\pi}^{L}_{{\cal{T}}(\cal{A})}(Q) = \sum_{[Q]} 
Q_{[1]} \otimes
Q_{[2]}\;,\label{eq:R115aLmod} \end{equation}
where  $Q_{[1]} \in {\cal{T}}(\cal{A})$ and $Q_{[1]} \in \cal{A}$ are such 
that
\begin{equation}  \sum_{[Q]} Q_{[1]}(a) \otimes Q_{[2]} = ((id
\otimes M)
\circ (\sigma \otimes S) \circ (id \otimes \sigma) \circ (id
\otimes \Delta) \circ (id \otimes Q) \circ \Delta)(a) 
\label{eq:R115bLmod} \end{equation}
for all $Q \in {\cal{T}}({\cal{A}})$ and all $a \in \cal{A}$. (Here the 
subspace ${\cal{T}}({\cal{A}})$ of ${\cal{L}}({\cal{A}})$ is defined as in
Section~VI.B.1, but with the right-hand side of (\ref{eq:R115bLmod}) 
replacing the right-hand side of
(\ref{eq:R115b}) in (\ref{eq:newT}). The motivation for this
definition is given in 
Appendix B.3.b). Then  the definition (\ref{eq:L63}) can be written
equivalently as
\begin{equation} 
\widetilde{\pi}^{L}_{{\cal{T}}(\cal{A})}({\widetilde{Q}}^{qL}_{j}) =
\sum_{k=1}^{d_{q}} {\widetilde{Q}}^{qL}_{k} \otimes
\pi^{q}_{kj}  
\label{eq:L59} \end{equation}  (for all $j = 1, 2, \ldots $), which then 
implies the consistency of
(\ref{eq:L63}).

\subsection{The identity operator as an irreducible tensor operator}

Suppose that $Q$ is the {\em{identity operator}} $id$ of 
${\cal{L}}({\cal{A}})$ (so that $Q(a) = a$
for all $a
\in \cal{A}$). Then, on using (\ref{eq:QD11}) and (\ref{eq:QD26}), it 
follows that
\begin{equation} (id \otimes M) \circ (\Delta
\otimes id) \circ (id \otimes S) \circ \Delta = id \otimes 1_{{\cal{A}}} ,
\label{eq:R117} \end{equation}  which, by (\ref{eq:R95}), leads to the 
agreeable conclusion that the
identity operator
$id$ is an {\em{ordinary}} irreducible tensor operator in the {\em{right}} 
regular corepresentation
formalism for the one-dimensional {\em{identity}} corepresentation whose 
sole matrix coefficient is
$1_{{\cal{A}}}$.
 
It is easily checked (using (\ref{eq:R132}), (\ref{eq:L36}), and 
(\ref{eq:L52}) in place of
(\ref{eq:R95})), that
$id$ is also a {\em{twisted}} irreducible tensor operator for the identity 
corepresentation in the
{\em{right}} regular corepresentation formalism , as well as being a both 
an {\em{ordinary}} and a
{\em{twisted}} irreducible tensor operator  for
the identity corepresentation in the {\em{left}} regular corepresentation 
formalism.

In fact, if one were to adopt the view that these results concerning the 
identity operator are an
{\em{essential}} requirement of any sensible definition of irreducible 
tensor operators, the fact
that they are not true if  $M$ is replaced by $M \circ \sigma$  but  $S$ is 
left unchanged, nor if
$M$ is left unchanged but $S$ is replaced by $S^{-1}$, then precludes 
further consideration of these
possibilities.

\subsection{Products as irreducible tensor operators}

 Suppose now that $\psi^{qR}_{j}$ and $\psi^{qL}_{j}$ are sets of basis 
functions for $\pi^{q}$ (as
defined in (\ref{eq:R79}) and (\ref{eq:L20}) respectively) and that the 
operators $Q^{qX}_{j}$ and
${\widetilde{Q}}^{qX}_{j}$ are {\em{defined}} by
\begin{eqnarray} \left. \begin{array}{ccc}  Q^{qR}_{j}(a) & = & 
M(\psi^{qR}_{j} \otimes a) , \\
{\widetilde{Q}}^{qR}_{j}(a) &  = & M(a \otimes \psi^{qR}_{j}) ,\\ 
Q^{qL}_{j}(a) & = & M(a \otimes
\psi^{qL}_{j}) , \\ {\widetilde{Q}}^{qL}_{j}(a) & = & M(\psi^{qL}_{j} 
\otimes a) , \end{array}
\right\}
\label{eq:R118,R136,L44.L62} \end{eqnarray}   for all $a \in \cal{A}$. 
Then the identities
(\ref{eq:QD9}) to (\ref{eq:QD27e}) imply that the operators $Q^{qR}_{j}$,
${\widetilde{Q}}^{qR}_{j}$, $Q^{qL}_{j}$, and
${\widetilde{Q}}^{qL}_{j}$ do indeed satisfy (\ref{eq:R95}), (\ref{eq:R132}), 
(\ref{eq:L36}), and
(\ref{eq:L52}) respectively, and so are irreducible tensor operators  
belonging to
$\pi^{q}$.

\subsection{Two useful identities for the ordinary irreducible tensor 
operators $Q^{qX}_{j}$
and 
${\widetilde{Q}}^{qX}_{j}$}

 If  $Q^{qX}_{k}$ is an {\em{ordinary}} irreducible tensor operator 
belonging to the unitary
irreducible right coaction $\pi^{q}$ of $\cal{A}$ and  $\psi^{pX}_{j}$ is a 
set of basis functions
for the unitary irreducible right coaction $\pi^{p}$ of $\cal{A}$, then
\begin{equation} \pi^{X}_{{\cal{A}}}(Q^{qX}_{k}(\psi^{pX}_{j})) =
\sum_{s=1}^{d_{p}} \sum_{t=1}^{d_{q}} (Q^{qX}_{t}(\psi^{pX}_{s}))
\otimes (M(\pi^{q}_{tk} \otimes \pi^{p}_{sj})) ,
\label{eq:R119,L47} \end{equation} for $X = R$ and $L$, for all $j = 1,2, 
\ldots, d_{p}$, and $k =
1,2, \ldots, d_{q}$. That is, by (\ref{eq:R155a}),
\begin{equation} \pi^{X}_{{\cal{A}}}(Q^{qX}_{k}(\psi^{pX}_{j})) =
\sum_{s=1}^{d_{p}} \sum_{t=1}^{d_{q}} (Q^{qX}_{t}(\psi^{pX}_{s}))
\otimes (\pi^{q} \sqtimes \pi^{p})_{ts,kj}\; ,
\label{eq:extra22} \end{equation}  for $X = R$ and $L$, for all $j = 1,2, 
\ldots, d_{p}$ and $k =
1,2,
\ldots, d_{q}$.

 By contrast, if  ${\widetilde{Q}}^{qX}_{k}$ is a {\em{twisted}} irreducible 
tensor operator
belonging
$\pi^{q}$, then
\begin{equation} 
\pi^{X}_{{\cal{A}}}({\widetilde{Q}}^{qX}_{k}(\psi^{pX}_{j})) =
\sum_{s=1}^{d_{p}} \sum_{t=1}^{d_{q}} 
({\widetilde{Q}}^{qX}_{t}(\psi^{pX}_{s}))
\otimes (M(\pi^{p}_{sj} \otimes \pi^{q}_{tk} )) ,
\label{eq:R134,L65} \end{equation}  for $X =R$ and $L$, for all $j = 1,2, 
\ldots, d_{p}$ and $k =
1,2,
\ldots, d_{q}$. It should be noted that the factors in the second term of the 
right-hand side of
(\ref{eq:R134,L65}) are interchanged relative to those of 
(\ref{eq:R119,L47}), which implies, by
(\ref{eq:R164sa}), that
\begin{equation} 
\pi^{X}_{{\cal{A}}}({\widetilde{Q}}^{qX}_{k}(\psi^{pX}_{j})) =
\sum_{s=1}^{d_{p}} \sum_{t=1}^{d_{q}} 
({\widetilde{Q}}^{qX}_{t}(\psi^{pX}_{s}))
\otimes (\pi^{q} \twsqtimes \pi^{p})_{ts,kj} ,
\label{eq:extra24} \end{equation}  for $X =R$ and $L$, for all $j = 1,2, 
\ldots, d_{p}$, and $k =
1,2, \ldots, d_{q}$, which involves the {\em{twisted}} tensor product.

The proof of (\ref{eq:R119,L47}) is as follows. Applying (\ref{eq:R79}) or 
(\ref{eq:L20})   to the
case $a = \psi^{pX}_{u}$ and invoking (\ref{eq:L46}) or (\ref{eq:L45}) (as 
appropriate) gives
\begin{equation} \sum_{v=1}^{d_{p}} (id \otimes 
M)(\pi^{X}_{{\cal{A}}}(Q^{qX}_{j}(\psi^{pX}_{v}))
\otimes S(\pi^{p}_{vu})) = \sum_{k=1}^{d_{q}} (Q^{qX}_{k}(\psi^{pX}_{u})) 
\otimes \pi^{q}_{kj}  .
\label{eq:R120} \end{equation}  However, for any $a \in \cal{A}$ 
\[ (id \otimes M)( \pi^{X}_{{\cal{A}}}(a) \otimes S(\pi^{p}_{vu})) = 
\sum_{[a]} a_{[1]} \otimes
M(a_{[2]} \otimes S(\pi^{p}_{vu})) ,\]   the right-hand side of which, on 
multiplication from the
right with $1_{{\cal{A}}} \otimes
\pi^{p}_{ui}$, and summing over
$u$, and using (\ref{eq:D.sect1.1.45ii}) and (\ref{eq:D.sect1.1.45iii}), 
reduces to $\delta_{iv}
\pi^{X}_{{\cal{A}}}(a)$. The desired result (\ref{eq:R119,L47}) is then 
obtained by  multiplying both
sides of (\ref{eq:R120}) from the right with $1_{{\cal{A}}} \otimes 
\pi^{p}_{ui}$ and summing over
$u$. The line of proof for (\ref{eq:R134,L65}) is similar.

\subsection{Products of operators}

If $Q$ and $Q^{\prime}$ are any two members of ${\cal{T}}({\cal{A}})$ 
and
$\pi^{X}_{{\cal{T}}(\cal{A})}$ is the {\em{ordinary}} right coaction defined 
in (\ref{eq:R115a})
and (\ref{eq:R115b}) (for
$X = R$) and (\ref{eq:R115aL}) and (\ref{eq:R115bL}) (for $X = L$), then
\begin{equation}  \pi^{X}_{{\cal{T}}(\cal{A})}(QQ^{\prime}) =
\pi^{X}_{{\cal{T}}(\cal{A})}(Q)\;\pi^{X}_{{\cal{T}}(\cal{A})}(Q^{\prime}) 
\label{eq:R120a,L49} \end{equation} for $X = R$ and $L$.  This can be re-
expressed as
\begin{equation}  (\pi^{X}_{{\cal{T}}(\cal{A})}\circ \widehat{M})(Q 
\otimes Q^{\prime}) =
((\widehat{M}
\otimes M) \circ (id \otimes
\sigma \otimes id) \circ (\pi^{X}_{{\cal{T}}(\cal{A})} \otimes 
\pi^{X}_{{\cal{T}}(\cal{A})}))(Q
\otimes Q^{\prime})
\label{eq:R120b,L50} \end{equation} for $X = R$ and $L$,  where the 
operator multiplication
operation $\widehat{M}$ is defined by
\begin{equation}  
\widehat{M}(Q \otimes Q^{\prime}) = Q \circ Q^{\prime}
\label{eq:R98} \end{equation}  for all $Q,Q^{\prime} \in 
{\cal{T}}({\cal{A}})$. 

By contrast, if $\widetilde{\pi}^{X}_{{\cal{T}}(\cal{A})}$ is the 
{\em{twisted}} right coaction
defined in (\ref{eq:R115amod}) and (\ref{eq:R115bmod})
 (for $X = R$) and (\ref{eq:R115aLmod}) and (\ref{eq:R115bLmod}) (for $X 
= L$), then
\begin{eqnarray}\begin{array}{l}  
(\widetilde{\pi}^{X}_{{\cal{T}}(\cal{A})}\circ \widehat{M})(Q
\otimes Q^{\prime}) = \\\;\;((\widehat{M}
\otimes M) \circ (id \otimes id \otimes \sigma) \circ (id \otimes
\sigma \otimes id) \circ (\widetilde{\pi}^{X}_{{\cal{T}}(\cal{A})} \otimes
\widetilde{\pi}^{X}_{{\cal{T}}(\cal{A})}))(Q \otimes Q^{\prime}) 
\end{array} \label{eq:R136a,L66}
\end{eqnarray}  for $X = R$ and $L$, whose right-hand side involves an 
extra twist factor $(id
\otimes id \otimes
\sigma)$ relative to the corresponding result (\ref{eq:R120b,L50}).  

The proofs of these statements just involve a straightforward application 
of the identities
(\ref{eq:QD8}) to (\ref{eq:QD27d}).

\section{Theorems of the Wigner-Eckart type}

If $\pi^{p}$, $\pi^{q}$, and $\pi^{r}$ are unitary irreducible 
corepresentations of
$\cal{A}$ of dimensions $d_{p}$, $d_{q}$, and $d_{r}$ respectively, 
$\phi^{pX}_{j}$ and
$\psi^{rX}_{\ell}$ are basis functions belonging to $\pi^{p}$ and 
$\pi^{r}$, and $Q^{qX}_{k}$ is an
{\em{ordinary}} irreducible tensor operator belonging to $\pi^{q}$, then
\begin{eqnarray} (\psi^{rX}_{\ell},Q^{qX}_{k}(\phi^{pX}_{j}))^{X} =
\sum_{\alpha=1}^{n_{qp}^{r}} \left(\begin{array}{ccc} r&,&
\alpha\\\ell & &    
\end{array} \right| \left. \begin{array}{cc}q & p\\ k & j   
\end{array}\right)  (r \mid Q^{qX} \mid
p)_{\alpha}\; ,
\label{eq:R168,L94} \end{eqnarray}  for $X = R$ and $L$, all $j = 1,2, 
\ldots, d_{p}$, all $k = 1,2,
\ldots, d_{q}$, and all $\ell = 1,2, \ldots, d_{r}$. Here the {\em{reduced 
matrix elements}} $(r \mid
Q^{qX} \mid p)_{\alpha}$ are given by
\begin{eqnarray} \begin{array}{ll}(r \mid Q^{qX} \mid p)_{\alpha} & = \\  
&\sum_{s=1}^{d_{p}}
\sum_{t=1}^{d_{q}}
\sum_{u,v=1}^{d_{r}} (\psi^{rX}_{u},Q^{qX}_{t}(\phi^{pX}_{s}))^{X} 
\left( \begin{array}{cc} q & p\\ t &  s \end{array} \right| \left. 
\begin{array}{ccc}r&,& \alpha\\v
 & &  \end{array}\right)  \\  &\times  \;
\{(({\bf{F}}^{r})^{-1})_{vu}/tr(({\bf{F}}^{r})^{-1})\}\end{array}
\label{eq:R171,L95} \end{eqnarray} for $\alpha = 1, 2, \ldots, 
{n_{qp}^{r}}$, and the inner products
$(\;,\;)^{R}$ and
$(\;,\;)^{L}$ are defined in (\ref{eq:R40}) and (\ref{eq:L14}). 

On the other hand, if ${\widetilde{Q}}^{qX}_{k}$ is a {\em{twisted}} 
irreducible tensor operator
belonging to $\pi^{q}$, then
\begin{eqnarray} 
(\psi^{rX}_{\ell},{\widetilde{Q}}^{qX}_{k}(\phi^{pX}_{j}))^{X} =
\sum_{\alpha=1}^{n_{pq}^{r}} \left(\begin{array}{ccc} r&,&
\alpha\\\ell & &    
\end{array} \right| \left. \begin{array}{cc}p & q\\ j & k   
\end{array}\right)  (r \mid
{\widetilde{Q}}^{qX} \mid p)_{\alpha} ,
\label{eq:R176,L101} \end{eqnarray}  for $X = R$ and $L$, all $j = 1,2, 
\ldots, d_{p}$, all $k = 1,2,
\ldots, d_{q}$, and all $\ell = 1,2, \ldots, d_{r}$, where the reduced matrix 
elements $(r \mid
{\widetilde{Q}}^{qX} \mid p)_{\alpha}$ are given by
\begin{eqnarray} \begin{array}{ll}(r \mid {\widetilde{Q}}^{qX} \mid 
p)_{\alpha} & = \\ 
&\sum_{s=1}^{d_{p}}
\sum_{t=1}^{d_{q}}
\sum_{u,v=1}^{d_{r}} 
(\psi^{rX}_{u},{\widetilde{Q}}^{qX}_{t}(\phi^{pX}_{s}))^{X} 
\left( \begin{array}{cc} p & q\\ s &  t \end{array} \right| \left. 
\begin{array}{ccc}r&,& \alpha\\v
 & &  \end{array}\right)  \\  &\times  \;
\{(({\bf{F}}^{r})^{-1})_{vu}/tr(({\bf{F}}^{r})^{-1})\}\end{array}
\label{eq:R177,L102} \end{eqnarray} for $\alpha = 1, 2, \ldots, 
{n_{pq}^{r}}$.

The results (\ref{eq:R168,L94}) and (\ref{eq:R176,L101}) exhibit the 
classic Wigner-Eckart theorem
behaviour, in that they show that the $j$, $k$, and $\ell$ dependences of 
the inner products
$(\psi^{rX}_{\ell},Q^{qX}_{k}\phi^{pX}_{j})^{X}$ and
$(\psi^{rX}_{\ell},{\widetilde{Q}}^{qX}_{k}\phi^{pX}_{j})^{X}$ are 
determined only by Clebsch-Gordan
coefficients, but it should be noted that in the general case in which 
$\cal{A}$ is non-commutative,
the inner products for the {\em{ordinary}} and {\em{twisted}} irreducible 
tensor operators involve
{\em{different}} sets of Clebsch-Gordan coefficients.

The proof of (\ref{eq:R168,L94}) in the case
$X = R$ is as follows. By (\ref{eq:R164v}) and (\ref{eq:R119,L47})
\[  \pi^{R}_{{\cal{A}}} (\psi^{rR*}_{\ell}Q^{qR}_{k}(\phi^{pR}_{j})) = 
\sum_{s=1}^{d_{p}}
\sum_{t=1}^{d_{q}} \sum_{u=1}^{d_{r}} 
(\psi^{rR*}_{u}Q^{qR}_{t}(\phi^{pR}_{s}))\otimes
(\pi^{r*}_{u\ell}\pi^{q}_{tk}\pi^{p}_{sj})
\] for all $j = 1, 2, \ldots, d_{p}$, $k = 1, 2, \ldots, d_{q}$, and $\ell = 1, 2, 
\ldots, d_{r}$.
Then, by (\ref{eq:extra15b})
\[  h(\psi^{rR*}_{\ell}Q^{qR}_{k}(\phi^{pR}_{j})) = \sum_{s=1}^{d_{p}}
\sum_{t=1}^{d_{q}} \sum_{u=1}^{d_{r}}
h(\psi^{rR*}_{u}Q^{qR}_{t}(\phi^{pR}_{s}))\;h(\pi^{r*}_{u\ell}\pi^{q}_{tk}
\pi^{p}_{sj})
\] for all $j = 1, 2, \ldots, d_{p}$, $k = 1, 2, \ldots, d_{q}$, and $\ell = 1, 2, 
\ldots, d_{r}$.
Invoking (\ref{eq:R40}) and (\ref{eq:R164}) then immediately gives 
(\ref{eq:R168,L94}) and
(\ref{eq:R171,L95}).

The proof of (\ref{eq:R168,L94}) in the case $X = L$ is similar. By 
(\ref{eq:R119,L47}),
(\ref{eq:L73}), and (\ref{eq:L90})
\[  \pi^{L}_{{\cal{A}}} 
(Q^{qL}_{k}(\phi^{pL}_{j})(S^{2}(\psi^{rL}_{\ell}))^{*}) =
\sum_{s=1}^{d_{p}}
\sum_{t=1}^{d_{q}} \sum_{u=1}^{d_{r}} 
(Q^{qL}_{t}(\phi^{pL}_{s})(S^{2}(\psi^{rL}_{u}))^{*}) \otimes
(\pi^{r*}_{u\ell}\pi^{q}_{tk}\pi^{p}_{sj})
\] for all $j = 1, 2, \ldots, d_{p}$, $k = 1, 2, \ldots, d_{q}$, and $\ell = 1, 2, 
\ldots, d_{r}$.
On applying  (\ref{eq:extra15b}), (\ref{eq:L14}), and (\ref{eq:R164}), the 
results
(\ref{eq:R168,L94}) and (\ref{eq:R171,L95}) are obtained for this case as 
well. 

The proof of (\ref{eq:R176,L101}) follows the same lines, but employs 
(\ref{eq:R134,L65})
 and (\ref{eq:R163}) in place of (\ref{eq:R119,L47}) and (\ref{eq:R164}).

\section{Products of irreducible tensor operators}

If $\pi^{p}$ and $\pi^{q}$ are unitary irreducible corepresentations of
$\cal{A}$ of dimensions $d_{p}$ and $d_{q}$ respectively, and 
$Q^{pX}_{j}$ and $Q^{qX}_{k}$ are
{\em{ordinary}} irreducible tensor operators belonging to $\pi^{p}$ and 
$\pi^{q}$, then
\begin{equation}  \pi^{X}_{{\cal{T}}({\cal{A}})} (Q^{pX}_{j}Q^{qX}_{k}) = 
\sum_{s=1}^{d_{p}}
\sum_{t=1}^{d_{q}} (Q^{pX}_{s}Q^{qX}_{t}) \otimes (\pi^{p} \sqtimes 
\pi^{q})_{st,jk}
\label{eq:R178} \end{equation} for $X = R$ and $L$ and for all $j = 1,2, 
\ldots, d_{p}$, and $k =
1,2,
\ldots, d_{q}$. Here the coactions $\pi^{R}_{{\cal{T}}({\cal{A}})}$ and
$\pi^{L}_{{\cal{T}}({\cal{A}})}$ are as defined in (\ref{eq:R115a}), 
(\ref{eq:R115b}),
(\ref{eq:R115aL}), and (\ref{eq:R115bL}), and the matrix coefficients of 
the tensor product
$\pi^{p} \sqtimes
\pi^{q}$ are given in (\ref{eq:R155a}). (The proof of (\ref{eq:R178}) just 
involves applying
(\ref{eq:R120b,L50}), (\ref{eq:R115}), (\ref{eq:L41a}), and (\ref{eq:R98}).) 

Because of the similarity in form between (\ref{eq:R178}) and 
(\ref{eq:R154}), it follows
immediately that 
\begin{equation} Q^{rX,\alpha}_{\ell} = \sum_{j=1}^{d_{p}} 
\sum_{k=1}^{d_{q}}
\left(
\begin{array}{cc} p & q\\ j &  k \end{array} \right| \left. 
\begin{array}{ccc}r&,& \alpha\\\ell & &  
\end{array}\right) Q^{pX}_{j}Q^{qX}_{k} ,
\label{eq:R179} \end{equation} for $\ell = 1, 2, \ldots, d_{r}$, and $\alpha 
= 1, 2, \ldots,
n_{pq}^{r}$. Here
$Q^{rX,\alpha}_{\ell}$ (for $\alpha = 1, 2, \ldots, n_{pq}^{r}$) are 
$n_{pq}^{r}$ ordinary
irreducible tensor operators belonging to $\pi^{r}$ that are, in general, all 
different. 

By contrast, if ${\widetilde{Q}}^{pX}_{j}$ and ${\widetilde{Q}}^{qX}_{k}$ 
are {\em{twisted}}
irreducible tensor operators belonging to $\pi^{p}$ and $\pi^{q}$, then
\begin{equation}  {\widetilde{\pi}}^{X}_{{\cal{T}}({\cal{A}})}
({\widetilde{Q}}^{pX}_{j}{\widetilde{Q}}^{qX}_{k}) =
\sum_{s=1}^{d_{p}}
\sum_{t=1}^{d_{q}} ({\widetilde{Q}}^{pX}_{s}{\widetilde{Q}}^{qX}_{t}) 
\otimes (\pi^{p} \twsqtimes
\pi^{q})_{st,jk}
\label{eq:R181} \end{equation} for $X = R$ and $L$ and for all $j = 1,2, 
\ldots, d_{p}$, and $k =
1,2,
\ldots, d_{q}$. Here the right coactions 
${\widetilde{\pi}}^{R}_{{\cal{T}}({\cal{A}})}$ and
${\widetilde{\pi}}^{L}_{{\cal{T}}({\cal{A}})}$ are as defined in 
(\ref{eq:R115amod}), 
(\ref{eq:R115bmod}), (\ref{eq:R115aLmod}), and (\ref{eq:R115bLmod}), 
and the matrix
coefficients of the {\em{twisted}} tensor product $\pi^{p} \twsqtimes
\pi^{q}$ are given in (\ref{eq:R164sa}). (This result (\ref{eq:R181}) is 
proved using
(\ref{eq:R136a,L66}), (\ref{eq:R133}), (\ref{eq:L59}), and (\ref{eq:R98}).) 

The analogue of (\ref{eq:R179}) for the {\em{twisted}} case is  
\begin{equation} {\widetilde{Q}}^{rX,\alpha}_{\ell} = \sum_{j=1}^{d_{p}} 
\sum_{k=1}^{d_{q}}
\left(
\begin{array}{cc} q & p\\ k &  j \end{array} \right| \left. 
\begin{array}{ccc}r&,& \alpha\\\ell & &  
\end{array}\right) {\widetilde{Q}}^{pX}_{j}{\widetilde{Q}}^{qX}_{k} ,
\label{eq:R182} \end{equation} for $\ell = 1, 2, \ldots, d_{r}$, and $\alpha 
= 1, 2, \ldots,
n_{qp}^{r}$. Here
${\widetilde{Q}}^{rX,\alpha}_{\ell}$ (for $\alpha = 1, 2, \ldots, 
n_{qp}^{r}$) are $n_{qp}^{r}$
twisted irreducible tensor operators belonging to $\pi^{r}$ that are, in 
general, all different. 

\section{Generalization to quantum homogeneous spaces}

\subsection{Quantum homogeneous spaces}

The definition and role  of quantum homogeneous spaces are best 
introduced by considering
the situation first in the very well established and  familiar context of a 
compact Lie group
$\cal{G}$.  The  homogeneous space formalism for $\cal{G}$ has
two essential features. Firstly, it is  equivalent to the theory in which 
$\cal{G}$ acts as a
transformation group on an external manifold $\cal{M}$, and, secondly, it 
is closely related to
the regular representation formalisms. Both of these aspects were 
reviewed briefly in
Section~I. 
 
With $\cal{G}$ taken to be a group of tranformations that act on an 
external
manifold $\cal{M}$, select some typical point of $\cal{M}$. Let ${\cal{H}}$ 
be the {\em{isotropy
subgroup}} of ${\cal{G}}$, which consists of all transformations of 
${\cal{G}}$ that send this
point into itself, and let ${\cal{M}}_{0}$ be the {\em{orbit}} of points of 
$\cal{M}$ that
are obtained by acting on this typical point with every transformation of 
${\cal{G}}$. In the case
in which $\cal{G}$ is the rotation group about $O$ and $\cal{M}$ is 
$\Re^{3}$, let ${\bf{r}}_{0}$ of
$\Re^{3}$ be this typical point. Then $\cal{H}$ is the subgroup of all 
rotations about the axis
from $O$ to the point  ${\bf{r}}_{0}$, and ${\cal{M}}_{0}$ is the sphere 
centred on  $O$
that contains the point  ${\bf{r}}_{0}$. With an appropriate choice of 
${\bf{r}}_{0}$,
${\cal{M}}_{0}$ can be parametrised by the spherical polar coordinates 
$\theta$ and
$\phi$. Effectively it is only the functional dependence on $\theta$ and 
$\phi$ that comes into
symmetry arguments, the dependence on the radial distance $r$ being 
immaterial. That is, only the
subspace ${\cal{M}}_{0}$ is actually relevant in the group theoretical 
calculations.
However, it is easily demonstrated that there is a one-to-one 
correspondence between the points of
${\cal{M}}_{0}$ and the set of {\em{left}} cosets $T{\cal{H}}$ of  $\cal{G}$ 
with respect to
$\cal{H}$. Thus the quantities of interest are the subset ${\cal{B}}$ of 
$R({\cal{G}})$ that
consists of those members of $R({\cal{G}})$ which are {\em{constant}} on 
each left coset
$T{\cal{H}}$. Then the operators acting on the members of ${\cal{B}}$ 
that correspond to the
operators $P(T)$ of (\ref{eq:1.17}) may be identified with the left regular 
operators
$\widehat{L}(T)$ of (\ref{eq:L8}), as {\em{restricted}} to act {\em{only}} 
on ${\cal{B}}$. Moreover
the only part of the integral (\ref{eq:1.19}) that is relevant to symmetry 
arguments is the part
involving $\theta$ and $\phi$, which is an integral over 
${\cal{M}}_{0}$, and hence is equivalent to the Haar integral of 
(\ref{eq:Haar}) applied to the
functions of
${\cal{B}}$. Finally, in the homogeneous space version, the irreducible 
tensor operators of
(\ref{eq:5.34}) become  mappings of  ${\cal{B}}$ into  ${\cal{B}}$.

Henceforth the  $\star$-Hopf algebra $R({\cal{G}})$  will be denoted by
$\cal{A}$. Then ${\cal{B}}$ becomes a $\star$-subalgebra of ${\cal{A}}$ 
and a {\em{left}}
coideal of ${\cal{A}}$. (The convention adopted here is that ${\cal{B}}$ is 
said to be a left coideal
of ${\cal{A}}$ if $\Delta(f) \in {\cal{A}} \otimes {\cal{B}}$ for all $f \in 
{\cal{B}}$, and
${\cal{B}}$ is said to be a right coideal of ${\cal{A}}$ if $\Delta(f) \in 
{\cal{B}} \otimes
{\cal{A}}$ for all $f \in {\cal{B}}$). It is also trivially true that ${\cal{B}}$ 
is
$S^{2}$-invariant.

There also exists a parallel version of this theory associated with the 
{\em{right}} regular
representation of ${\cal{G}}$, the operators $\widehat{R}(T)$ of which are 
defined by
$\widehat{R}(T)\,f(T^{\prime}) = f(T^{\prime}T)$  for all $f$ and for all 
$T,T^{\prime} \in
\cal{G}$. Then, for example, if $\cal{G}$ is the group of all rotations in 
this space about $O$ and
$\cal{M}$ is $\Re^{3}$, in place of the transformations of (\ref{eq:1.2}) 
one could define another
set in which 
\begin{equation} ({\bf{r}}^{\prime})^{T} = {\bf{r}}^{T}{\bf{R}}(T) ,
\label{eq:1.2a}\end{equation}  
where ${\bf{r}}^{T}$ denotes the transpose of ${\bf{r}}$. Then, for a typical 
point
${\bf{r}}_{0}$, there is a one-to-one correspondence between the points of 
the orbit ${\cal{M}}_{0}$
of
${\bf{r}}_{0}$ and the set of {\em{right}} cosets
${\cal{H}}T$ of 
$\cal{G}$ with respect to the corresponding isotropy subgroup $\cal{H}$ 
of
${\bf{r}}_{0}$. In this case the quantity of interest is the set ${\cal{B}}$ of 
representative
functions of $\cal{G}$ which are constant on each {\em{right}} coset 
${\cal{H}}T$. Then ${\cal{B}}$
is a $\star$-subalgebra of ${\cal{A}}$ and a {\em{right}} coideal of 
${\cal{A}}$, and the
analogues of operators $P(T)$ of (\ref{eq:1.17}) are the $\widehat{R}(T)$ 
restricted to
${\cal{B}}$.  It is again trivially true that ${\cal{B}}$ is $S^{2}$-invariant.

There are various ways in which these ideas can be generalized to produce 
quantum homogeneous
spaces$^{23-26, 31-36}$, but the present development follows the work of 
Dijkhuizen and
Koornwinder$^{23-26}$. In this formulation one works with a $\star$-
Hopf algebra ${\cal{A}}$ (which
is in general both non-commutative and non-cocommutative), and with a 
$\star$-subalgebra ${\cal{B}}$
of
${\cal{A}}$ that is either a right coideal of ${\cal{A}}$ or is a left coideal of 
${\cal{A}}$. (The
explicit discussion in Refs. 23-26 is actually given for the former situation, 
but clearly the
formulation can also be re-expressed  for the latter situation).  Dijkhuizen 
and
Koornwinder$^{23-26}$ have discussed various other algebraic objects 
that are associated with
${\cal{B}}$, and have shown that in the case of the quantum SU(2) group 
there exists a
one-parameter family of such spaces (called `quantum 2-spheres') which 
are mutually
non-isomorphic, and they have related these to the work of 
Podle\'{s}$^{31}$.   

For the case in which ${\cal{B}}$ is a {\em{left}} coideal of ${\cal{A}}$ it 
will be assumed,
for reasons that will become clear in due course, that ${\cal{B}}$ is 
$S^{2}$-invariant. However,
when ${\cal{B}}$ is a {\em{right}} coideal of ${\cal{A}}$ there is no need to 
make this assumption
when investigating the irreducible tensor operators. Whether this 
assumption is needed  in this
case for other purposes is a matter that has been discussed by Dijkhuizen 
and
Koornwinder$^{23-24}$.
 
\subsection{The restricted right and left regular coactions}

In the special case in which $\cal{A}$ is the dual of a group algebra, the 
operations of $\cal{A}$
corresponding to the right and left regular actions of the group algebra are 
the right and left
regular {\em{coactions}} of $\cal{A}$. Consequently the restrictions of 
right and left regular
actions of the group algebra to $\cal{B}$ correspond to the right and left
regular coactions of $\cal{A}$ restricted to $\cal{B}$.  These are  not only 
the
relevant operations of the classical homogeneous space formulation but 
they are also the basic
operations of the {\em{quantum}} homogeneous space formulation. 

Explicitly, the right regular
coaction   $\pi^{R}_{\cal{A}}$ and the left regular coaction 
$\pi^{L}_{\cal{A}}$ for a general
compact quantum group algebra $\cal{A}$ are defined (in (\ref{eq:R1}) 
and (\ref{eq:L1}) by 
\[ \pi^{R}_{\cal{A}} = \Delta\;,\;\pi^{L}_{\cal{A}} = \sigma \circ (S 
\otimes id) \circ
\Delta .\]
Both are {\em{right}} coactions of $\cal{A}$ with carrier space $\cal{A}$. 
The corresponding
{\em{restricted right regular coaction}}
$\pi^{R}_{\cal{B}}$ and {\em{restricted left regular coaction}} 
$\pi^{L}_{\cal{B}}$ may then be
defined by 
\begin{equation} \pi^{R}_{\cal{B}} = \pi^{R}_{\cal{A}|\cal{B}} =
\Delta_{|\cal{B}}\;,\;\pi^{L}_{\cal{B}} = \pi^{L}_{\cal{A}|\cal{B}} =
\sigma \circ (S \otimes id) \circ \Delta_{|\cal{B}} . 
\label{eq:R1,L1res}\end{equation}
In the context of the restricted right regular coaction it is being assumed 
that $\cal{B}$ is a right
coideal of $\cal{A}$, whereas in the restricted left regular coaction context
$\cal{B}$ is assumed to be a left coideal of $\cal{A}$. Because of the extra 
twist factor $\sigma$
in the definition of $\pi^{L}_{\cal{B}}$, it follows that {\em{both}} 
$\pi^{R}_{\cal{B}}$ {\em{and}}
$\pi^{L}_{\cal{B}}$ are {\em{right}} coactions of $\cal{A}$ with carrier 
space $\cal{B}$. (The role
of $\pi^{R}_{\cal{B}}$ as a transitive $\star$-coaction corresponding to 
the transitive action of
the quantum group associated with $\cal{A}$ on the quantum 
homogeneous space associated with
$\cal{B}$ has been described by Dijkhuizen and Koornwinder$^{23,24}$).

The restriction of the Haar functional $h$ of $\cal{A}$ to $\cal{B}$ 
provides a
positive definite integral for $\cal{B}$ with the invariance properties
\begin{equation} (\MCA \circ (h \otimes id) \circ \pi^{X}_{\cal{B}})(b) = 
(\MAC \circ (id \otimes h)
\circ \pi^{X}_{\cal{B}})(b) = h(b)\;1_{\cal{A}}  
\label{eq:R78a,R78bres} \end{equation}
for all $b \in {\cal{B}}$ and for both $X = R$ and $X = L$ (c.f. 
(\ref{eq:R78a,R78b}) and
(\ref{eq:L19a,L19b})).

The restricted right and left regular corepresentations are both 
{\em{unitary}}, provided that the
inner products on the carrier space $\cal{B}$ are chosen in the following 
way:
\begin{enumerate}
\item for the restricted {\em{right}} regular corepresentation take
\begin{equation}  \langle b , b^{\prime} \rangle_{\cal{B}} = 
(b,b^{\prime})^{R} = h(M(b^{\star}
\otimes b^{\prime}))\;\;{\mbox{for all $b,b^{\prime} \in \cal{B}$}}; 
\label{eq:R40res} \end{equation}
\item for the restricted {\em{left}} regular corepresentation take
\begin{equation}  \langle b , b^{\prime} \rangle_{\cal{B}} = 
(b,b^{\prime})^{L} = h(M(b \otimes
(S^{2}(b^{\prime}))^{\star}))\;\; {\mbox{for all $b,b^{\prime} \in 
\cal{B}$}} .
\label{eq:L14res} \end{equation}
\end{enumerate}
It should be noted that $(S^{2}(b^{\prime}))^{\star} \in {\cal{B}}$ if 
${\cal{B}}$ is
$S^{2}$-invariant. (The proofs of these unitary properties are essentially 
the same as for the
corresponding results (\ref{eq:R40}) and (\ref{eq:L14}) for the unrestricted 
regular coactions).

The effects of the restricted right and left regular coactions on products 
are also essentially the
same as for the unrestricted coactions (c.f.(\ref{eq:R164v}) and 
(\ref{eq:L73})). 
   
\subsection{Basis functions for the restricted right and left regular 
coactions}

Suppose that $\pi^{p}_{jk}$ are the matrix coefficients of a 
corepresentation $\pi^{p}$ of
$\cal{A}$ of finite dimension $d_{p}$. Then the {\em{basis functions}} 
$\psi^{pR}_{j}$ {\em{of}}
$\pi^{p}$ {\em{with respect to the restricted right regular coaction}} may 
be defined to be a set of
$d_{p}$ elements of
$\cal{B}$ that have the property that   
\begin{equation} \pi^{R}_{{\cal{B}}}(\psi^{pR}_{j}) = \sum_{k=1}^{d_{p}} 
\psi^{pR}_{k} \otimes
\pi^{p}_{kj}
\label{eq:R79res} \end{equation} for all $j = 1, 2, \ldots, d_{p}$. Similarly 
the {\em{basis
functions}} $\psi^{pL}_{j}$ {\em{of}} $\pi^{p}$ {\em{with respect to the 
restricted left regular
coaction}} may be defined as a set of $d_{p}$ elements of $\cal{B}$ that 
have the property that   
\begin{equation} \pi^{L}_{{\cal{B}}}(\psi^{pL}_{j}) = \sum_{k=1}^{d_{p}} 
\psi^{pL}_{k} \otimes
\pi^{p}_{kj}
\label{eq:L20res} \end{equation} for all $j = 1, 2, \ldots, d_{p}$. In both 
cases the matrix
coefficients $\pi^{p}_{kj}$ are elements of $\cal{A}$, and need not be 
members of $\cal{B}$.

By contrast with the unrestricted coaction situation, in neither case is 
there any guarantee
that for a given corepresentation $\pi^{p}$ there actually exist basis 
functions for restricted
coactions. For example, in the restricted right regular coaction case, a set 
of basis
functions is provided (for any fixed choice of $\ell = 1, 2, \ldots$) by
\begin{equation} \psi^{pR}_{j} = \pi^{p}_{\ell j}
\label{eq:R80res} \end{equation} for all $j = 1, 2, \ldots, d_{p}$ {\em{only 
if}}  $\pi^{p}_{\ell j}$
is an element of $\cal{B}$ for {\em{some}} $j = 1, 2, \ldots, d_{p}$. (Then 
(I.43) and the fact that
$\cal{B}$ is assumed in this context to be a right coideal of $\cal{A}$ 
imply that $\pi^{p}_{\ell j}
\in \cal{B}$ for {\em{all}} $j = 1, 2, \ldots, d_{p}$). Similarly, in the 
restricted left regular
coaction case, an example is provided (for any fixed choice of $\ell = 1, 2, 
\ldots$) by
\begin{equation} \psi^{pL}_{j} = S^{-2}(\pi^{p\star}_{j\ell})
\label{eq:L22res} \end{equation} for all $j = 1, 2, \ldots, d_{p}$ {\em{only 
if}}  $\pi^{p}_{j \ell}$
is an element of $\cal{B}$ for {\em{some}} $j = 1, 2, \ldots, d_{p}$. (Then 
(I.43) and the fact that
$\cal{B}$ is assumed in this context to be a $S^{2}$-invariant left coideal 
and $\star$-subalgebra
of $\cal{A}$ imply that $S^{-2}(\pi^{p\star}_{j\ell}) \in \cal{B}$ for 
{\em{all}} $j = 1, 2,
\ldots, d_{p}$).

One very useful result, proved in the same way as the corresponding 
unrestricted identity
(\ref{eq:L90}), is that if $\psi^{qL}_{k}$ exists then
\begin{equation} \pi^{L}_{{\cal{B}}}((S^{2}(\psi^{qL}_{k}))^{\star}) = 
\sum_{t=1}^{d_{p}}
(S^{2}(\psi^{qL}_{t}))^{\star} \otimes \pi^{q\star}_{tk} 
\label{eq:L90res} \end{equation} for all $k = 1, 2, \ldots, d_{p}$.

The orthogonality properties of basis functions are the essentially the 
same as for those
for the unrestricted case that have been given in (\ref{eq:R81c}), 
(\ref{eq:R81d}), and
(\ref{eq:R81g}). (The only qualification is that now it has to be assumed 
that the relevant matrix
corepresentation coefficients are members of
$\cal{B}$).   

If the basis functions $\psi^{qX}_{k}$ and $\phi^{pX}_{j}$ are members of 
$\cal{B}$ (so that they
are basis functions for the restricted coactions), then their products 
$M(\psi^{qX}_{k}
\otimes \phi^{pX}_{j})$ are also members of $\cal{B}$. Consequently the 
analysis of Subsection V.D 
goes through without modification, except that all basis functions involved 
(including the
$\theta^{r,\alpha X}_{\ell}$) are in $\cal{B}$ and one can always replace 
the unrestricted
coactions $\pi^{X}_{{\cal{A}}}$ by the restricted coactions 
$\pi^{X}_{{\cal{B}}}$. The rest of
discussion of Section V on tensor products and Clebsch-Gordan coefficients 
still
applies in its entirety. 

\subsection{The irreducible tensor operators in the restricted 
corepresentation
formalisms}

\subsubsection{Introduction}

Let $\pi^{q}$ be a unitary irreducible right coaction of $\cal{A}$ of 
dimension
$d_{q}$ with matrix coefficients $\pi^{q}_{jk}$. It will be shown that  
within both the
restricted right and the restricted left regular corepresentation formalisms 
there exist two types of
irreducible tensor operators that both belong to this corepresentation
$\pi^{q}$. These will be denoted by
$Q^{qX}_{j\cal{B}}$ and ${\widetilde{Q}}^{qX}_{j\cal{B}}$ (for $j = 1, 2, 
\ldots, d_{q}$ and for $X
= R$ or $L$), and will be called {\em{ordinary}} and {\em{twisted}} 
irreducible tensor operators
associated with $\cal{B}$ respectively. These operators are members of 
$\cal{L}(\cal{B})$, which is
the set of linear operators that map $\cal{B}$ into $\cal{B}$. Naturally 
the two types of irreducible
tensor operators coincide in the special case in which $\cal{A}$ is 
commutative.   As much of the
analysis is the same as for the unrestricted case, all the proofs will be 
omitted.

\subsubsection{Definition of the ordinary irreducible tensor operators
$Q^{qX}_{j\cal{B}}$ and twisted irreducible tensor operators
${\widetilde{Q}}^{qX}_{j\cal{B}}$}

The {\em{ordinary irreducible tensor operators}}
$Q^{qX}_{j\cal{B}}$ belonging to the unitary irreducible right coaction 
$\pi^{q}$ of
$\cal{A}$ are {\em{defined}} (for $X = L,R$) to be members of 
$\cal{L}(\cal{B})$ that satisfy the
condition 
\begin{equation}  ((id \otimes M) \circ (\pi^{X}_{\cal{B}} \otimes id) \circ 
(Q^{qX}_{j\cal{B}}
\otimes S) \circ \pi^{X}_{\cal{B}})(b) = \sum_{k=1}^{d_{q}} 
Q^{qX}_{k\cal{B}}(b) \otimes
\pi^{q}_{kj}
\label{eq:L46res} \end{equation}  for all $b \in \cal{B}$ and all $j = 1, 2, 
\ldots, d_{q}$.

The {\em{twisted irreducible tensor operators}} 
${\widetilde{Q}}^{qX}_{j\cal{B}}$ belonging to the
unitary irreducible right coaction $\pi^{q}$ of $\cal{A}$ are 
{\em{defined}} (for $X = L,R$) to be
members of
${\cal{L}}({\cal{B}})$ that satisfy the condition
\begin{equation}  ((id \otimes M) \circ (id \otimes \sigma) \circ 
(\pi^{X}_{\cal{B}} \otimes id)
\circ ({\widetilde{Q}}^{qX}_{j\cal{B}}
\otimes S^{-1}) \circ \pi^{X}_{\cal{B}})(b) = \sum_{k=1}^{d_{q}} 
{\widetilde{Q}}^{qX}_{k\cal{B}}(b)
\otimes
\pi^{q}_{kj}
\label{eq:extra21res} \end{equation}  for all $b \in \cal{B}$ and all $j = 1, 
2, \ldots, d_{q}$.
(This definition (\ref{eq:extra21res}) differs from the corresponding 
definition (\ref{eq:L46res})
only in the replacement of $M$ by $M \circ \sigma$ {\em{and}} $S$ by 
$S^{-1}$ (neither of which have
any effect in the special case in which $\cal{A}$ is commutative)).

\subsubsection{Properties of the irreducible tensor operators associated 
with quantum
homogeneous spaces}

\begin{enumerate}

\item Suppose that $Q$ is the {\em{identity operator}} $id$ of 
${\cal{L}}({\cal{B}})$ (so
that
$Q(b) = b$ for all $b
\in \cal{B}$). Then  $id$ is both an {\em{ordinary}} and a {\em{twisted}} 
irreducible tensor operator
for the one-dimensional {\em{identity}} corepresentation of $\cal{A}$ 
(whose sole matrix coefficient
is $1_{{\cal{A}}}$) in the restricted {\em{right}} regular corepresentation 
formalism, as well
as being a both an {\em{ordinary}} and a {\em{twisted}} irreducible tensor 
operator  for
this identity corepresentation in the restricted {\em{left}} regular 
corepresentation formalism.

\item Suppose now that $\psi^{qR}_{j}$ and $\psi^{qL}_{j}$ are sets of 
basis functions for $\pi^{q}$,
as defined in (\ref{eq:R79res}) and (\ref{eq:L20res}) respectively (so that 
they are members of
$\cal{B}$), and suppose that the operators
$Q^{qX}_{j\cal{B}}$ and ${\widetilde{Q}}^{qX}_{j\cal{B}}$ are 
{\em{defined}} by
\begin{eqnarray} \left. \begin{array}{ccc}  Q^{qR}_{j\cal{B}}(b) & = & 
M(\psi^{qR}_{j} \otimes b) ,
\\ {\widetilde{Q}}^{qR}_{j\cal{B}}(b) &  = & M(b \otimes \psi^{qR}_{j}) 
,\\ Q^{qL}_{j\cal{B}}(b) & =
& M(b
\otimes
\psi^{qL}_{j}) , \\ {\widetilde{Q}}^{qL}_{j\cal{B}}(b) & = & M(\psi^{qL}_{j} 
\otimes b) , \end{array}
\right\}
\label{eq:R118,R136,L44.L62res} \end{eqnarray}   for all $b \in \cal{B}$. 
Then $Q^{qR}_{j\cal{B}}$,
${\widetilde{Q}}^{qR}_{j\cal{B}}$, $Q^{qL}_{j\cal{B}}$, and
${\widetilde{Q}}^{qL}_{j\cal{B}}$   are
irreducible tensor operators  belonging to $\pi^{q}$.

\end{enumerate}

\subsection{Wigner-Eckart type theorems associated with quantum 
homogeneous
spaces}

If $\pi^{p}$, $\pi^{q}$, and $\pi^{r}$ are unitary irreducible 
corepresentations of
$\cal{A}$ of dimensions $d_{p}$, $d_{q}$, and $d_{r}$ respectively, 
$\phi^{pX}_{j}$ and
$\psi^{rX}_{\ell}$ are basis functions belonging $\pi^{p}$ and $\pi^{r}$ 
(with $\phi^{pX}_{j}$ and
$\psi^{rX}_{\ell}$ being assumed here to be members of $\cal{B}$), and 
$Q^{qX}_{k\cal{B}}$ is an
{\em{ordinary}} irreducible tensor operator belonging to $\pi^{q}$ (with 
respect to the
relevant restricted  regular coaction), then
\begin{eqnarray} (\psi^{rX}_{\ell},Q^{qX}_{k\cal{B}}(\phi^{pX}_{j}))^{X} =
\sum_{\alpha=1}^{n_{qp}^{r}} \left(\begin{array}{ccc} r&,&
\alpha\\\ell & &    
\end{array} \right| \left. \begin{array}{cc}q & p\\ k & j   
\end{array}\right)  (r \mid
Q^{qX}_{\cal{B}} \mid p)_{\alpha}\; ,
\label{eq:R168,L94res} \end{eqnarray}  for $X = R$ and $L$, all $j = 1,2, 
\ldots, d_{p}$, all $k =
1,2,
\ldots, d_{q}$, and all $\ell = 1,2, \ldots, d_{r}$. Here the {\em{reduced 
matrix elements}} $(r \mid
Q^{qX}_{\cal{B}} \mid p)_{\alpha}$ are given by
\begin{eqnarray} \begin{array}{ll}(r \mid Q^{qX}_{\cal{B}} \mid 
p)_{\alpha} & = \\ 
&\sum_{s=1}^{d_{p}}
\sum_{t=1}^{d_{q}}
\sum_{u,v=1}^{d_{r}} (\psi^{rX}_{u},Q^{qX}_{t\cal{B}}(\phi^{pX}_{s}))^{X} 
\left( \begin{array}{cc} q & p\\ t &  s \end{array} \right| \left. 
\begin{array}{ccc}r&,& \alpha\\v
 & &  \end{array}\right)  \\  &\times  \;
\{(({\bf{F}}^{r})^{-1})_{vu}/tr(({\bf{F}}^{r})^{-1})\}\end{array}
\label{eq:R171,L95res} \end{eqnarray} for $\alpha = 1, 2, \ldots, 
{n_{qp}^{r}}$, the inner products
$(\;,\;)^{R}$ and $(\;,\;)^{L}$ are defined in (\ref{eq:R40res}) and 
(\ref{eq:L14res}), and
${\bf{F}}^{r}$ is  defined in (\ref{eq:R148}).

On the other hand, if ${\widetilde{Q}}^{qX}_{k\cal{B}}$ is a {\em{twisted}} 
irreducible tensor
operator belonging to $\pi^{q}$ (with respect to the
relevant restricted  regular coaction), then
\begin{eqnarray} 
(\psi^{rX}_{\ell},{\widetilde{Q}}^{qX}_{k\cal{B}}(\phi^{pX}_{j}))^{X} =
\sum_{\alpha=1}^{n_{pq}^{r}} \left(\begin{array}{ccc} r&,&
\alpha\\\ell & &    
\end{array} \right| \left. \begin{array}{cc}p & q\\ j & k   
\end{array}\right)  (r \mid
{\widetilde{Q}}^{qX}_{\cal{B}} \mid p)_{\alpha} ,
\label{eq:R176,L101res} \end{eqnarray}  for $X = R$ and $L$, all $j = 1,2, 
\ldots, d_{p}$, all $k =
1,2,
\ldots, d_{q}$, and all $\ell = 1,2, \ldots, d_{r}$, where the reduced matrix 
elements $(r \mid
{\widetilde{Q}}^{qX}_{\cal{B}} \mid p)_{\alpha}$ are given by
\begin{eqnarray} \begin{array}{ll}(r \mid {\widetilde{Q}}^{qX}_{\cal{B}} 
\mid p)_{\alpha} & = \\ 
&\sum_{s=1}^{d_{p}}
\sum_{t=1}^{d_{q}}
\sum_{u,v=1}^{d_{r}} 
(\psi^{rX}_{u},{\widetilde{Q}}^{qX}_{t\cal{B}}(\phi^{pX}_{s}))^{X} 
\left( \begin{array}{cc} p & q\\ s &  t \end{array} \right| \left. 
\begin{array}{ccc}r&,& \alpha\\v
 & &  \end{array}\right)  \\  &\times  \;
\{(({\bf{F}}^{r})^{-1})_{vu}/tr(({\bf{F}}^{r})^{-1})\}\end{array}
\label{eq:R177,L102res} \end{eqnarray} for $\alpha = 1, 2, \ldots, 
{n_{pq}^{r}}$.

These results (\ref{eq:R168,L94res}) and (\ref{eq:R176,L101res}) 
demonstrate that again the $j$, $k$,
and
$\ell$ dependences of the inner products
$(\psi^{rX}_{\ell},Q^{qX}_{k\cal{B}}\phi^{pX}_{j})^{X}$ and
$(\psi^{rX}_{\ell},{\widetilde{Q}}^{qX}_{k\cal{B}}\phi^{pX}_{j})^{X}$ are 
determined only by
Clebsch-Gordan coefficients, and so they have the same form as in the 
classic Wigner-Eckart theorem.
(However, it should be noted that in the general case in which $\cal{A}$ is 
non-commutative, the
inner products for the {\em{ordinary}} and {\em{twisted}} irreducible 
tensor operators involve
{\em{different}} sets of Clebsch-Gordan coefficients).

As the proofs of (\ref{eq:R168,L94res}) and (\ref{eq:R176,L101res}) follow 
the same lines as in the
unrestricted case considered in Section VII, they will be omitted here.

\subsection{Products of irreducible tensor operators associated with 
quantum homogeneous
spaces}

The arguments of Section~VIII can be applied (with $\cal{A}$ replaced by 
$\cal{B}$) to show that
\begin{equation} Q^{rX,\alpha}_{\ell\cal{B}} = \sum_{j=1}^{d_{p}} 
\sum_{k=1}^{d_{q}}
\left(
\begin{array}{cc} p & q\\ j &  k \end{array} \right| \left. 
\begin{array}{ccc}r&,& \alpha\\\ell & &  
\end{array}\right) Q^{pX}_{j\cal{B}}Q^{qX}_{k\cal{B}} ,
\label{eq:R179res} \end{equation} for $\ell = 1, 2, \ldots, d_{r}$, and 
$\alpha = 1, 2, \ldots,
n_{pq}^{r}$. Here
$Q^{rX,\alpha}_{\ell\cal{B}}$ (for $\alpha = 1, 2, \ldots, n_{pq}^{r}$) are 
$n_{pq}^{r}$ ordinary
irreducible tensor operators belonging to $\pi^{r}$ that are, in general, all 
different. Moreover,  
\begin{equation} {\widetilde{Q}}^{rX,\alpha}_{\ell\cal{B}} = 
\sum_{j=1}^{d_{p}} \sum_{k=1}^{d_{q}}
\left(
\begin{array}{cc} q & p\\ k &  j \end{array} \right| \left. 
\begin{array}{ccc}r&,& \alpha\\\ell & &  
\end{array}\right) 
{\widetilde{Q}}^{pX}_{j\cal{B}}{\widetilde{Q}}^{qX}_{k\cal{B}} ,
\label{eq:R182res} \end{equation} for $\ell = 1, 2, \ldots, d_{r}$, and 
$\alpha = 1, 2, \ldots,
n_{qp}^{r}$, where
${\widetilde{Q}}^{rX,\alpha}_{\ell\cal{B}}$ (for $\alpha = 1, 2, \ldots, 
n_{qp}^{r}$) are
$n_{qp}^{r}$ twisted irreducible tensor operators belonging to $\pi^{r}$ 
that are, in general, all
different.

\acknowledgements

This work was started while the author was on research leave at the 
Sektion Physik of the
Ludwig-Maximilians-Universit\"{a}t, M\"{u}nchen. The author is grateful 
to Professor J. Wess and his
colleagues for their hospitality and to the Deutscher Akademischer 
Austauschdienst for its financial
support.

\appendix

\section{Introduction}

The purpose of this Appendix is to {\em{motivate}} the definitions that are 
given in the main body of
the paper for the irreducible tensor operators and the projection operators. 
This will be done by
considering the simple special case in which the Hopf algebra $\cal{A}$ is 
the set of functions
defined on a {\em{finite}} group $\cal{G}$ of order $g$, so that the dual 
$\cal{A}^{\prime}$ of
$\cal{A}$ is the group algebra of $\cal{G}$. In this situation both 
$\cal{A}$ and $\cal{A}^{\prime}$
are of finite dimension $g$. Of course, as $\cal{A}$ is commutative in this 
special case, the
resulting expressions are to some extent ambiguous, in that in this special 
case
$M$ is indistinguishable from $M \circ \sigma$ and $S$ is 
indistinguishable from
$S^{-1}$. The demonstration of the correctness, consistency, and 
usefulness of the definitions that
are actually employed for the {\em{general}} case are the subject matter of 
the self-contained
arguments of the main body of this paper.         

To proceed with this motivation, it is necessary to start with some well-
known facts concerning the
relationship of $\cal{A}$ and $\cal{A}^{\prime}$. It is easily shown that
$\cal{A}^{\prime}$ is also a Hopf algebra, whose multiplication operation 
$M_{\cal{A}^{\prime}}$,
comultiplication operation $\Delta_{\cal{A}^{\prime}}$, and antipode 
$S_{\cal{A}^{\prime}}$ are
related to those of $\cal{A}$ by
\begin{equation} \langle M_{\cal{A}^{\prime}}(a^{\prime} \otimes 
b^{\prime}), a \rangle  = \langle
(a^{\prime} \otimes b^{\prime}), \Delta(a) \rangle
\label{eq:QD33} \end{equation} for all $a \in \cal{A}$ and all 
$a^{\prime},b^{\prime} \in
\cal{A}^{\prime}$,
\begin{equation} \langle \Delta_{\cal{A}^{\prime}}(a^{\prime}), (a 
\otimes b) \rangle  =
\langle a^{\prime}, M(a \otimes b) \rangle
\label{eq:QD37} \end{equation} for all $a,b \in \cal{A}$ and all 
$a^{\prime} \in \cal{A}^{\prime}$,
and 
\begin{equation} \langle S_{\cal{A}^{\prime}}(a^{\prime}), a \rangle  =
\langle a^{\prime}, S(a) \rangle
\label{eq:QD47} \end{equation} for all $a \in \cal{A}$ and all $a^{\prime} 
\in \cal{A}^{\prime}$.

 Now suppose that $\pi_{V}$ is a right coaction of $\cal{A}$ with carrier 
space $V$. Then there
exists a corresponding {\em{left action}} $\pi^{\prime}_{V}$ of 
$\cal{A}^{\prime}$ with the
{\em{same}} carrier space $V$. This is a linear mapping from 
$\cal{A}^{\prime} \otimes
V$ to $V$ such that
\begin{equation} \pi^{\prime}_{V} \circ (id \otimes \pi^{\prime}_{V}) = 
\pi^{\prime}_{V}\circ
(M_{\cal{A}^{\prime}} \otimes id) 
\label{eq:R15} \end{equation} and
\begin{equation} \pi^{\prime}_{V} \circ (u_{\cal{A}^{\prime}} \otimes id) 
= \MCV ,
\label{eq:R16} \end{equation} where $u_{\cal{A}^{\prime}}$ is the unit 
operator of
$\cal{A}^{\prime}$, and  $\MCV(z \otimes v) = zv$ for all $v \in V$ and 
all $z \in \C$. The
relationship between $\pi_{V}$ and
$\pi^{\prime}_{V}$ can be expressed as
\begin{equation} \pi^{\prime}_{V}(a^{\prime} \otimes v) = \sum_{[v]} 
v_{[1]} \;\langle a^{\prime} ,
v_{[2]} \rangle
\label{eq:R13} \end{equation} for all $a^{\prime} \in \cal{A}^{\prime}$ 
and all $v \in V$, where the
notation of (\ref{eq:R12}) has been employed, or equivalently as
\begin{equation} \pi^{\prime}_{V}(a^{\prime} \otimes v) = (\MVC \circ 
(id \otimes ev) \circ (\sigma
\otimes id) \circ (id \otimes \pi_{V}))(a^{\prime} \otimes v) 
\label{eq:R14} \end{equation} for all $a^{\prime} \in \cal{A}^{\prime}$ 
and all $v \in V$, where the
evaluation map $ev$ of (\ref{eq:extra10}) has been used.  The inverse of 
these is
\begin{equation} \pi_{V}(v) = \sum_{j} \pi^{\prime}_{V}(a^{j} \otimes v) 
\otimes a_{j} ,
\label{eq:R26b} \end{equation} for all $v \in V$, where the basis of 
$\cal{A}$ and the dual basis of
$\cal{A}^{\prime}$ appear, and are assumed to be such that 
(\ref{eq:extra11}) holds. If $V$ is of
dimension $d$ with basis elements $v_{1}, v_{2}, \ldots, v_{d}$, the matrix 
elements
$\pi^{\prime}(a^{\prime})_{jk}$ of the representation are such that
\begin{equation} \pi^{\prime}_{V}(a^{\prime} \otimes v_{j}) =
\sum_{k=1}^{d}\pi^{\prime}(a^{\prime})_{kj}\; v_{k} ,
\label{eq:R26c} \end{equation} for all $a^{\prime} \in \cal{A}^{\prime}$ 
and $j = 1, 2, \ldots, d$.
It then follows from (\ref{eq:R26b}) that  these representation matrix 
elements
$\pi^{\prime}(a^{\prime})_{jk}$ are related to the corepresentation 
matrix coefficients of
$\pi_{jk}^{V}$ of (\ref{eq:R5}) by
\begin{equation} \pi^{\prime}(a^{\prime})_{jk} = \langle a^{\prime} , 
\pi_{jk}^{V}\rangle ,
\label{eq:R26d} \end{equation} for all $a^{\prime} \in \cal{A}^{\prime}$ 
and $j,k = 1, 2, \ldots, d$.

	The {\em{right regular action}} $R$ of $\cal{G}$ is defined by
\begin{equation} (R(x \otimes f))(y) = f(yx) ,
\label{eq:R30} \end{equation} for all elements $x,y \in \cal{G}$ and all 
functions $f$ defined on
$\cal{G}$, where $yx$ is evaluated using the group multiplication 
operation of $\cal{G}$.  This
definition (\ref{eq:R30}) can be immediately extended to all $x,y \in 
\cal{A}^{\prime}$, with $f$
being any element of
$\cal{A}$. It is then easy to verify that $R$ is a {\em{left}} action of 
$\cal{A}^{\prime}$ with
carrier space $\cal{A}$, and, using  (\ref{eq:R13}) (or its equivalents), that 
it is the left action
that corresponds to the right regular coaction $\pi^{R}_{\cal{A}}$ of 
$\cal{A}$ that was defined in
(\ref{eq:R8}). Then, by (\ref{eq:R14}), 
\begin{equation} R(x \otimes f) = (\MAC \circ (id \otimes ev) \circ 
(\sigma
\otimes id) \circ (id \otimes \Delta))(x \otimes f) 
\label{eq:R32} \end{equation} for all $x \in \cal{A}^{\prime}$ and all $f 
\in \cal{A}$, and also, by
(\ref{eq:R8}) and (\ref{eq:R13}),
\begin{equation} \langle y ,R(x \otimes f) \rangle = \langle y \otimes x , 
\Delta(f) \rangle =
\langle y \otimes x , \pi^{R}_{\cal{A}}(f) \rangle
\label{eq:R32prime} \end{equation} for all $x,y \in \cal{A}^{\prime}$ and 
all $f \in \cal{A}$. The
content of (\ref{eq:R30}) can usefully be re-expressed as
\begin{equation} \widehat{R}(x)f(y) = f(yx) ,
\label{eq:R29} \end{equation} by introducing an operator 
$\widehat{R}(x)$ for each $x \in
\cal{A}^{\prime}$.

	Similarly, the {\em{left regular action}} $L$ of $\cal{G}$ is defined 
by
\begin{equation} (L(x \otimes f))(y) = f(x^{-1}y) ,
\label{eq:L9} \end{equation} for all elements $x,y \in \cal{G}$ and all 
functions $f$ defined on
$\cal{G}$, which immediately extends to all $x,y \in \cal{A}^{\prime}$ 
and all $f \in \cal{A}$. Then
$L$ is the left action of $\cal{A}^{\prime}$ that corresponds to the left 
regular coaction
$\pi^{L}_{\cal{A}}$ of $\cal{A}$ that was defined in (\ref{eq:L1}). Thus, by 
(\ref{eq:R14}), 
\begin{equation} L(x \otimes f) = (\MAC \circ (id \otimes ev) \circ 
(\sigma
\otimes id) \circ (id \otimes \pi^{L}_{{\cal{A}}}))(x \otimes f) 
\label{eq:L11} \end{equation} for all $x \in \cal{A}^{\prime}$ and all $f 
\in \cal{A}$, and also, by
(\ref{eq:L1}) and (\ref{eq:R13}),
\begin{equation} \langle y ,L(x \otimes f) \rangle = \langle y \otimes x , 
(\sigma \circ (S
\otimes id) \circ \Delta)(f) \rangle =
\langle y \otimes x , \pi^{L}_{\cal{A}}(f) \rangle
\label{eq:L12} \end{equation} for all $x,y \in \cal{A}^{\prime}$ and all $f 
\in \cal{A}$. Moreover,
(\ref{eq:L9}) can be re-written as
\begin{equation} \widehat{L}(x)f(y) = f(S_{{\cal{A}}^{\prime}}(x)y) ,
\label{eq:L8} \end{equation} by introducing an operator $\widehat{L}(x)$ 
for each $x \in
\cal{A}^{\prime}$.

\section{Irreducible tensor operators}

\subsection{Outline of argument}

The first stage is to recall  the definition of irreducible tensor operators in 
the group
theoretical context. The next stage is to cast these considerations into the 
language of Hopf
algebras, and the final stage is to put them into a form in which they 
involve quantities belonging
{\em{only}} to $\cal{A}$. The argument will be given first for the 
{\em{right}} regular situation,
and then for the {\em{left}} regular situation. In each case the 
{\em{definitions}} of irreducible
tensor operators will be deduced first, and motivation for the {\em{right 
coactions}} that appear in
the consistency arguments will follow.  (The diagrammatic method that is 
described, for example, by
Majid$^{37}$, was employed to deduce the proofs that follow, but for 
typographical
convenience these proofs have been transcribed here into the usual purely 
symbolic form).

\subsection{Irreducible tensor operators in the group theoretical right 
regular representation
formalism}

\subsubsection{Derivations of the conditions for the $Q^{qR}_{j}$ and
${\widetilde{Q}}^{qR}_{j}$}

 Let
${\bf{\Gamma}}^{q}$ be a
$d_{q}$ dimensional representation of
$\cal{G}$. Then the irreducible tensor operators $Q^{qR}_{j}$ may be 
defined to act on functions
defined on
$\cal{G}$ in such a way that
\begin{equation} \widehat{R}(x) Q^{qR}_{j} \widehat{R}(x^{-1}) = 
\sum_{k=1}^{d_{q}}
\Gamma^{q}_{kj}(x) Q^{qR}_{k} ,
\label{eq:R82} \end{equation} for all $x \in \cal{G}$ and $j = 1, 2, \dots, 
d_{q}$, or, more
explicitly, such that
\begin{equation} (\widehat{R}(x) Q^{qR}_{j} \widehat{R}(x^{-1}))(f(y)) = 
\sum_{k=1}^{d_{q}}
\Gamma^{q}_{kj}(x) (Q^{qR}_{k} (f))(y) ,
\label{eq:R83} \end{equation} for all $x,y \in \cal{G}$, for all functions 
$f$ defined on $\cal{G}$,
and $j = 1, 2, \dots, d_{q}$. Now define 
$\widehat{\pi}^{R\;\prime}_{{\cal{L}}(\cal{A})}(x)$ by
\begin{equation} \widehat{\pi}^{R\;\prime}_{{\cal{L}}(\cal{A})}(x)(Q) = 
\widehat{R}(x) Q
\widehat{R}(x^{-1})
\label{eq:R84} \end{equation} for all $x \in \cal{G}$ and for all linear 
operators
$Q$ that act on functions defined on $\cal{G}$, so that (\ref{eq:R82}) 
becomes
\begin{equation} 
\widehat{\pi}^{R\;\prime}_{{\cal{L}}(\cal{A})}(x)(Q^{qR}_{j}) =
\sum_{k=1}^{d_{q}}
\Gamma^{q}_{kj}(x) Q^{qR}_{k} ,
\label{eq:R85} \end{equation} for all $x \in \cal{G}$ and $j = 1, 2, \dots, 
d_{q}$. As discussed
previously in Section~I (in only a slightly different context), the 
consistency of the
definition (\ref{eq:R85}) is consequence of the assumption that
${\bf{\Gamma}}^{q}$ is a representation of $\cal{G}$ and the fact that
\begin{equation} \widehat{\pi}^{R\;\prime}_{{\cal{L}}(\cal{A})}(xy) =
\widehat{\pi}^{R\;\prime}_{{\cal{L}}(\cal{A})}(x)\;
\widehat{\pi}^{R\;\prime}_{{\cal{L}}(\cal{A})}(y)
\label{eq:R86} \end{equation} for all $x,y \in \cal{G}$.

As  
\begin{equation} \Delta_{\cal{A}^{\prime}}(x) = x \otimes x 
\;,\;S_{\cal{A}^{\prime}}(x) = x^{-1}
\label{eq:R87,R88} \end{equation} for all $x \in \cal{G}$, it follows from 
(\ref{eq:R29}) and
(\ref{eq:R30}) that
\begin{eqnarray} \begin{array}{l}
(\widehat{\pi}^{R\;\prime}_{{\cal{L}}(\cal{A})}(x)(Q^{qR}_{j}))f(y) = \\ 
(ev \circ (id \otimes R)
\circ (id \otimes id \otimes Q^{qR}_{j}) \circ (id \otimes id \otimes R) 
\circ (id \otimes id \otimes
S_{\cal{A}^{\prime}}\otimes id) \\ 
\circ (id \otimes \Delta_{\cal{A}^{\prime}} \otimes id))(y \otimes x 
\otimes f) ,
\end{array} \label{eq:R88prime} \end{eqnarray} which is now well-
defined for all $x,y \in
\cal{A}^{\prime}$, for all $f \in \cal{A}$, and 
$Q^{qR}_{j} \in {\cal{L}}({\cal{A}})$. Thus, by (\ref{eq:R32}),
\begin{eqnarray}
\begin{array}{l}(\widehat{\pi}^{R\;\prime}_{{\cal{L}}(\cal{A})}(x)(Q^{qR}
_{j}))f(y) = \\ (ev \circ
(id \otimes \MAC) \circ (id \otimes id \otimes ev) \circ (id \otimes 
\sigma
\otimes id) \circ (id \otimes id \otimes \Delta) \\ 
\circ (id \otimes id \otimes Q^{qR}_{j}) \circ (id \otimes id \otimes 
\MAC) \circ (id \otimes id
\otimes id
\otimes ev) \\ 
\circ (id \otimes id \otimes \sigma \otimes id) 
\circ (id \otimes id \otimes id \otimes \Delta)
\circ (id \otimes id \otimes S_{\cal{A}^{\prime}}\otimes id)  \\
\circ (id \otimes \Delta_{\cal{A}^{\prime}} \otimes id))(y \otimes x 
\otimes f) ,
\end{array} \label{eq:R88primea} \end{eqnarray} which reduces, by 
(\ref{eq:QD47}), to
\begin{eqnarray}
\begin{array}{l}(\widehat{\pi}^{R\;\prime}_{{\cal{L}}(\cal{A})}(x)(Q^{qR}
_{j}))f(y) = \\ (\MC \circ
(ev \otimes \MC) \circ (id \otimes \sigma \otimes id) \circ (id \otimes ev 
\otimes id
\otimes id) \\ 
\circ \{id \otimes id \otimes (\sigma \circ \Delta \circ Q^{qR}_{j}) 
\otimes id\}  \circ (id \otimes
id \otimes \sigma) \circ (id \otimes id \otimes ev \otimes id) \\
\circ (id \otimes id \otimes id \otimes S \otimes id)
\circ (id \otimes  \Delta_{\cal{A}^{\prime}} \otimes \sigma) \circ (id 
\otimes id \otimes
\Delta))(y \otimes x \otimes f) ,
\end{array} \label{eq:R89} \end{eqnarray} However, (\ref{eq:QD37}) 
implies that
\begin{equation} (ev \circ (id \otimes M) \circ (id \otimes \sigma))(x 
\otimes a \otimes b) = (\MC
\circ (ev \otimes ev) \circ (\Delta_{\cal{A}^{\prime}} \otimes id \otimes 
id)(x \otimes a \otimes b) 
\label{eq:R90} \end{equation} for all $a,b \in \cal{A}$ and all $x \in 
\cal{A}^{\prime}$, so
(\ref{eq:R89}) reduces to
\begin{eqnarray}
\begin{array}{l}(\widehat{\pi}^{R\;\prime}_{{\cal{L}}(\cal{A})}(x)(Q^{qR}
_{j}))f(y)   \\ =  (\MC
\circ (ev \otimes ev) \circ (id \otimes \sigma \otimes M) \circ (id \otimes 
id \otimes
\Delta
\otimes id) \\ 
 \;\circ (id \otimes id \otimes  Q^{qR}_{j} \otimes S)  \circ (id \otimes id 
\otimes
\Delta))(y \otimes x \otimes f) \\ 
\end{array}  \end{eqnarray} and hence
\begin{eqnarray}
\begin{array}{l}(\widehat{\pi}^{R\;\prime}_{{\cal{L}}(\cal{A})}(x)(Q^{qR}
_{j}))f(y) =  \\ (\MC \circ
(ev \otimes ev) \circ (id \otimes \sigma \otimes id) \circ (id \otimes id 
\otimes id
\otimes M) \circ (id \otimes id \otimes
\Delta
\otimes id) \\
\circ (id \otimes id \otimes  Q^{qR}_{j} \otimes S) 
   \circ (id \otimes id \otimes
\Delta) \circ (id \otimes id \otimes \MAC))(y \otimes x \otimes f \otimes 
z) 
\end{array} \label{eq:R91a} \end{eqnarray} for all $x,y \in 
\cal{A}^{\prime}$, all $f \in \cal{A}$, 
$Q^{qR}_{j} \in {\cal{L}}({\cal{A}})$, and all $z \in \C$.

Turning to the right-hand side of (\ref{eq:R85}), by (\ref{eq:R26d}),
\[ \sum_{k=1}^{d_{q}}\Gamma^{q}_{kj}(x) (Q^{qR}_{k}(f))(y) = 
\sum_{k=1}^{d_{q}} \langle x,
\pi^{q}_{kj} \rangle \langle y ,Q^{qR}_{k} \rangle ,
\] for all $x,y \in \cal{A}^{\prime}$, and $j = 1 ,2, \dots, d_{q}$, and 
hence
\begin{eqnarray} \begin{array}{l} \sum_{k=1}^{d_{q}}\Gamma^{q}_{kj}(x)( 
Q^{qR}_{k}(f))(y) = \\
\sum_{k=1}^{d_{q}} (\MC
\circ (ev \otimes ev) \circ (id \otimes \sigma \otimes id) \circ (id \otimes 
id \otimes Q^{qR}_{k}
\otimes
\pi^{q}_{kj}) \\
\circ (id \otimes id \otimes id \otimes u))(y \otimes x \otimes f \otimes 
z)  \label{eq:R92}
\end{array}  \end{eqnarray} for all $x,y \in \cal{A}^{\prime}$, all $f \in 
\cal{A}$, $Q^{qR}_{k} \in
{\cal{L}}({\cal{A}})$, and all $z \in
\C$. 

Now (\ref{eq:R85}) implies that the expressions on the right-hand sides of 
(\ref{eq:R91a}) and
(\ref{eq:R92}) may be equated. As the first three factors of each, namely 
$\MC
\circ (ev \otimes ev) \circ (id \otimes \sigma \otimes id)$, are the same, 
the equality holds with
these removed. However, on both sides of this new equality, the factor $y 
\otimes x$ is only acted
on by a succession of identity operators of the form $id \otimes id$. 
Consequently both $y
\otimes x$ and these identity operators can be removed, leaving the result 
that (\ref{eq:R82})  is
equivalent to the condition  
\begin{eqnarray} \begin{array}{l} ((id \otimes M) \circ (\Delta \otimes 
id) \circ (Q^{qR}_{j}
\otimes S) \circ \Delta \circ \MAC)(f \otimes z) = \\
\;\;\sum_{k=1}^{d_{q}} ((Q^{qR}_{k} \otimes
\pi^{q}_{kj}) \circ (id \otimes u))(f \otimes z)
\end{array}  \end{eqnarray} for all  $f \in \cal{A}$, all $j = 1, 2, \ldots, 
d_{q}$, and all $z \in
\C$. This can be rewritten as 
\begin{equation}  ((id \otimes M) \circ (\Delta \otimes id) \circ 
(Q^{qR}_{j}
\otimes S) \circ \Delta)(f) = \sum_{k=1}^{d_{q}} Q^{qR}_{k}(f) \otimes
\pi^{q}_{kj}
 \end{equation}  for all $f \in \cal{A}$ and all $j = 1, 2, \ldots, d_{q}$, 
which is the condition
(\ref{eq:R95}).

Because $M$ is indistinguishable from $M \circ \sigma$ and $S$ is 
indistinguishable from
$S^{-1}$ in the situation being considered here, the above arguments 
would equally well apply with
each of the following 3 substitutions:
\begin{enumerate}
\item replace $M$ by $M \circ \sigma$, but leave $S$ unchanged; 
\item leave $M$ unchanged, but replace $S$ by $S^{-1}$;
\item replace $M$ by $M \circ \sigma$ {\em{and}} replace $S$ by $S^{-
1}$.
\end{enumerate} However, in the general case in which $\cal{A}$ is non-
commutative,  the
possibilities (1) and (2) are {\em{excluded}} because with them the identity 
operator would not be
an irreducible tensor operator belonging to the identity corepresentation. 
Of course, with the
substitution (3), (\ref{eq:R95}) changes into (\ref{eq:R132}), which is the 
defining condition for a
{\em{twisted}} irreducible tensor operator ${\widetilde{Q}}^{qR}_{j}$. 

\subsubsection{Derivations of the right coactions 
$\pi^{R}_{{\cal{L}}(\cal{A})}$ and
$\widetilde{\pi}^{R}_{{\cal{L}}(\cal{A})}$}

First recast (\ref{eq:R84}) as
\begin{equation} \pi^{R\;\prime}_{{\cal{L}}(\cal{A})}(x \otimes Q) = 
\widehat{R}(x) Q
\widehat{R}(x^{-1})
\label{eq:R96} \end{equation} for all $x \in \cal{G}$ and for all linear 
operators $Q$ that act on
functions defined on
$\cal{G}$, where $\pi^{R\;\prime}_{{\cal{L}}(\cal{A})}$ is a mapping 
from ${\cal{A}}^{\prime}
\otimes {\cal{L}}({\cal{A}})$ into ${\cal{L}}({\cal{A}})$. In Hopf algebra 
language, this can be
re-expressed as
\begin{eqnarray} \begin{array}{lll} \pi^{R\;\prime}_{{\cal{L}}(\cal{A})}(x 
\otimes Q) & = &
(\widehat{M}\circ (id \otimes \widehat{M}) \circ (\widehat{R} \otimes id 
\otimes \widehat{R}) \circ
\\
 & &  \circ (id \otimes \sigma) \circ (id \otimes 
S_{\cal{A}^{\prime}}\otimes id) \circ
(\Delta_{\cal{A}^{\prime}} \otimes id))(x \otimes Q) ,
\end{array} \label{eq:R97} \end{eqnarray} where $\widehat{M}$ is the 
operator multiplication
operator defined in (\ref{eq:R98}). It is then easily shown that
$\pi^{R\;\prime}_{{\cal{L}}(\cal{A})}$ is a left action of 
$\cal{A}^{\prime}$ with carrier space
${\cal{L}}({\cal{A}})$. 

This expression for $\pi^{R\;\prime}_{{\cal{L}}(\cal{A})}$ can be re-
expressed in terms of the
structure constants introduced in Section~II with respect to the basis 
$a_{1}, a_{2},
\ldots $ of
$\cal{A}$, and the basis $a^{1}, a^{2}, \ldots $ of its dual 
$\cal{A}^{\prime}$. First define
the  operators ${\cal{P}}^{k}_{j}$  by
\begin{equation} {\cal{P}}^{k}_{j}(a) = \langle a^{k} , a \rangle a_{j}   
\label{eq:R103} \end{equation}  (for all $a \in \cal{A}$ and all $j,k = 1, 2, 
\ldots $), and then
define the matrix elements
$q^{k}_{j}$ of $Q$ by
\begin{equation} q^{k}_{j} = \langle a^{k} , Q(a_{j}) \rangle 
\label{eq:R107} \end{equation}  (for all $j,k = 1, 2, \ldots $). Clearly the 
operators
${\cal{P}}^{k}_{j}$ are members of ${\cal{L}}({\cal{A}})$, and any operator 
$Q$ of
${\cal{L}}({\cal{A}})$ can be expressed as $Q = \sum_{j,k} Q^{j}_{k} 
{\cal{P}}^{k}_{j}$. Then, by
(\ref{eq:R32prime}) and (\ref{eq:R103}), 
\begin{equation} \widehat{R}(a^{m}) = \sum_{j,k =1}^{g} \mu^{jm}_{k} 
{\cal{P}}^{k}_{j}
\label{eq:R104} \end{equation} for all $m = 1, 2\ \ldots, g$. On 
substituting (\ref{eq:R104}) into
(\ref{eq:R97}), and using (\ref{eq:QD37}) and (\ref{eq:QD47}), it follows 
that
\begin{equation} \pi^{R\;\prime}_{{\cal{L}}(\cal{A})}(a^{m} \otimes Q) = 
\sum_{i,j,k,\ell,u,v,w =
1}^{g} (m^{m}_{uv} s^{v}_{w}
\mu^{ju}_{i} \mu^{\ell w}_{k} q^{i}_{\ell})\;{\cal{P}}^{k}_{j} ,
\label{eq:R108} \end{equation}  for all $Q \in {\cal{L}}({\cal{A}})$ and all 
$m = 1, 2\ \ldots, g$.

 Using (\ref{eq:R26b}), the corresponding right coaction 
$\pi^{R}_{{\cal{L}}(\cal{A})}$ of
$\cal{A}$ with the same carrier space
${\cal{L}}({\cal{A}})$ is  given by
\begin{equation} \pi^{R}_{{\cal{L}}(\cal{A})}(Q) =
\sum_{m=1}^{g}\pi^{R\;\prime}_{{\cal{L}}(\cal{A})}(a^{m}
\otimes Q) \otimes a_{m} ,
\label{eq:R109} \end{equation} for all $Q \in {\cal{L}}({\cal{A}})$. Thus 
\[ \pi^{R}_{{\cal{L}}(\cal{A})}(Q) =
\sum_{i,j,k,\ell,m,u,v,w=1}^{g} (m^{m}_{uv} s^{v}_{w}
\mu^{ju}_{i} \mu^{\ell w}_{k} q^{i}_{\ell}) \;( {\cal{P}}^{k}_{j} \otimes 
a_{m})\; ,
 \]  for all $Q \in {\cal{L}}({\cal{A}})$. The final stage is to re-express this 
in a basis free
form. This can be done by writing 
\[ \pi^{R}_{{\cal{L}}(\cal{A})}(Q) = \sum_{[Q]} Q_{[1]} \otimes
Q_{[2]}\;, \]
where  $Q_{[1]} \in {\cal{L}}(\cal{A})$ and $Q_{[1]} \in \cal{A}$ are such 
that
\[  \sum_{[Q]} Q_{[1]}(a) \otimes Q_{[2]} = ((id \otimes M) \circ
(\Delta \otimes id) \circ (Q \otimes S) \circ \Delta)(a) \]
for all $Q \in {\cal{L}}({\cal{A}})$ and all $a \in \cal{A}$.  

The right coaction $\widetilde{\pi}^{R}_{{\cal{L}}(\cal{A})}$ is obtained 
from
$\pi^{R}_{{\cal{L}}(\cal{A})}$ by replacing $M$ by $M \circ \sigma$ and 
replacing $S$ by
$S^{-1}$. The right coactions $\pi^{R}_{{\cal{T}}(\cal{A})}$ and
$\widetilde{\pi}^{R}_{{\cal{T}}(\cal{A})}$ of Sections VI.B.1 and VI.B.2 
have essentially the same
definitions as  $\pi^{R}_{{\cal{L}}(\cal{A})}$ and 
$\widetilde{\pi}^{R}_{{\cal{L}}(\cal{A})}$, except
that their domains are restricted to the appropriate subspaces 
${\cal{T}}(\cal{A})$.

\subsection{Irreducible tensor operators in the group theoretical left 
regular representation
formalism}

\subsubsection{Derivations of the conditions for the $Q^{qL}_{j}$ and
${\widetilde{Q}}^{qL}_{j}$}

The argument for the {\em{left}} regular formalism follows exactly the 
same line as that for the
right regular case given above up to (\ref{eq:R88prime}), the only 
differences being that the
operators $\widehat{R}(x)$ must be replaced by the operators 
$\widehat{L}(x)$, the left action $R$
must be replaced by the left action $L$, and the label $R$ must be 
replaced by $L$ on the irreducible
tensor operators $Q^{qR}_{j}$, on the left action 
$\pi^{R\;\prime}_{{\cal{L}}(\cal{A})}$, and on the
corresponding right coaction $\pi^{R}_{{\cal{L}}(\cal{A})}$. Thus, by 
(\ref{eq:L11}), the analogue
of (\ref{eq:R88primea}) is
\begin{eqnarray}
\begin{array}{l}(\widehat{\pi}^{L\;\prime}_{{\cal{L}}(\cal{A})}(x)(Q^{qL}
_{j}))f(y)\;\; (\; =
\widehat{L}(x) Q^{qL}_{j} \widehat{L}(x^{-1})\;) \;=\\
 (ev \circ (id \otimes  \MAC) \circ (id \otimes id \otimes ev) \circ (id 
\otimes \sigma \otimes id)
\circ (id \otimes id \otimes \sigma) \\\circ (id \otimes id \otimes S 
\otimes id)  \circ  (id
\otimes id \otimes \Delta)  
\circ (id \otimes id \otimes Q^{qL}_{j}) \circ (id \otimes id \otimes 
\MAC) \\ \circ (id \otimes id
\otimes id \otimes ev)  
\circ (id \otimes id \otimes \sigma \otimes id) \circ (id \otimes id 
\otimes id \otimes \sigma)\\
\circ (id \otimes id \otimes id \otimes S \otimes id) \circ
 (id \otimes id \otimes id \otimes \Delta) \\
\circ (id \otimes id \otimes S_{\cal{A}^{\prime}}\otimes id)  
\circ (id \otimes \Delta_{\cal{A}^{\prime}} \otimes id))(y \otimes x 
\otimes f) ,
\end{array}\label{eq:L33a}\end{eqnarray} which reduces, by 
(\ref{eq:QD47}), (\ref{eq:R90}), and
(\ref{eq:QD21}) to
\begin{eqnarray}
\begin{array}{l}(\widehat{\pi}^{L\;\prime}_{{\cal{L}}(\cal{A})}(x)(Q^{qL}
_{j}))f(y) =  \\ (\MC \circ
(ev \otimes ev) \circ (id \otimes \sigma \otimes id) \circ (id \otimes id 
\otimes \sigma)
\circ (id \otimes id \otimes S \otimes id) \\\circ (id \otimes id \otimes 
M \otimes id) 
\circ (id \otimes id \otimes id \otimes \Delta) 
\circ (id \otimes id \otimes S \otimes Q^{qL}_{j})\\ 
   \circ (id \otimes id \otimes
\Delta) \circ (id \otimes id \otimes \MAC))(y \otimes x \otimes f \otimes 
z) 
\end{array} \label{eq:L34} \end{eqnarray} for all $x,y \in 
\cal{A}^{\prime}$, all $f \in \cal{A}$, 
$Q^{qL}_{j} \in {\cal{L}}({\cal{A}})$, and all $z \in \C$.

The right-hand side of the irreducible tensor operator definition that 
corresponds to this is
\begin{eqnarray}
\begin{array}{l}
\sum_{k=1}^{d_{q}}\Gamma^{q}_{kj}(x)( Q^{qL}_{k}(f))(y) = \\
\sum_{k=1}^{d_{q}} (\MC
\circ (ev \otimes ev) \circ (id \otimes \sigma \otimes id) \circ (id \otimes 
id \otimes Q^{qL}_{k}
\otimes
\pi^{q}_{kj}) \\
\circ (id \otimes id \otimes id \otimes u))(y \otimes x \otimes f \otimes 
z)  \label{eq:L35}
\end{array}  \end{eqnarray} for all $x,y \in \cal{A}^{\prime}$, all $f \in 
\cal{A}$, $Q^{qL}_{k} \in
{\cal{L}}({\cal{A}})$, and all $z \in
\C$. 

Equating the right-hand sides (\ref{eq:L34}) and(\ref{eq:L35}), removing 
the common first three
factors ($\MC \circ (ev \otimes ev) \circ (id \otimes \sigma \otimes id)$) 
of each, and removing the
factor $y \otimes x$ and the succession of identity operators of the form 
$id \otimes id$ that act
on $y \otimes x$, it follows that the defining condition becomes  
\begin{eqnarray} \begin{array}{l}  (\sigma \circ (S \otimes id) \circ (M 
\otimes id) \circ (id
\otimes \Delta) \circ (S \otimes  Q^{qL}_{j}) \circ \Delta \circ \MAC)(f 
\otimes z)  \\
 =\sum_{k=1}^{d_{q}} ((Q^{qL}_{k} \otimes
\pi^{q}_{kj}) \circ (id \otimes u))(f \otimes z)
\end{array}  \end{eqnarray} for all  $f \in \cal{A}$, all $j = 1, 2, \ldots, 
d_{q}$, and all $z \in
\C$. This can be rewritten as 
\begin{equation}  (\sigma \circ (S \otimes id) \circ (M \otimes id) \circ 
(id \otimes \Delta) \circ
(S \otimes  Q^{qL}_{j}) \circ \Delta)(f) = \sum_{k=1}^{d_{q}} Q^{qL}_{k}(f) 
\otimes
\pi^{q}_{kj}
 \end{equation}  for all $f \in \cal{A}$ and all $j = 1, 2, \ldots, d_{q}$, 
which is the condition
(\ref{eq:L36}).

Because $\Delta_{{\cal{A}}^{\prime}}$ is indistinguishable from $\sigma 
\circ
\Delta_{{\cal{A}}^{\prime}}$ and $S_{{\cal{A}}^{\prime}}$ is 
indistinguishable from
$S_{{\cal{A}}^{\prime}}^{-1}$ in the situation being considered here, the 
above arguments would
equally well apply to (\ref{eq:L33a}) with each of the following 3 
substitutions:
\begin{enumerate}
\item replace $\Delta_{{\cal{A}}^{\prime}}$ by $\sigma \circ
\Delta_{{\cal{A}}^{\prime}}$, but leave $S_{{\cal{A}}^{\prime}}$ 
unchanged; 
\item leave $\Delta_{{\cal{A}}^{\prime}}$ unchanged, but replace 
$S_{{\cal{A}}^{\prime}}$ by
$S_{{\cal{A}}^{\prime}}^{-1}$;
\item (c) replace $\Delta_{{\cal{A}}^{\prime}}$ by $\sigma \circ
\Delta_{{\cal{A}}^{\prime}}$ {\em{and}} replace $S_{{\cal{A}}^{\prime}}$ 
by
$S_{{\cal{A}}^{\prime}}^{-1}$.
\end{enumerate} In each case the $S$ factor in (\ref{eq:L33a}) should be 
left unchanged because it
comes from the definition (\ref{eq:L1}) of the left regular (right) coaction. 
(Replacing $S$ by
$S^{-1}$ in (\ref{eq:L1}) would give another right coaction, but the 
original one is merely the
double contragredient of this). In the general case in which
$\cal{A}$ is non-commutative, the possibilities (1) and (2) are again 
{\em{excluded}} because with
them the identity operator would not be an irreducible tensor operator 
belonging to the identity
corepresentation. However, with the substitution (3), the analogue of 
(\ref{eq:L36}) is
(\ref{eq:L52}), which is the defining condition for a {\em{twisted}} 
irreducible tensor operator
${\widetilde{Q}}^{qL}_{j}$. 

\subsubsection{Derivations of the right coactions 
$\pi^{L}_{{\cal{L}}(\cal{A})}$ and
$\widetilde{\pi}^{L}_{{\cal{L}}(\cal{A})}$}

The {\em{left}} regular analogues of (\ref{eq:R96}) and (\ref{eq:R97}) are
\[ \pi^{L\;\prime}_{{\cal{L}}(\cal{A})}(x \otimes Q) = \widehat{L}(x) Q 
\widehat{L}(x^{-1}) \] and
\begin{eqnarray} \begin{array}{lll} \pi^{L\;\prime}_{{\cal{L}}(\cal{A})}(x 
\otimes Q) & = &
(\widehat{M}\circ (id \otimes \widehat{M}) \circ (\widehat{L} \otimes id 
\otimes \widehat{L}) \circ
\\
 & &  \circ (id \otimes \sigma) \circ (id \otimes 
S_{\cal{A}^{\prime}}\otimes id) \circ
(\Delta_{\cal{A}^{\prime}} \otimes id))(x \otimes Q) ,
\end{array} \label{eq:extra30} \end{eqnarray} where 
$\pi^{L\;\prime}_{{\cal{L}}(\cal{A})}$ is a left
action of $\cal{A}^{\prime}$ with carrier space ${\cal{L}}({\cal{A}})$. 
However, by (\ref{eq:L12})
and (\ref{eq:R103}), and with the basis of
$\cal{A}^{\prime}$ defined above,
\begin{equation} \widehat{L}(a^{m}) = \sum_{j,k,\ell =1}^{g} \mu^{\ell 
j}_{k} s^{m}_{\ell}
{\cal{P}}^{k}_{j}
\label{eq:L37} \end{equation} for all $m = 1, 2\ \ldots, g$. On 
substituting (\ref{eq:L37}) into
(\ref{eq:extra30}), and using (\ref{eq:QD37}) and (\ref{eq:QD47}), it 
follows that
\[ \pi^{L\;\prime}_{{\cal{L}}(\cal{A})}(a^{m} \otimes Q) = 
\sum_{i,j,k,\ell,n,u,v,w = 1}^{g}
(m^{v}_{wu} s^{m}_{v} s^{w}_{n}
\mu^{uj}_{i} \mu^{n\ell}_{k} q^{i}_{\ell})\;{\cal{P}}^{k}_{j} ,\]   for all $Q 
\in
{\cal{L}}({\cal{A}})$ and all $m = 1, 2\ \ldots, g$. 

Then, using (\ref{eq:R26b}), the
corresponding right coaction
$\pi^{L}_{{\cal{L}}(\cal{A})}$ of $\cal{A}$ with the same carrier space
${\cal{L}}({\cal{A}})$ is  given by
\[ \pi^{L}_{{\cal{L}}(\cal{A})}(Q) = 
\sum_{m=1}^{g}\pi^{L\;\prime}_{{\cal{L}}(\cal{A})}(a^{m}
\otimes Q)
\otimes a_{m} ,\] for all $Q \in {\cal{L}}({\cal{A}})$. Thus 
\[ \pi^{L}_{{\cal{L}}(\cal{A})}(Q) =
\sum_{i,j,k,\ell,m,n,u,v,w = 1}^{g} (m^{v}_{wu} s^{m}_{v} s^{w}_{n}
\mu^{uj}_{i} \mu^{n\ell}_{k} q^{i}_{\ell})\;({\cal{P}}^{k}_{j} \otimes 
a_{m})\; , \]   for all $Q \in
{\cal{L}}({\cal{A}})$. The final stage is to re-express this in a basis free
form, which can be done by writing 
\[ \pi^{L}_{{\cal{L}}(\cal{A})}(Q) = \sum_{[Q]} Q_{[1]} \otimes
Q_{[2]}\;, \]
where  $Q_{[1]} \in {\cal{L}}(\cal{A})$ and $Q_{[1]} \in \cal{A}$ are such 
that
\[  \sum_{[Q]} Q_{[1]}(a) \otimes Q_{[2]} = (\sigma \circ (S
\otimes id)
\circ (M \otimes id) \circ (id \otimes \Delta) \circ (S \otimes Q) \circ 
\Delta)(a) \]
for all $Q \in {\cal{L}}({\cal{A}})$ and all $a \in \cal{A}$.

For the corresponding twisted coaction 
$\widetilde{\pi}^{L}_{{\cal{L}}(\cal{A})}$ the argument is
similar. With the substitutions $\Delta_{\cal{A}^{\prime}} \rightarrow
\sigma \circ \Delta_{\cal{A}^{\prime}}$ and 
$S_{\cal{A}^{\prime}}\rightarrow
S_{\cal{A}^{\prime}}^{-1}$, (\ref{eq:extra30}) gives
\begin{eqnarray} \begin{array}{lll} 
\widetilde{\pi}^{L\;\prime}_{{\cal{L}}(\cal{A})}(x \otimes Q) &
= & (\widehat{M}\circ (id \otimes \widehat{M}) \circ (\widehat{L} 
\otimes id \otimes \widehat{L})
\circ
\\
 & &  \circ (id \otimes \sigma) \circ (id \otimes S_{\cal{A}^{\prime}}^{-
1}\otimes id) \circ
((\sigma \otimes \Delta_{\cal{A}^{\prime}}) \otimes id))(x \otimes Q) ,
\end{array} \label{eq:L54} \end{eqnarray} where 
$\widetilde{\pi}^{L\;\prime}_{{\cal{L}}(\cal{A})}$ is
another left action of $\cal{A}^{\prime}$ with carrier space 
${\cal{L}}({\cal{A}})$. By
(\ref{eq:L37}) and (\ref{eq:QD22}), this gives
\[ \widetilde{\pi}^{L\;\prime}_{{\cal{L}}(\cal{A})}(a^{m} \otimes Q) = 
\sum_{i,j,k,\ell,n,u,v =
1}^{g} (m^{m}_{nv} s^{v}_{u} \mu^{uj}_{i} \mu^{n\ell}_{k} 
q^{i}_{\ell})\;{\cal{P}}^{k}_{j} ,\]   for
all $Q
\in {\cal{L}}({\cal{A}})$ and all $m = 1, 2\ \ldots, g$. Using 
(\ref{eq:R26b}), the
corresponding right coaction
$\widetilde{\pi}^{L}_{{\cal{L}}(\cal{A})}$ of $\cal{A}$ with the same 
carrier space
${\cal{L}}({\cal{A}})$ is  given by
\[ \widetilde{\pi}^{L}_{{\cal{L}}(\cal{A})}(Q) =
\sum_{m=1}^{g}\widetilde{\pi}^{L\;\prime}_{{\cal{L}}(\cal{A})}(a^{m} 
\otimes Q)
\otimes a_{m} ,\] for all $Q \in {\cal{L}}({\cal{A}})$. Thus 
\[ \widetilde{\pi}^{L}_{{\cal{L}}(\cal{A})}(Q) =
\sum_{i,j,k,\ell,m,n,u,v = 1}^{g} (m^{m}_{nv} s^{v}_{u} \mu^{uj}_{i} 
\mu^{n\ell}_{k}
q^{i}_{\ell})\;({\cal{P}}^{k}_{j} \otimes a_{m})\; , \]   for all $Q \in 
{\cal{L}}({\cal{A}})$. 
This can be re-expressed in a basis free
form by writing 
\[ \widetilde{\pi}^{L}_{{\cal{L}}(\cal{A})}(Q) = \sum_{[Q]} Q_{[1]} \otimes
Q_{[2]}\;, \]
where  $Q_{[1]} \in {\cal{L}}(\cal{A})$ and $Q_{[1]} \in \cal{A}$ are such 
that
\[  \sum_{[Q]} Q_{[1]}(a) \otimes Q_{[2]} = ((id
\otimes M)
\circ (\sigma \otimes S) \circ (id \otimes \sigma) \circ (id
\otimes \Delta) \circ (id \otimes Q) \circ \Delta)(a) \]
for all $Q \in {\cal{L}}({\cal{A}})$ and all $a \in \cal{A}$.

The right coactions $\pi^{L}_{{\cal{T}}(\cal{A})}$ and
$\widetilde{\pi}^{L}_{{\cal{T}}(\cal{A})}$ of Sections VI.C.1 and VI.C.2 
have essentially the same
definitions as  $\pi^{L}_{{\cal{L}}(\cal{A})}$ and 
$\widetilde{\pi}^{L}_{{\cal{L}}(\cal{A})}$, except
that their domains are restricted to the appropriate subspaces 
${\cal{T}}(\cal{A})$.

\section{Projection operators}

The {\em{right}} regular formalism will be considered first. If 
${\bf{\Gamma}}^{p}$ is a unitary
irreducible representation of dimension
$d_{p}$ of a finite group $\cal{G}$ of order $g$, the projection operators in 
the right regular
formalism are defined by
\[ {\cal{P}}^{pR}_{mn} = (d_{p}/g) \sum_{x \in {\cal{G}}} 
\Gamma^{p}(x)^{*}_{mn} \widehat{R}(x) \]
for all $m,n = 1, 2\ \ldots, d_{p}$.  This can be re-written as
\[ {\cal{P}}^{pR}_{mn} = (d_{p}/g) \sum_{x \in {\cal{G}}} \Gamma^{p}(x^{-
1})^{*}_{nm} \widehat{R}(x)
,\] and hence, by (\ref{eq:R26d}), (\ref{eq:R29}), (\ref{eq:R30}), and 
(\ref{eq:R87,R88}),
\[ \langle y ,{\cal{P}}^{pR}_{mn}f \rangle = (d_{p}/g) \sum_{x \in 
{\cal{G}}} \langle
S_{\cal{A}^{\prime}}(x) ,
\pi^{p}_{nm} \rangle \langle y, R(x \otimes f) \rangle \] for all $m,n = 1, 
2\ \ldots, d_{p}$. Here
the $\pi^{p}_{nm}$ are the matrix coefficients of the corepresentation 
$\pi^{p}$ of
$\cal{A}$ that is dual to ${\bf{\Gamma}}^{p}$. Then, by 
(\ref{eq:R32prime}), for all$x,y \in
\cal{A}^{\prime}$ and all $f \in \cal{A}$,
\[ \langle y ,{\cal{P}}^{pR}_{mn}f \rangle = (d_{p}/g) \sum_{x \in 
{\cal{G}}} \langle ((id \otimes
id \otimes S_{\cal{A}^{\prime}}) \circ (id \otimes 
\Delta_{\cal{A}^{\prime}})) (y
\otimes x) , (\pi^{R}_{\cal{A}}(f) \otimes \pi^{p}_{nm}) \rangle ,\] and so, 
by
(\ref{eq:D.sect1.1.45ii}) and (\ref{eq:R47}), 
\begin{equation} \langle y ,{\cal{P}}^{pR}_{mn}f \rangle = (d_{p}/g) 
\sum_{x \in {\cal{G}}} \langle
(y \otimes x) ,((id \otimes M) \circ (id \otimes id \otimes S)) 
(\pi^{R}_{\cal{A}}(f) \otimes
\pi^{p}_{nm}) \rangle .\label{eq:R137} \end{equation} But the Haar 
functional is such that
\[ h(a) = (1/g) \sum_{x \in {\cal{G}}} \langle x , a \rangle \] for all $a \in 
\cal{A}$, so
(\ref{eq:R137}), (\ref{eq:extra15}), and (\ref{eq:D.sect1.1.45ii}) imply that 
\begin{equation} {\cal{P}}^{pR}_{mn}f  = d_{p} \sum_{[f]} f^{R}_{[1]} 
h(M(f^{R}_{[2]} \otimes
\pi^{p*}_{mn}))  ,\label{eq:R139} \end{equation} As multiplication is 
commutative in this special
case, this could equally well be written as
\begin{equation} {\cal{P}}^{pR}_{mn}f  = d_{p} \sum_{[f]} f^{R}_{[1]} 
h(M(\pi^{p*}_{mn} \otimes
f^{R}_{[2]}))  \label{eq:R140prime} \end{equation} for all $m,n = 1, 2\ 
\ldots, d_{p}$. In the
general case the two formulae (\ref{eq:R139}) and (\ref{eq:R140prime}) 
are different, but the
arguments given in Subsection~IV.B show that (\ref{eq:R140prime}) (i.e.
(\ref{eq:R143a})) is actually the correct choice.

The argument in the {\em{left}} regular formalism follows exactly the 
same line, and can be obtained
by merely replacing the label $R$ by $L$ at each stage.

\end{document}